\documentclass[review]{elsarticle}

\usepackage[a4paper, total={6in,9in}]{geometry}

\usepackage{lipsum}
\usepackage{xcolor, hyperref}
\hypersetup{
    colorlinks = true,
    linkcolor = red,
    linkbordercolor = {white}
}

\usepackage{lipsum}
\usepackage{amsfonts}
\usepackage{graphicx}
\usepackage{epstopdf}
\usepackage{algorithmic}
\usepackage{amsmath}
\usepackage{amsopn}
\usepackage{multirow}
\usepackage{xcolor}
\usepackage[title]{appendix}

\newcommand{\refereeSecond}[1]{\textcolor{black}{#1}}
\newcommand{\refereeThird}[1]{\textcolor{black}{#1}}

\DeclareMathOperator*{\argmin}{arg\,min}
\newtheorem{lemma}{Lemma}

\begin{document}

\begin{frontmatter}

    \title{An unfitted RBF-FD method in a least-squares setting for elliptic PDEs on
        complex geometries}

    \author[igor]{Igor Tominec}
    \ead{igor.tominec@it.uu.se}
    \address[igor]{Uppsala University, Department of Information Technology, Division of Scientific Computing}

    \author[eva]{Eva Breznik}
    \ead{eva.breznik@it.uu.se}
    \address[eva]{Uppsala University, Department of Information Technology, Centre for Image Analysis}



    \begin{abstract}
        Radial basis function generated finite difference (RBF-FD) methods  for PDEs require a set of interpolation points which
        conform to the computational domain $\Omega$. One of the requirements leading to approximation robustness is to place
        the interpolation points with a locally uniform distance around the boundary of $\Omega$. However generating
        interpolation points with such properties is a cumbersome problem.
        Instead, the interpolation points can be extended
        over the boundary and as such completely decoupled from the shape of $\Omega$.  In this paper we present a modification to
        the least-squares RBF-FD method which allows the interpolation points to be placed in a box that encapsulates $\Omega$.
        This way, the node placement over a complex domain in 2D and 3D is greatly simplified.
        Numerical experiments on solving an elliptic model PDE over complex 2D geometries show that our approach is robust.
        Furthermore it performs better in terms of the approximation error and the runtime vs. error compared with the classic RBF-FD methods.
        It is also possible to use our approach in 3D, which we indicate
        by providing convergence results of a solution over a thoracic diaphragm.
    \end{abstract}


\end{frontmatter}


\section{Introduction}
\begin{figure}[h!]
    \centering
\includegraphics[width=0.35\linewidth]{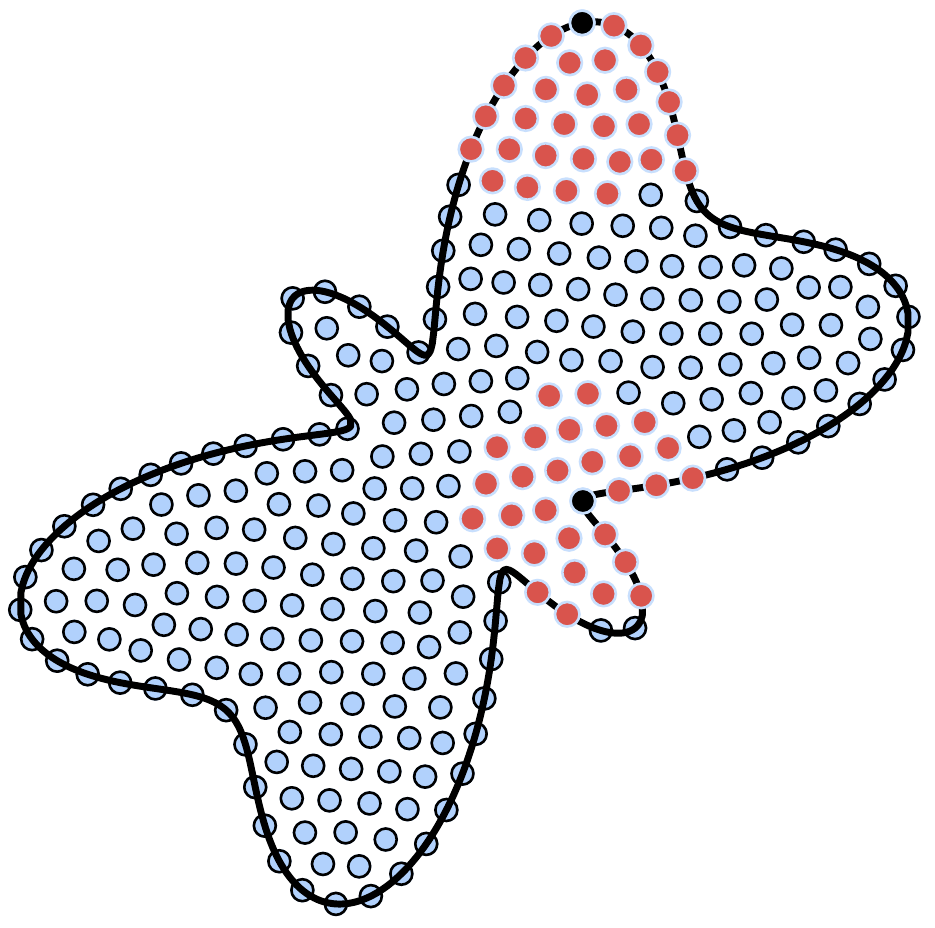}
\includegraphics[width=0.43\linewidth]{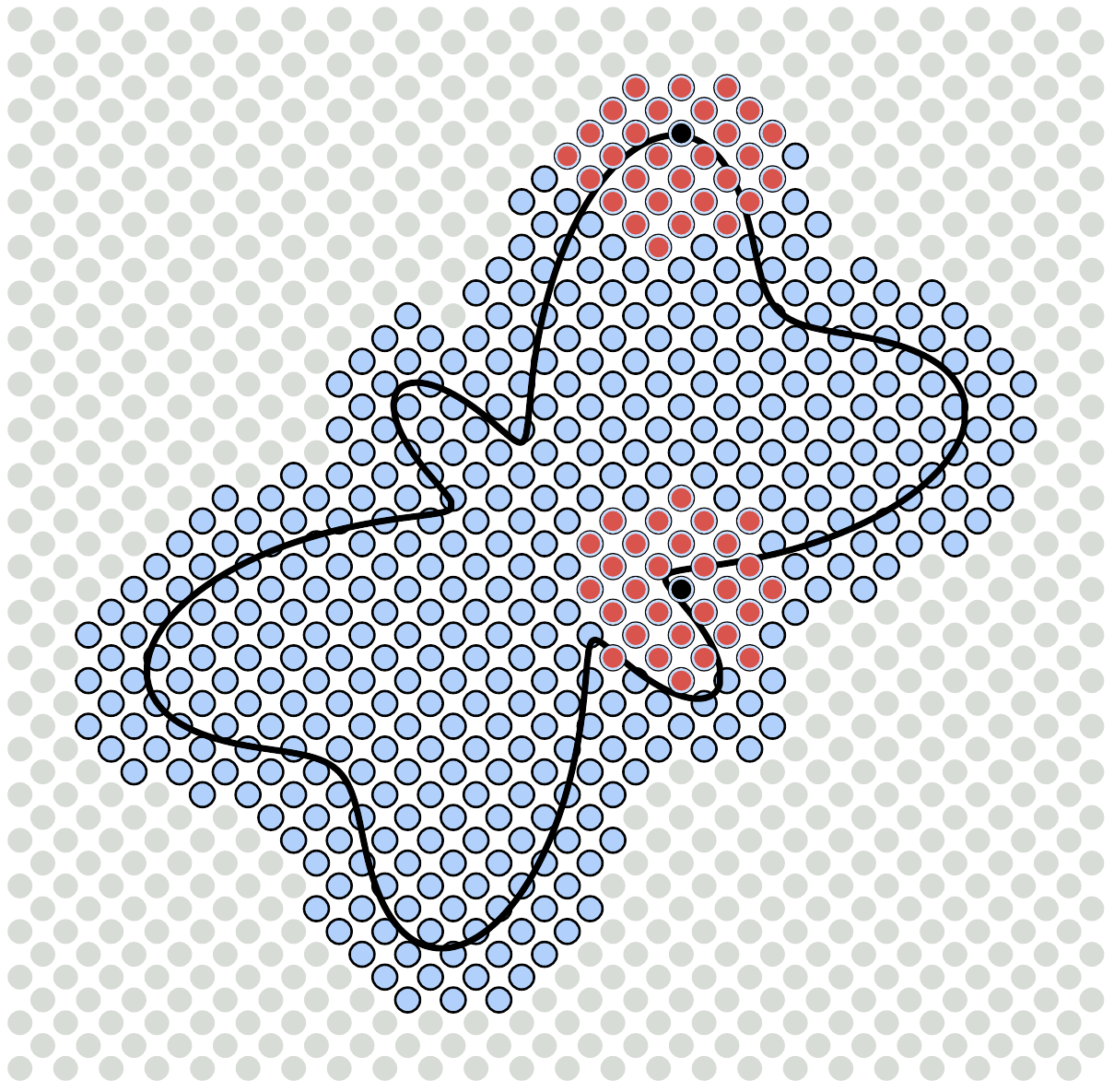}
\caption{Node distributions over a butterfly. Left: a classical distribution of interpolation points over $\Omega$ for RBF-FD in the least-squares 
setting with two skewed stencils on the boundary. Here the black point is the stencil center and the red points are the members of the stencil. 
Right: a node distribution
over $\Omega$ for the unfitted RBF-FD method in the least-squares setting with a less skewed stencil on the boundary.}
\label{fig:intro:butterfly}
\end{figure}

Most localized radial basis function (RBF) methods for computing solutions to partial differential equations (PDEs)
 (for example \cite{Tolstykh, shankar_overlapped, milovanovic18})
require a set of interpolation points that
conforms to a computational domain $\Omega$.
Normally, a localized RBF method uses a collection of local interpolation problems over the subsets of interpolation points \refereeSecond{(e.g. stencils or patches)} to generate 
the compactly supported cardinal functions which are then employed to solve a PDE.
It is well known that locally non-uniform node distances between interpolation points increase the conditioning of the local
interpolation problem and cause unwanted growth of the cardinal functions \cite{BF10,shankar_overlapped}.
While it is trivial to place uniform points in the very interior of $\Omega$, it is on the other hand
challenging to place them in the vicinity of a boundary of $\Omega$ and at the same time maintain their uniformity. 
This is especially difficult 
in three dimensions. 
\refereeSecond{One possibility to avoid these difficulties is to improve the node generation techniques \cite{slakandkosec,shankarnodes,BFF_nodes,zamolonobile}. 
Another possibility to circumvent the aforementioned difficulties is to build a PDE discretization over a set of interpolation points decoupled from the shape of $\Omega$ (see Figure \ref{fig:intro:butterfly}): this is what we focus on in this work.}

\refereeSecond{The first steps towards solving PDEs using a computational grid decoupled from the computational domain were made by Peskin in 
1972 \cite{Peskin2, Peskin1}, 
who introduced the \emph{Immersed boundary method} (IB) in order to 
simplify the simulations of fluid-structure interactions (FSI). 
The method uses a fixed, structured Eulerian grid to move the flow, 
and a Lagrangian grid to move a curvilinear immersed boundary. 
The two grids are coupled through a special forcing term that is added to the continuous PDE. 
A difficulty in IB is that the 
FSI solution quantities (pressure and stress) become discontinuous at the interface of the fluid and the solid, 
leading to low-order accuracy. A remedy is the \emph{Immersed interface method} \cite{Leveque}, 
which directly imposes the jump condition 
over the interface, and by that enables high-order accuracy.
The Cartesian cut cell methods \cite{cutcell1, cutcell2} also belong to the related work, where a finite difference method (FD) is employed on a 
background Cartesian grid, cut by an arbitrarily shaped immersed boundary. The interior cells are treated as in classical FD, while the 
cut cells get a different treatment, depending on the application.
}

The so called \emph{unfitted methods} \cite{Barett, Hansbo, Burman}
which are a part of the finite element
methods address the decoupling of $\Omega$ and a mesh:
the boundary of $\Omega$ is enclosed in a box with a background mesh,
where only the elements which have a non-empty intersection with $\Omega$ are taken to be active.
Those methods tend to suffer from ill-conditioning in the presence of small and irregular cuts close to
the boundary of $\Omega$. This has for example been addressed in \cite{Sticko1, Sticko2, elfverson2018cutiga}. 
An additional challenge is the enforcement of Dirichlet boundary conditions,
which has been addressed by introducing a penalty term over the boundary elements \cite{Nitsche, Babuska}.

Another related approach is the placement of \emph{ghost points} (also \emph{fictitious points}) in finite difference methods, where additional
points (unknowns) are placed outside of $\Omega$
in order to enforce Neumann type boundary conditions in a more accurate way.
This is a widely used concept, an example can be found in \cite{bf_ghostpoints}, where for every added point,
an additional equation is generated in order to maintain a square linear system of equations.
In \cite{rbf_ghostpoints2, rbf_ghostpoints} the authors introduce ghost points for a global radial basis function (RBF) collocation
method in order to decouple interpolation points from $\Omega$.
Their computational study shows that this is a feasible approach and that the error under node refinement
tends to be smaller compared with the fitted method.
However the study is limited to using basis functions with a global support, and the study 
does not provide an insight into how many ghost points to use and how that affects the stability properties.

A partition of unity based RBF method in a least-squares setting (RBF-PUM-LS) \cite{rbf-pum-ls} 
enables a decoupling of the interpolation points 
from $\Omega$ by placing a set of overlapping patches over $\Omega$, 
where every patch contains interpolation points independent of the shape of $\Omega$. The authors provided 
numerical evidence that RBF-PUM-LS is an accurate and robust method to solve an elliptic model problem, 
but have not studied the effects of the patches that extend outside $\Omega$.

A recently introduced
\emph{RBF-FD method in a least-squares setting} (RBF-FD-LS) \cite{tominec2020squares} was proven to be significantly
more robust compared to the same method in the collocation setting (RBF-FD-C), especially 
in the presence of Neumann-type boundary conditions \cite{tominec2020squares}. However the interpolation points 
are required to conform to $\Omega$. Another study leading to a least-squares RBF-FD was introduced in \cite{Davydov19}.
\begin{figure}[h!]
    \centering
\includegraphics[width=0.35\linewidth]{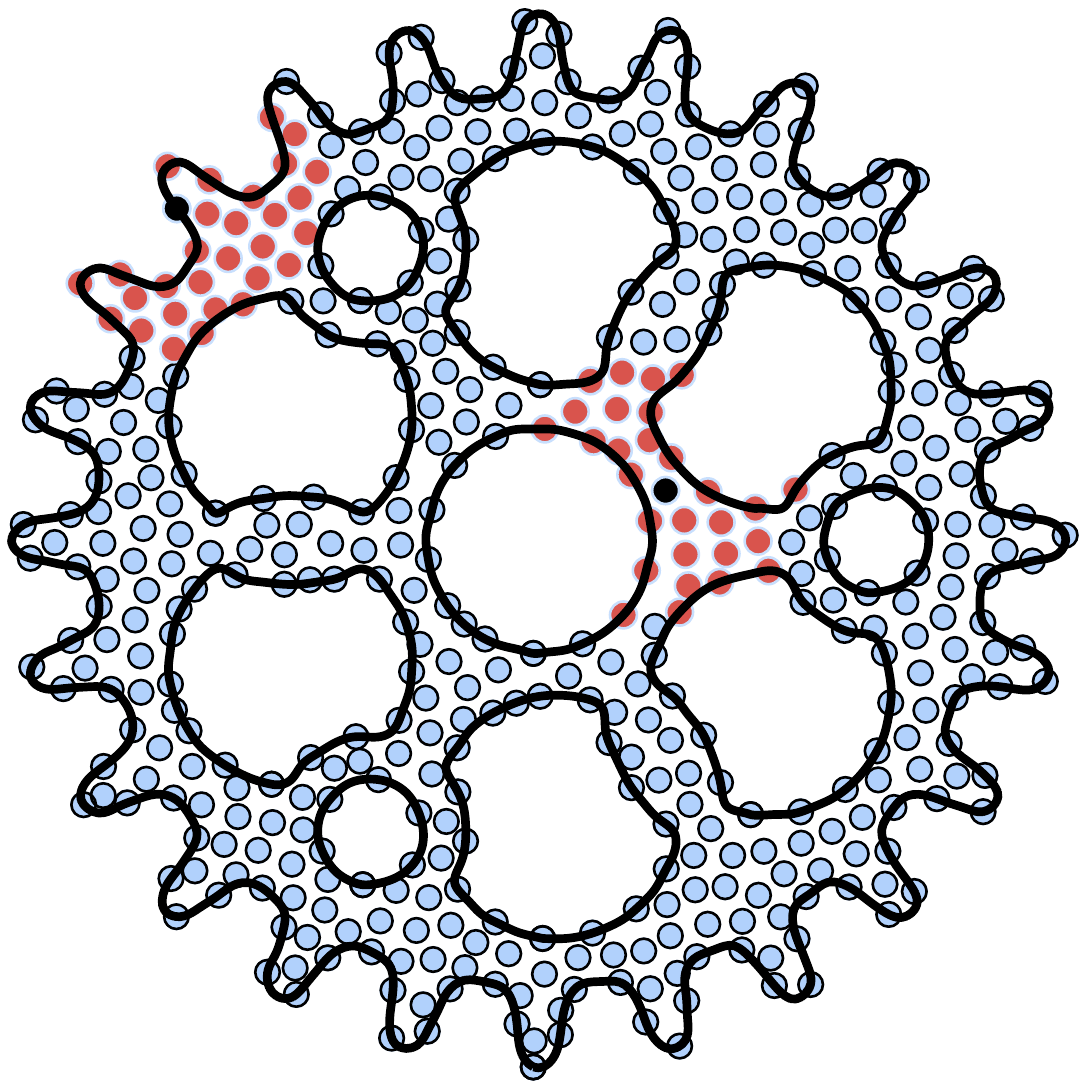}
\includegraphics[width=0.43\linewidth]{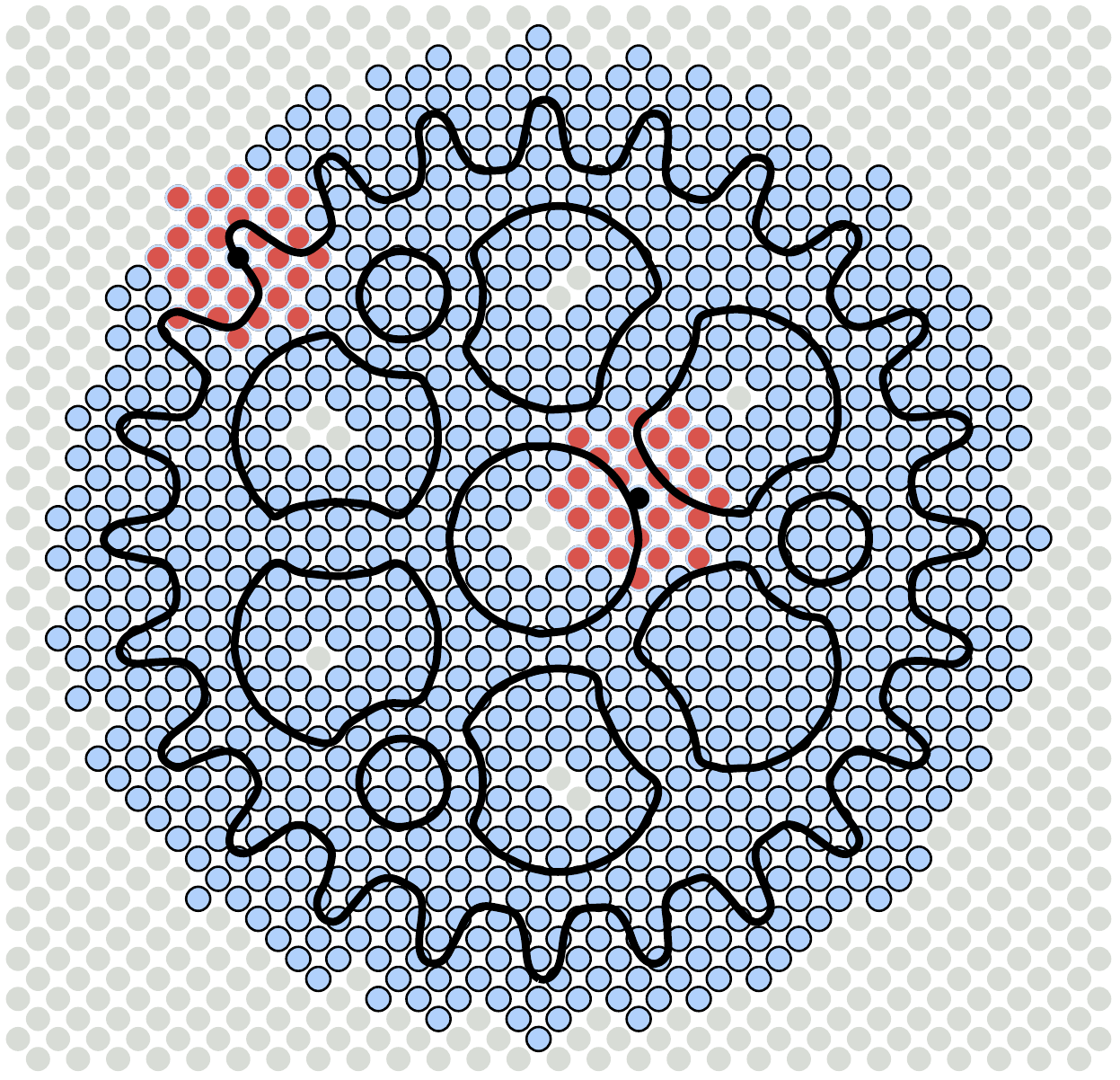}
\caption{Node distributions over a drilled 24-tooth sprocket. Left: a classical distribution of interpolation points over $\Omega$ for RBF-FD in the 
least-squares setting with two skewed stencils on the boundary. Here the black point is the stencil center and the red points are the members of the stencil. Right: a node distribution
over $\Omega$ for the unfitted RBF-FD method in the least-squares setting with a less skewed stencil on the boundary.}
\label{fig:intro:sprocket}
\end{figure}

In this paper we use the RBF-FD-LS method \cite{tominec2020squares} and for that method introduce an approach to decouple interpolation nodes from $\Omega$. 
A computational domain $\Omega$ is enclosed in a box that contains a set of
regularly spaced interpolation nodes with reasonably good interpolation properties (see Figure \ref{fig:intro:butterfly}
and Figure \ref{fig:intro:sprocket}).
The solution unknowns are determined by solving a system where every equation is an evaluation of a PDE
in a point $y \in \bar \Omega$. Here every $y$ picks the closest interpolation point which is used as a center of an approximation stencil 
built over a set of neighboring interpolation points in a box. The shape of such a stencil is given in Figure \ref{fig:intro:butterfly}
and Figure \ref{fig:intro:sprocket}.
An additional strength of the \refereeSecond{present} approach is that the approximation stencils around the boundary become less skewed
which is found to reduce the magnitude of the error around the boundary, especially when the stencil size is large.


The paper is organized as follows.
In Section \ref{section:modelproblem} we state \refereeSecond{the model} problem.
In Section \ref{section:method} we provide a description of
the unfitted RBF-FD-LS method, together with the formulas for computing local differentiation weights, the global differentiation matrices and the discretization of \refereeSecond{the model} problem.
In Section \ref{sec:method:linearindep} we study linear independence of the cardinal functions
as the interpolation points move away from the boundary of
$\Omega$ and develop a heuristic criterion to keep the linear independence of cardinal functions unchanged. This is a necessary condition 
for the well-posedness of the discrete PDE problem. 
In Section \ref{section:erroranalysis} we study the relation between the discrete solution and the analytic solution by
deriving a discrete error estimate without an a-priori bound on the stability norm.
The 2D experiments are presented in Sections \ref{section:experiments_butterfly} and \ref{section:experiments_sprocket}.
In the former section we consider a butterfly domain and numerically investigate the behaviours of the error against true solution, stability norm and the condition numbers
under node and polynomial degree refinements. The results are compared against
RBF-FD-LS and RBF-FD-C.
In the latter we use a drilled 24-tooth sprocket as a computational domain
and study the effects of the unfitted discretization on the spatial distribution of the error for a fixed internodal distance and
several polynomial degrees.
In Section \ref{section:experiments_diaphragm} we numerically study the convergence under node refinement in 3D, where we use 
a thoracic diaphragm geometry extracted from medical images.
Lastly, Section \ref{section:finalremarks} concludes the paper and offers directions for further work.

\section{The model problem}
\label{section:modelproblem}
We choose to evaluate our method by solving the Poisson equation on an open and bounded domain $\Omega \subset \mathbb{R}^d$ ($d=2,3$)
with mixed boundary conditions.
\begin{eqnarray}
	\label{eq:model:Poisson}
	  \Delta u(y) &=& f_2(y),\quad y \in \Omega, \nonumber\\
	  u(y) &=& f_0(y),\quad  y \in \partial\Omega_0, \nonumber\\
	  \nabla u(y) \cdot n &=& f_1(y),\quad y \in \partial\Omega_1,
  \end{eqnarray}
where $\partial\Omega_0$ and $\partial\Omega_1$ are two disjoint parts of a smooth boundary $\partial\Omega$. The solution $u$ is throughout the paper 
assumed to be smooth.
In the theoretical parts of the paper we prefer to work with the following formulation of the same problem:
\begin{equation}
    \label{section:modelproblem_system}
    Du(y) = F(y),
\end{equation}
where:
\begin{eqnarray}
    \label{eq:model:Poisson_system}
    Du(y) =
    \left\{
    \begin{array}{ll}
        \Delta u(y), & y \in \Omega, \\
        u(y), & y \in \partial\Omega_0, \\
        \nabla u(y) \cdot n(y), & y \in \partial\Omega_1, \\
    \end{array}
    \right. 
     \mbox{ and }
    F(y) =
    \left\{
    \begin{array}{ll}
        f_2(y), & y \in \Omega, \\
        f_0(y), & y \in \partial\Omega_0, \\
        f_1(y), & y \in \partial\Omega_1 \\
    \end{array}.
    \right.
\end{eqnarray}
The numerical solution is going to be sought using:
\begin{equation}
    \label{eq:model:ansatz}
u_h(y) = \sum_{i=1}^N u_h(x_i) \Psi_i(y),
\end{equation}
where $\Psi_i(y)$, $i=1,..,N$ are the RBF-FD cardinal functions and $u_h(x_i)$ are the nodal values of the solution. 
\section{The unfitted RBF-FD method}
\label{section:method}
In this section we discuss the choice of point sets that discretize the computational domain, the generation of the cardinal functions 
\eqref{eq:model:ansatz} using the RBF-FD method,
the assembly of the evaluation and differentiation matrices and the discretization of the model problem \eqref{eq:model:Poisson}.

\subsection{The point sets}
Two sets of computational points are distributed over $\Omega$:
\begin{itemize}
	\item  \emph{The interpolation point set} $X=\{x_i\}_{i=1}^N$ for generating the cardinal functions.
	\item \emph{The evaluation point set} $Y=\{y_i\}_{j=1}^M$ for sampling the PDE \eqref{eq:model:Poisson}.
\end{itemize}
We noted in \cite{tominec2020squares} that the $X$ points are supposed to be distributed such that the internodal distance
is as uniform as possible, since the Lebesgue constants associated with the cardinal functions then stay fairly small. On the other
hand the evaluation point set \refereeThird{$Y$} does not influence the magnitude of the Lebesgue constants but is instead important for
the implicit integration that occurs when solving a discretized system of equations in the least-squares sense \cite{tominec2020squares}.
Thus the constraints for placing the \refereeThird{Y} points are far more forgiving: \refereeThird{as long as $\Omega$ is overall well covered with $Y$ points, some 
of them can lie very close to each other, also in the vicinity of the boundary.}

We choose the interpolation point set such that it does not conform to the computational domain $\Omega$
(see Figure \ref{fig:intro:butterfly} and Figure \ref{fig:intro:sprocket}). \refereeSecond{Throughout the paper we
take $X$ to be a tilted Cartesian grid with spacing $h$ for 2D cases. For the 3D case we use a point set obtained using the algorithm from \cite{Kiera19}.} The motivation for that is a simplified point generation,
and the benefit of a polynomial unisolvency on those points. The latter is important for the well-posedness of an interpolation
problem over a stencil.

\begin{figure}
\centering
	\includegraphics[width=0.5\linewidth]{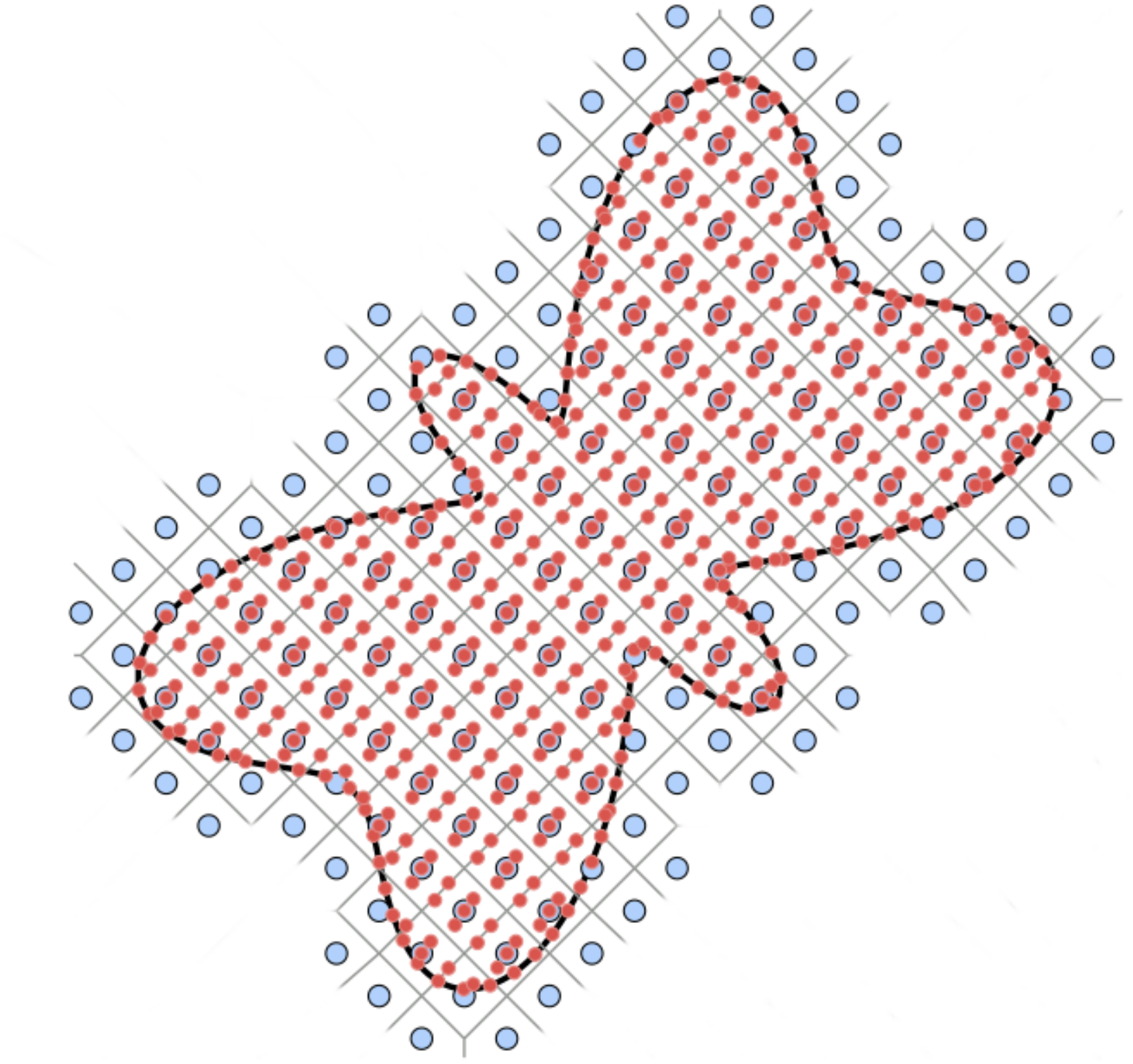}
	\includegraphics[width=0.4\linewidth]{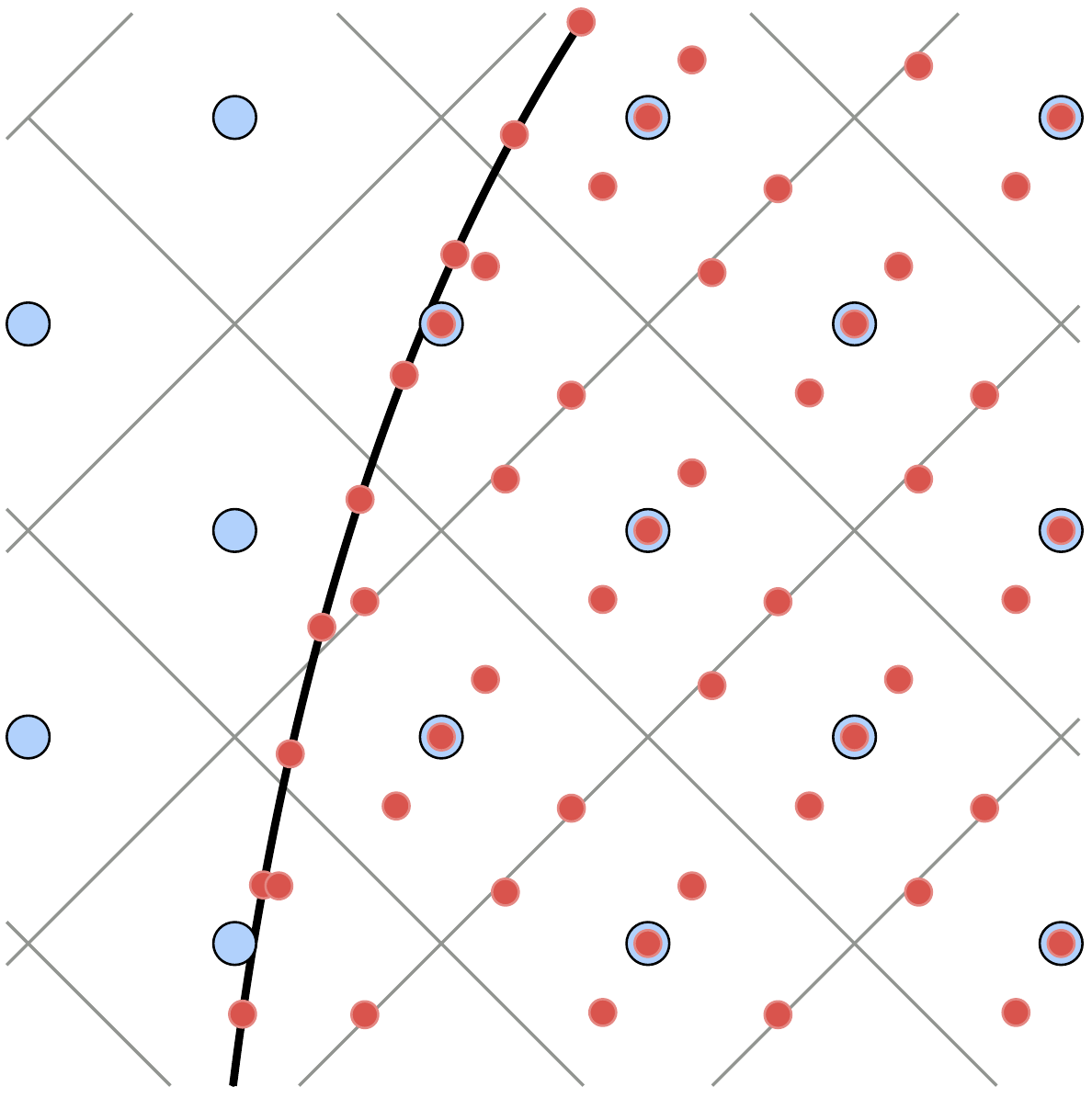}
	\caption{The image on the left displays a butterfly domain with evaluation points (smaller red markers). \refereeSecond{Evaluation points are formed such that local templates of $q=5$ points are 
	placed in Voronoi cells (grey lines) 
	centered around every interpolation
	point (larger blue markers). The image on the right is a closer view over the boundary where $q=5$ points are placed in every Voronoi cell.}}
	\label{fig:method:evalpts_voronoi}
\end{figure}

The evaluation point set $Y$ conforms to $\Omega$.
\refereeSecond{We can subdivide $Y$ into the points on the boundary and the points in the interior of the domain. 
We choose the interior points such that there are $q$ points placed in every
Voronoi cell centered around each $x\in X$. Here $q$ is the oversampling parameter.} 
The boundary points are then placed on
$\partial\Omega$ with a uniform distance that corresponds to the distance between the interior points.
For a visual representation see Figure \ref{fig:method:evalpts_voronoi}.
With such a relation between the $X$- and the $Y$-points it follows that the cardinality of those sets very closely matches the relation:
$M \approx qN$, where $M$ is the number of $Y$-points and $N$ is the number of $X$-points. 
\refereeThird{In \cite{tominec2020squares} we studied the effect of $q$ on the stability properties of the oversampled discretization and to the accuracy of the final approximation. 
We found out that in 2D, the stability and accuracy improved until $q=3$, but the improvement after that point was not significant anymore. In our experience $q>3$ 
was always a good choice of the oversampling parameter for solving the steady-state PDEs in $2$D. 
In 3D, the oversampling can be determined by the relation $q_{\text{3D}} \approx q_{\text{2D}}\sqrt{q_{\text{2D}}}$, 
based on a reasonable assumption that sampling with $q_{\text{1D}}$ points in $1D$ is equivalent to sampling with 
$q_{\text{2D}} = q_{\text{1D}}^2$ points in $2D$, which is in $3D$ 
equivalent to sampling with $q_{\text{3D}} = q_{\text{1D}}^3 = q_{\text{1D}}\, q_{\text{1D}}^2 = \sqrt{q_{\text{2D}}}\, q_{\text{2D}}$ points.}

\refereeThird{Note that in Figure \ref{fig:method:evalpts_voronoi} the points placed into every Voronoi region are obtained using the algorithm introduced in \cite{BFF_nodes}. It is also possible 
to avoid placing $q$ points into every Voronoi region by using a global point set $Y$ that conforms to $\Omega$, which we tested in \cite{tominec2020squares}. The only requirement of the global $Y$ point set is 
that every Voronoi cell gets a reasonable sampling, with approximately $q$ points. Halton point layout would be a good candidate for forming a global $Y$. 
In this paper, however, we use the Voronoi cell point placement technique. 
}

\subsection{The RBF-FD trial space}
Let $\Omega_i$ be a subdomain holding a collection of points $X_{\Omega_i} = \{x_j\}_{j=1}^n$ that are a subset of the interpolation points
placed on top of the computational domain $\Omega$. Then every stencil is defined by a tuple $(\Omega_i, x_i)$, where
$x_i \in X_{\Omega_i}$ is the stencil center point (see Figure \ref{fig:intro:butterfly} for a visual representation of a stencil).
\refereeSecond{The solution $u_h(x)$ over a stencil is spanned by a combination of cubic polyharmonic splines
$\phi_l(x) = ||x-x_l||^3$ and multivariate monomials $\{\refereeSecond{\bar p}_k\}_{k=1}^{m}$ of degree $D_m$. The relation between $m$ and $D_m$ is: $m = \binom{D_m+d}{d}$. The linear combination then reads:}
\begin{eqnarray}
	\label{eq:method:phspolybasis}
	u_h(x)&=&\sum_{l=1}^n c_l \phi_l(x) + \sum_{k=1}^{m} \beta_k \refereeSecond{\bar p}_k(x),\\
	&&\text{ subject to  } \sum_{l=1}^n c_l \refereeSecond{\bar p}_k(x_l)=0,\quad k=1,..,m, \nonumber
\end{eqnarray}
where $n$ is the stencil size, $c_l$ are the interpolation coefficients and $\beta_k$ are the Lagrange multipliers. 
\refereeSecond{The approach with appending the monomial basis to the polyharmonic spline basis was studied in \cite{Barnett15,Bayona19,BFFB17}}.
The coefficients $c_l$ and $\beta_k$ from \eqref{eq:method:phspolybasis} are computed by requiring the
interpolation conditions
$u_h(X_{\Omega_i}) = u(X_{\Omega_i})$, where $u(X_{\Omega_i}):=\underline{u}^{(i)}$ are \refereeSecond{the stencil data (stencil nodal values)}.
This results in a square system
of equations:
\begin{equation}
	\label{eq:M}
	\underbrace{\begin{pmatrix}
		A^{\refereeSecond{(i)}} & P^{\refereeSecond{(i)}} \\
		(P^{\refereeSecond{(i)}})^T & 0
	\end{pmatrix}}_{:=\tilde A^{(i)}}
	\begin{pmatrix}
		\underline{c}^{\refereeSecond{(i)}}\\
		\underline{\beta}^{\refereeSecond{(i)}}
        \end{pmatrix}
	=
        \begin{pmatrix}
          \underline{u}^{\refereeSecond{(i)}}\\
		0
	\end{pmatrix}.
\end{equation}
Here $A_{jl}^{\refereeSecond{(i)}} = \phi_l(x_j)$ for indices $j,l=1,..,n$ and $P_{jk}^{\refereeSecond{(i)}} = \refereeSecond{\bar p}_k(x_j)$ for index $k=1,..,m$ and $u_j^{\refereeSecond{(i)}}=u^{\refereeSecond{(i)}}(x_j)$.

The stencil-based solution in any point $z \in \Omega_i$ is expressed by reusing the computed coefficients
$\underline{c}^{\refereeSecond{(i)}}$ and $\underline{\beta}^{\refereeSecond{(i)}}$ from \eqref{eq:M} within the linear combination \eqref{eq:method:phspolybasis}:
\begin{eqnarray}
	\label{eq:method:stencil_cardinal}
	u_h^{(i)}(\refereeSecond{z}) &=&\left(
	\underbrace{
	  \begin{pmatrix}
		\phi_1(\refereeSecond{z}),..,\phi_n(\refereeSecond{z}), & \refereeSecond{\bar p}_1(\refereeSecond{z}),..,\refereeSecond{\bar p}_m(\refereeSecond{z})
	\end{pmatrix}}_{:=b^{(i)}(\refereeSecond{z})}
	\begin{pmatrix}
		  \underline{c}^{\refereeSecond{(i)}}\\
		\underline{\beta}^{\refereeSecond{(i)}}
	\end{pmatrix}\right)_{1:n}
	\nonumber \\
	&=&\left(b^{(i)}(\refereeSecond{z}) ((\tilde A^{\refereeSecond{(i)}})^{-1}\right)_{1:n}\, \underline{u}^{(i)} \equiv(\refereeThird{\psi}^{(i)}_1(\refereeSecond{z})\cdots \refereeThird{\psi}^{(i)}_n(\refereeSecond{z}))\,\underline{u}^{(i)}\equiv \underline w^{(i)}(\refereeSecond{z})\,\underline{u}^{(i)},
	\label{eq:w2}
\end{eqnarray}
where
$\refereeThird{\psi}_1^{\refereeSecond{(i)}}(\refereeSecond{z}),..,\refereeThird{\psi}_n^{\refereeSecond{(i)}}(\refereeSecond{z})$ are the stencil-based cardinal functions 
with the Kronecker delta property \refereeThird{($\psi^{(i)}_k(z) = 1$ if $z=x_k$, otherwise $0$)} and $\underline w^{\refereeSecond{(i)}}(\refereeSecond{z})$ are the
local stencil weights \refereeSecond{for evaluating the stencil-based solution at a point $\refereeSecond{z}$}.

The next step is to use the formulation from \eqref{eq:method:stencil_cardinal}
to represent the solution over the whole computational domain $\Omega$.
\refereeSecond{First, the inverse $(\tilde A^{\refereeSecond{(i)}})^{-1}$} from \eqref{eq:M} \refereeSecond{is} computed for every stencil $(\Omega_i,x_i)$. 
\refereeSecond{Among all ($N$) available stencils, we now associate every evaluation point $y \in Y$ with an index of the closest stencil center point defined as:}
\begin{equation}
	\label{eq:method:stencilSelection}
	\rho(y) = \argmin_{i}\|y-x_i\|,\quad i=1,..,N.
\end{equation}
\refereeSecond{By construction, we have that $y \in \Omega_{\rho(y)}$. Thus, we can use \eqref{eq:method:stencil_cardinal} and employ the stencil selection criterion \eqref{eq:method:stencilSelection} 
to express the local solution at every 
$y \in Y$:
\begin{equation}
	\label{eq:method:localAnsatz}
u^{\rho(y)}_h(y) = \underline w^{\rho(y)}(y)\,\underline{u}^{\rho(y)} = \sum_{j=1}^n w_{j}^{\rho(y)}(y)\, u^{\rho(y)}(x_{j}),\quad x_{j} \in X_{\Omega_{\rho(y)}}.
\end{equation}
The local solution is using local indexing over the vector of weights and over the vector of the stencil nodal values. To solve a PDE, we need a global representation of the solution, that is, 
a linear combination of the global nodal values and the global weights (global cardinal functions).
}
\refereeSecond{Using \eqref{eq:method:localAnsatz},} the global solution at any $y\in \Omega$ is formally written as:
\begin{eqnarray}
	\label{eq:cardinalF}
u_h(y) = \sum_{i=1}^N u(x_i) \Psi_i(y) \nonumber &=& \sum_{i=1}^N u(x_i) \sum_{j=1}^n w^{\rho(y)}_j\refereeSecond{(y)}\, \delta_{i,\Gamma(j,\rho(y))} \nonumber \\
&=& \sum_{i=1}^N u(x_i) \sum_{j=1}^n \left[b^{\refereeSecond{\rho(y)}}(y)\tilde{A}_{\rho(y)}^{-1}\right]_j \delta_{i,\Gamma(j,\rho(y))},
\end{eqnarray}
where $\Psi_i$  are the global cardinal functions, $w_j^{\rho(y)}$ are the local weights over the stencil centered at $x_{\rho(y)}$.
Furthermore the operator
$\Gamma(j, \rho(y)): \mathcal{I}[1,n] \to \mathcal{I}[1,N]$
 is an index mapping from the $j-th$
local weight of a stencil with index $\rho(y)$ to its global equivalent.
\refereeThird{Here the sum $ \sum_{j=1}^n w^{\rho(y)}_j\refereeSecond{(y)}\, \delta_{i,\Gamma(j,\rho(y))}$ is using the Kroenecker delta function to search 
for one local weight $w^{\rho(y)}_j\refereeSecond{(y)}$ which is multiplied only with the corresponding 
global nodal solution $u(x_i)$ inside the linear combination \eqref{eq:cardinalF}. 
In this way we are mimicking the matrix-vector product of one row of the global evaluation matrix and the vector of weights generated 
specifically for $y$. This ansatz could be used for theoretical studies of the RBF-FD method in the future.}
\refereeSecond{Although the global RBF-FD cardinal functions are discontinuous along the interfaces of the Voronoi regions \cite{tominec2020squares}, 
a derivative of the unknown solution is determined in an unique way: for every given $y$ we choose an unique stencil representation \eqref{eq:method:localAnsatz} with an index $\rho(y)$ defined in \eqref{eq:method:stencilSelection}. 
Within that stencil we are always able to smoothly evaluate a derivative.}
\refereeSecond{A representation of a differential operator $\mathcal{L}$} in $y$ is then written as:
\begin{eqnarray}
	\label{eq:cardinalF_diff}
\mathcal{L} u_h(y) &=& \sum_{i=1}^N u(x_i) \mathcal{L}\Psi_i(y) \nonumber \\
&=& \sum_{i=1}^N u(x_i) \sum_{j=1}^n \left[\mathcal{L}b^{(i)}(y)\,\tilde{A}_{\rho(y)}^{-1}\right]_j \delta_{i,\Gamma(j,\rho(y))},
\end{eqnarray}
\subsection{Evaluation and differentiation matrices}
Setting $y=Y$ in equation \eqref{eq:cardinalF} we arrive at the discrete representation of a solution in the $Y$-points:
\begin{equation}
\label{eq:evalMatrix}
u_h(Y) = E_h(Y,X) u_h(X),
\end{equation}
where $E_h(Y,X)$ is a rectangular evaluation matrix interpolating $u_h(X)$ from $X$ to $Y$. Its components are $(E_h)_{ik} = \Psi_i(y_k)$.

The discrete representation of a differential operator $\mathcal{L}$ is obtained by setting $y=Y$ in \eqref{eq:cardinalF_diff}:
\begin{equation}
    \label{eq:diffMatrix}
    \mathcal{L} u_h(Y) = D_h^{\mathcal{L}}(Y,X) u_h(X),
    \end{equation}
where $D_h^{\mathcal{L}}$ is a rectangular differentiation matrix with components $(D^{\mathcal{L}}_h)_{ik} = \mathcal{L}\Psi_i(y_k)$. 
MATLAB code for generating matrices $E_h(Y,X)$ and $D_h^{\mathcal{L}}(Y,X)$ is available in \cite{Tominec_rbffdcode2021}.

\subsection{The unfitted discretization of a PDE}
The model problem \eqref{eq:model:Poisson_system} is discretized with the RBF-FD operators $E_h$ and $D_h^{\mathcal{L}}$
given in
\eqref{eq:evalMatrix} and \eqref{eq:diffMatrix} respectively. The result is the semi-discrete matrix $D_h(y,X)$ and the semi-discrete vector $F(y)$,
where:
\begin{eqnarray}
\label{eq:unfitted:Poisson_system}
D_h(y,X) =
\left\{
\begin{array}{ll}
	\beta_2 D_h^\Delta(y,X), & y \in \Omega \\
	\beta_0 E_h(y,X), & y \in \partial\Omega_0 \\
	\beta_1 D_h^{\nabla\cdot n(y)}(y,X), & y \in \partial\Omega_1 \\
\end{array}
\right. 
 \mbox{ and }
F(y) =
\left\{
\begin{array}{ll}
	\beta_2 f_2(y), & y \in \Omega \\
	\beta_0 f_0(y), & y \in \partial\Omega_0 \\
	\beta_1 f_1(y), & y \in \partial\Omega_1 \\
\end{array}
\right. \nonumber \\
\end{eqnarray}
Here $\beta_2, \beta_0, \beta_1$ are the scalings of the PDE and the boundary conditions.
Setting $y=Y$ we obtain a rectangular linear system of size $M \times N$:
\begin{equation}
	\label{eq:method:PDE_discrete}
D_h(Y,X) u_h(X) = F(Y).
\end{equation}
For that system choose the scalings:
$$\beta_2 = \frac{1}{\sqrt{M_2}},\quad \beta_0 = h^{-1} \frac{1}{\sqrt{M_0}},\quad \beta_1 = \frac{1}{\sqrt{M_1}},$$
where $h$ is the average node distance in the node set $X$ and where $M_2$, $M_0$ and $M_1$ are the number of evaluation points placed 
over $\Omega$, $\partial\Omega_0$ and $\partial\Omega_1$ respectively.

The motivation to use this scaling is two fold. Firstly, the factors $\frac{1}{\sqrt{M_2}}$, $\frac{1}{\sqrt{M_0}}$ and
$\frac{1}{\sqrt{M_1}}$ relate every component of $D_h^T D_h$, which
are discrete inner products, to continuous inner products plus a first order integration error \cite{tominec2020squares}.
Secondly, the factor $h^{-1}$ is used to impose the Dirichlet condition in a weak sense such that the matrix \refereeSecond{$D_h^T D_h$} is nonsingular: 
this is a classical approach in those finite element methods
which use a solution space that does not exactly satisfy the Dirichlet condition.
In order to prove uniqueness of the solution in such a setup, a mesh
dependent parameter $h^{-1}$ has to be introduced via inverse inequalities, which is at the end multiplying the
added Dirichlet penalty term \cite{Babuska}. We also note that we are not able to impose a Dirichlet condition exactly in an efficient way
due to using $X$-points that are unfitted with respect to the boundary.

The numerical solution is obtained by solving \eqref{eq:method:PDE_discrete} for the solution coefficients $u_h(X)$ and then interpolating
this data onto the evaluation points $Y$.
This can be written as:
\begin{equation}
	\label{eq:method:PDE_solution}
u_h(Y) = E_h(Y,X)\, D_h^+(Y,X) F(Y),
\end{equation}
where $D_h^+(Y,X) = (D_h^T D_h)^{-1} D_h^T$ is \refereeSecond{a pseudoinverse, which, in practice, is} computed using the QR decomposition.

Once $u_h(X)$ is computed, the residual is given by:
\begin{equation}
	\label{eq:method:residual}
r(Y) = D_h(Y,X)u_h(X) - F(Y).
\end{equation}
The least-squares residual is by definition orthogonal to the column-space of $D_h(Y,X)$, which implies the relation:
\begin{equation}
	\label{eq:method:residual_ortho}
D_h^T r(Y) = 0\quad \Rightarrow \quad D_h^+ r(Y) = (D_h^T D_h)^{-1} D_h^T r(Y) = 0.
\end{equation}

\refereeSecond{Throughout this paper we solve \eqref{eq:method:PDE_solution} using the \emph{mldivide()} function in Matlab, which uses a sparse 
QR decomposition as an intermediate step in order to obtain the solution to the 
rectangular system. In our experience this was a reasonably fast approach when the number of unknowns was -- roughly speaking -- smaller or equal to $10^5$. For problems with larger amounts of unknowns, 
the reader could use a conjugate-gradient method on the normal system $D_h^T D_h = D_h^T F$. Unfortunately this approach is prone to numerical instabilities when $D_h$ is ill-conditioned, since 
the already large condition number is then squared: $\kappa(D_h^T D_h) = \kappa(D_h)^2$. However, there exists a numerically more stable algorithm, called \emph{lsqr} \cite{lsqr}, 
which is a conjugate-gradient like iterative solver designed for solving large rectangular systems of equations.}

%
\section{Linear independence of cardinal functions}
\label{sec:method:linearindep}
In this section we address the difficulties related to the linear independence of the cardinal functions which arise when using the unfitted discretization.
We provide a criterion upon which a certain amount
of the interpolation points that extend outside of the computational domain is removed. A similar study, but in a context of
the isogeometric finite element method is performed in
\cite{elfverson2018cutiga}.

\refereeSecond{Throughout Section \ref{section:method} we outlined that we are using a solution ansatz 
\eqref{eq:cardinalF} plugged into the PDE problem \eqref{eq:model:Poisson_system} to then solve the discretized PDE problem 
\eqref{eq:method:PDE_discrete} for the unknown nodal values 
$u(x_i),\, i=1,..,N$. It is not possible to solve \eqref{eq:method:PDE_discrete} for these nodal values, 
unless the cardinal functions $\Psi_i(y), i=1,..,N$ in \eqref{eq:cardinalF} form a basis for $\mathbb{R^N}$, 
i.e. unless they are linearly independent with each other.}
Since \refereeSecond{the columns} of the matrix $E_h$ from \eqref{eq:evalMatrix} contain \refereeSecond{the whole set of sampled cardinal functions}, 
we investigate the linear independence
of those columns by computing the smallest singular value $\sigma_{\min}(E_h)$.
When \refereeSecond{$\sigma_{\min}(E_h)=0$} we have that the columns have a nonzero nullspace and
thus the sampled family of cardinal functions is linearly dependent.
\refereeSecond{In this section we do not focus on studying $\sigma_{\min}(D_h)$, as this approach would in addition carry information on 
the well-posedness of the PDE problem. For example, if we had observed $\sigma_{\min}(D_h)=0$, 
then this could well be due to using an unfitted discretization, but also other reasons such as: 
ill-posedness of the PDE problem, wrong imposition of boundary conditions, etc.}

Whereas the RBF-FD cardinal functions \refereeSecond{in the fitted setting} are indeed linearly independent when
the interpolation points conform to $\Omega$ and the interpolation points $X$ are a subset of the evaluation points 
$Y$ \cite{tominec2020squares}, it is important to check whether this is true in the unfitted case
as well.

A cardinal function $\Psi_k(y)$ has a Kronecker delta property in the $X$-points:
\begin{eqnarray*}
	\Psi_k(y) =
	\left\{
	\begin{array}{ll}
		1, & y=x_k \subset X \\
		0 & y\in X\setminus x_k,
		\end{array}
	\right.
\end{eqnarray*}
which guarantees linear independence as long as $X \subseteq Y$ and $X\subseteq \Omega$ since in this case,
there is always at least one point in $\Omega$ for every cardinal function (e.g. $x_k$ for $\Psi_k$)
 where $\Psi_k$ is one, but all other $\Psi_j$ for indices $k \neq j $ are $0$.
A problem when using the unfitted discretization
can occur due to the compact support of $\Psi_k$. When an external $x_k$ is placed such that its corresponding $\Psi_k$
vanishes before it reaches the interior of $\Omega$, then $\Psi_k(y)=0$ for every $y \in \Omega$ and the basis function becomes linearly
dependent (on $\Omega$) with all others. 
\refereeSecond{Note that when $\Psi_k(y)=0$ for every $y \in \Omega$, we also have that all derivatives of $\Psi_k(y)$ are $0$ for every $y \in \Omega$. Thus 
the $k$-th column of the PDE matrix $D_h$, which is built upon a combination of $\Psi_k$ and its derivatives, is a $0$ vector: this implies that $D_h$ in this case does not have a full column rank.}

\begin{figure}[h!]
	\centering
\includegraphics[width=0.49\linewidth]{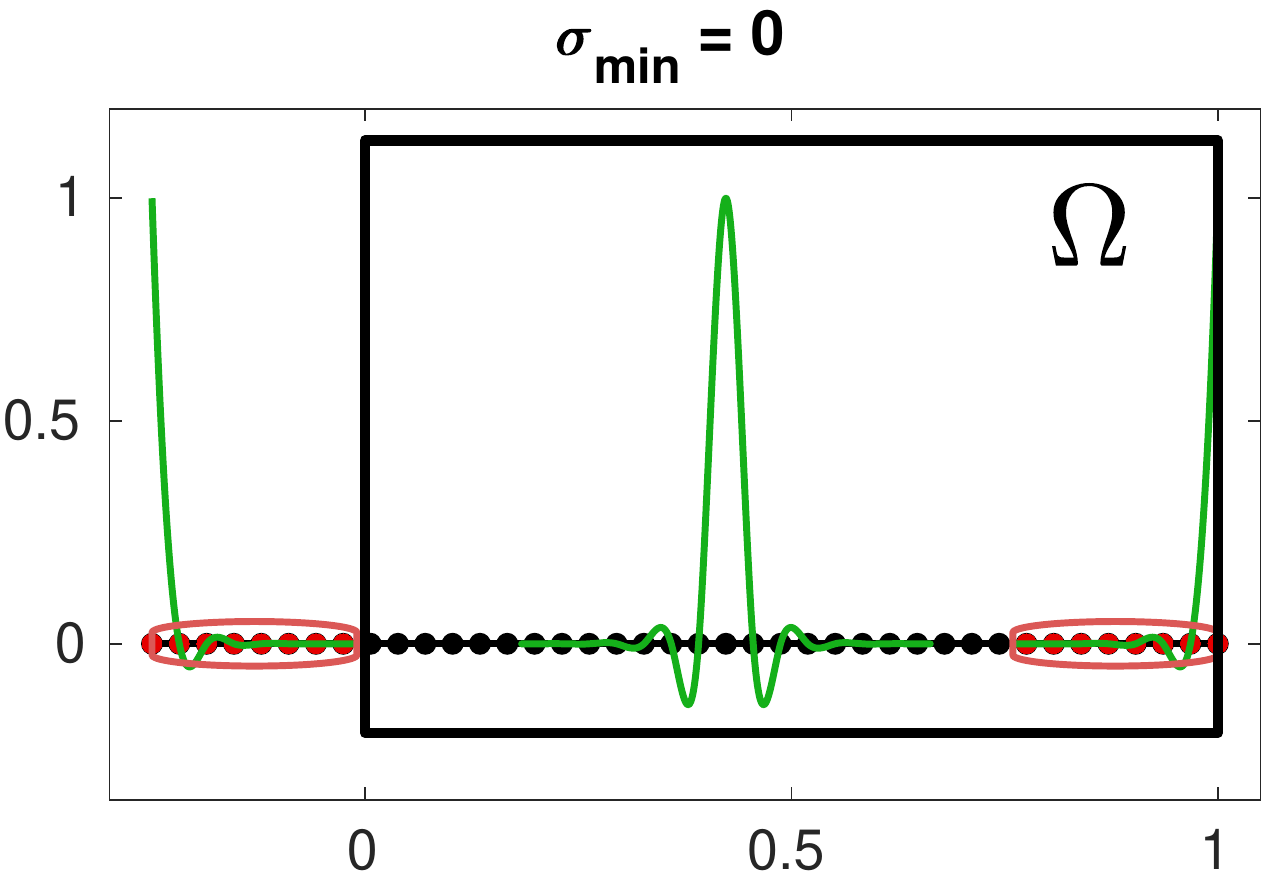} \\
\includegraphics[width=0.49\linewidth]{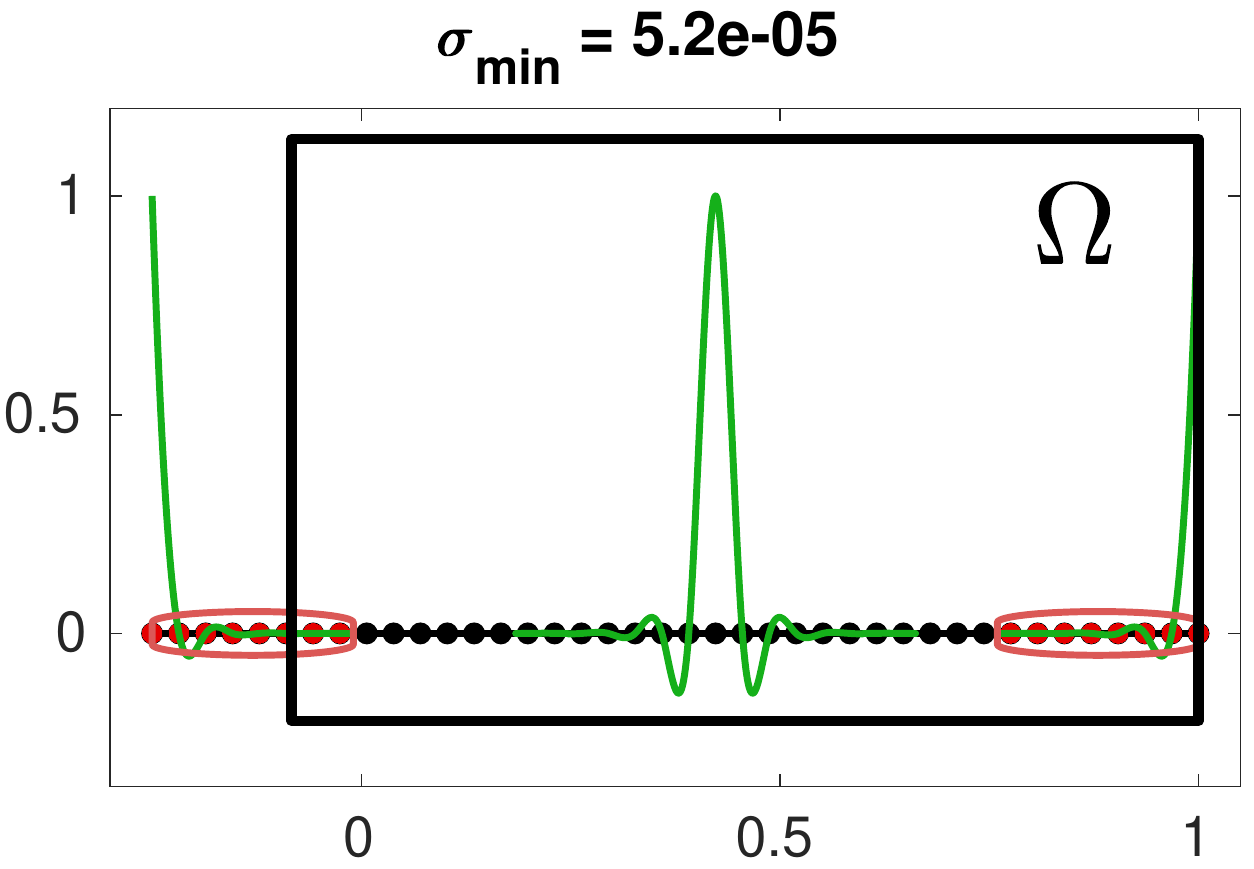}
\includegraphics[width=0.49\linewidth]{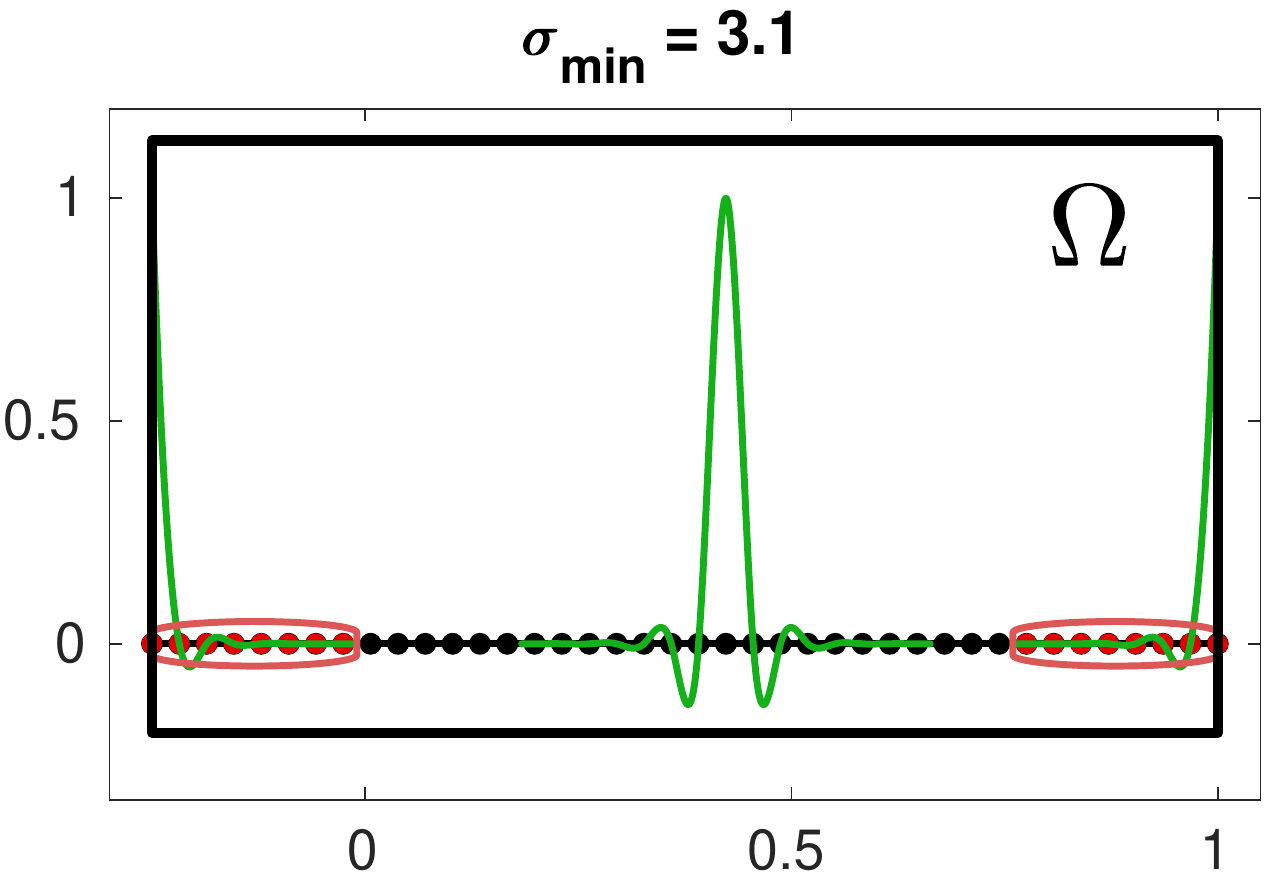}
\caption{Setup for measuring the smallest singular value of the evaluation matrix $E_h$ as the size of the computational domain $\Omega$ grows.
The green lines represent three cardinal functions. The encircled red points represent the support of the two outmost cardinal functions.}
\label{fig:method:cardinalTest_graphics_1d}
\end{figure}

In Figure \ref{fig:method:cardinalTest_graphics_1d} we can see a one-dimensional setup, where in the top plot, the $X$-points
are placed outside of $\Omega$ such that the left-most cardinal function $\Psi_{\mathrm{left}}$ does not have a support (red ellipse) in $\Omega$, which
results in a singular $E_h$ ($\sigma_{\min}=0$). As the support of $\Psi_{\mathrm{left}}$ enters $\Omega$, then $E_h$ becomes non-singular
($\sigma_{\min}=5.2\cdot 10^{-5}$),
and when the support is fully contained inside $\Omega$ then the smallest singular value gets considerably larger ($\sigma_{\min}=3.1$).

\begin{figure}[h!]
	\centering
\includegraphics[width=0.325\linewidth]{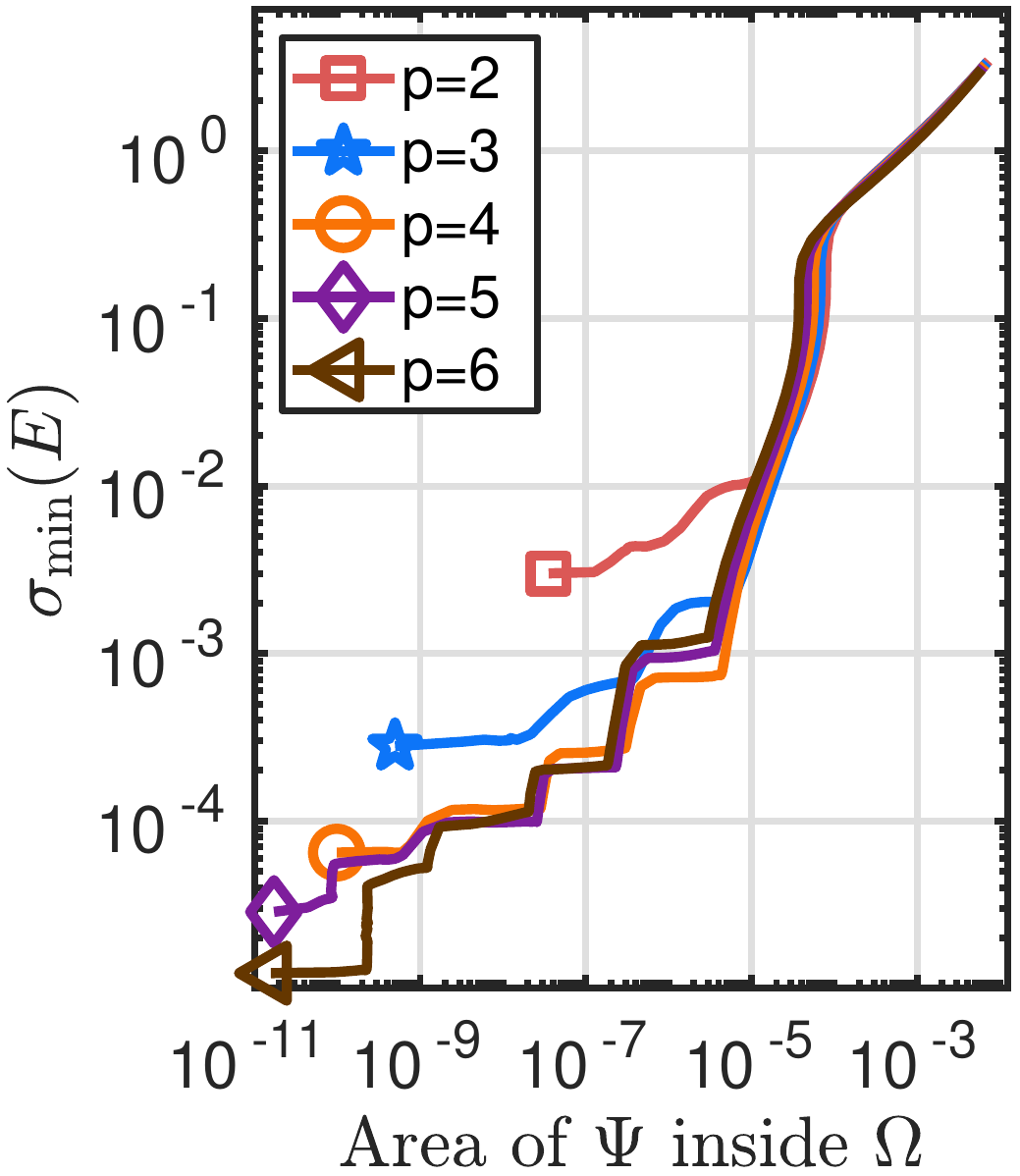}
\includegraphics[width=0.325\linewidth]{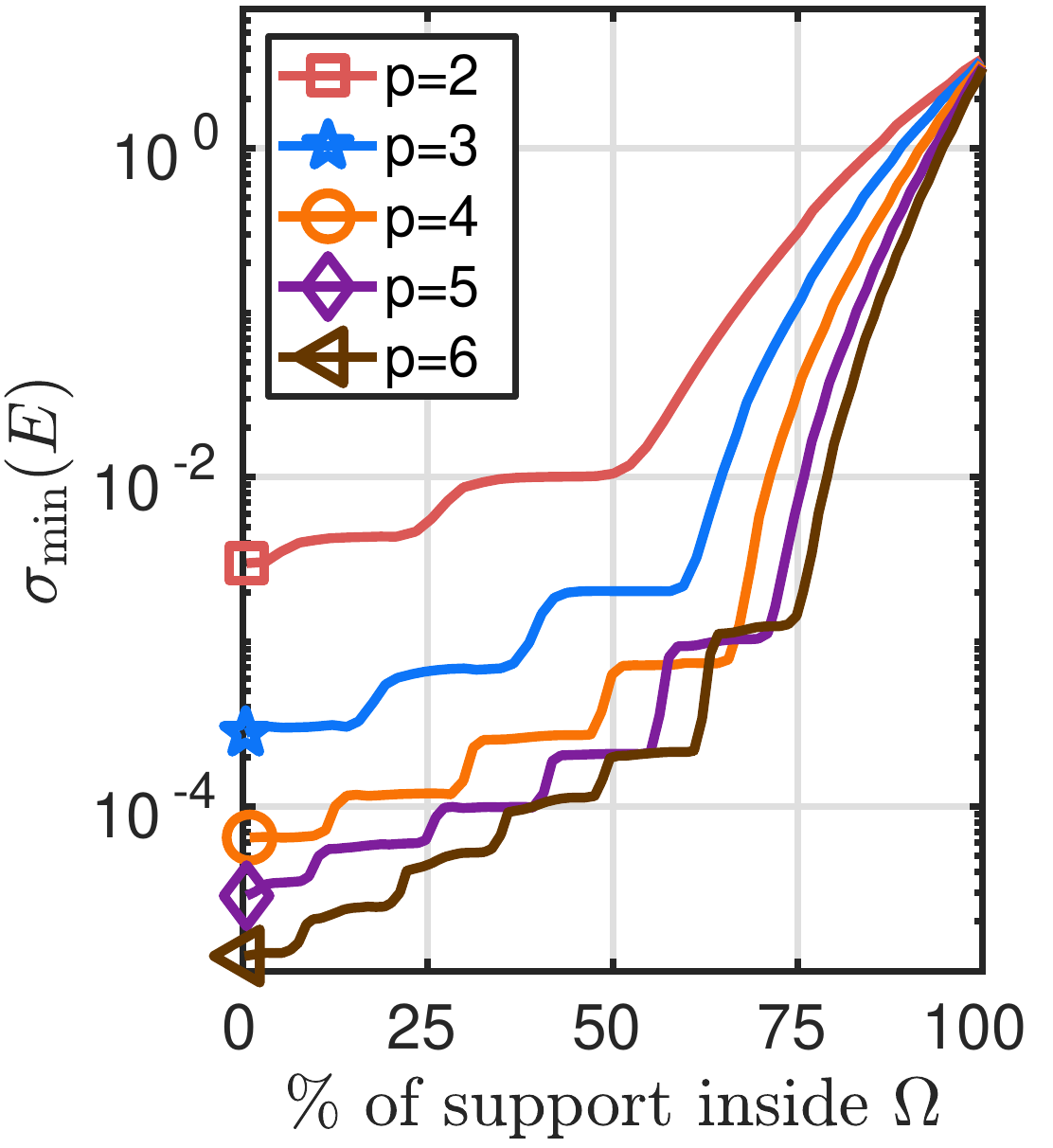}
\includegraphics[width=0.325\linewidth]{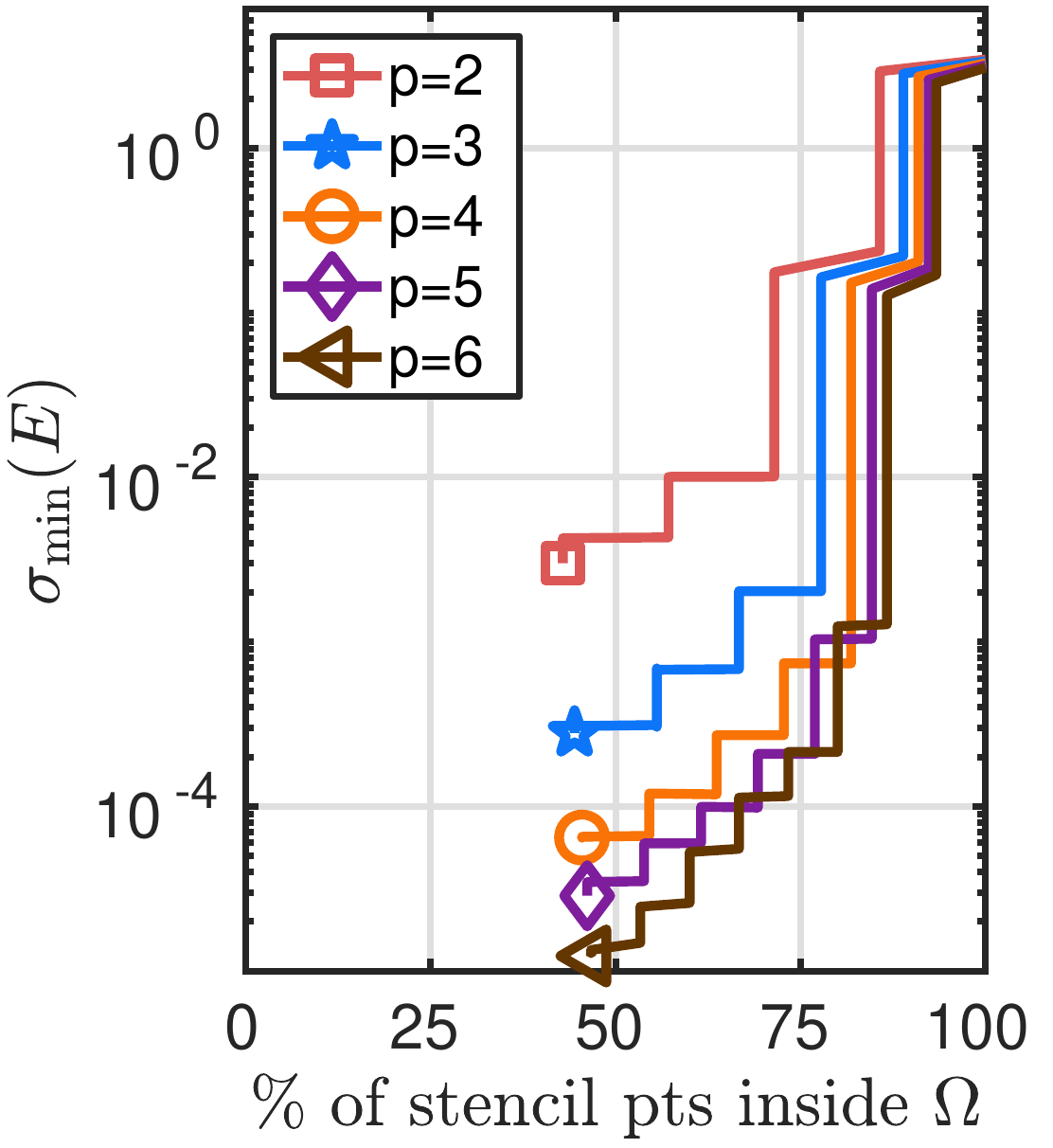}
\caption{One-dimensional case: The smallest singular value as a function of (i) the approximate area of the outmost cardinal function $\Psi_{\text{outmost}}$ 
that penetrates inside $\Omega$,
(ii) the support of $\Psi_{\text{outmost}}$ inside $\Omega$, (iii) the percentage of stencil points inside $\Omega$, which are a part of the outmost stencil.}
\label{fig:method:cardinalTest_1d}
\end{figure}

The same setup is used in Figure \ref{fig:method:cardinalTest_1d}, where for different polynomial degrees $p$, $\sigma_{\min}(E_h)$ is
computed as a function of:
\begin{itemize}
\item the approximate area under $\Psi_{\mathrm{left}}$ inside $\Omega$,
\item the percentage of the compact support of $\Psi_{\mathrm{left}}$ inside $\Omega$,
\item the percentage of stencil points of the left-most stencil inside $\Omega$.
\end{itemize}

From Figure \ref{fig:method:cardinalTest_1d} we can see that as long as the percentage of support inside $\Omega$ is larger
than $0$, $\sigma_{\min}$ is also larger than
$0$. When the percentage of support is gradually increased $\sigma_{\min}$ is also increased, approximately in the same way
as the area of $\Psi_{\mathrm{left}}$ inside $\Omega$. In the third plot we can see that when the percentage of the stencil
points inside $\Omega$ is more than $50\%$, $E_h$ is always non-singular given that the cardinal functions
(columns of $E_h$) are well sampled. We note that any cardinal function $\Psi_k$ centered in $x_k$ is not genereated by a
single stencil, but by several stencils which have $x_k$ as a neighboring point. In this sense the compact support of $\Psi_k$
is decoupled from the support of one stencil centered around $x_k$. This is the reason why the measurements of $\sigma_{\min}$ do
not coincide when considered as a function of a compact support inside $\Omega$ and stencil points inside $\Omega$.

\begin{figure}[h!]
	\centering
\includegraphics[width=0.325\linewidth]{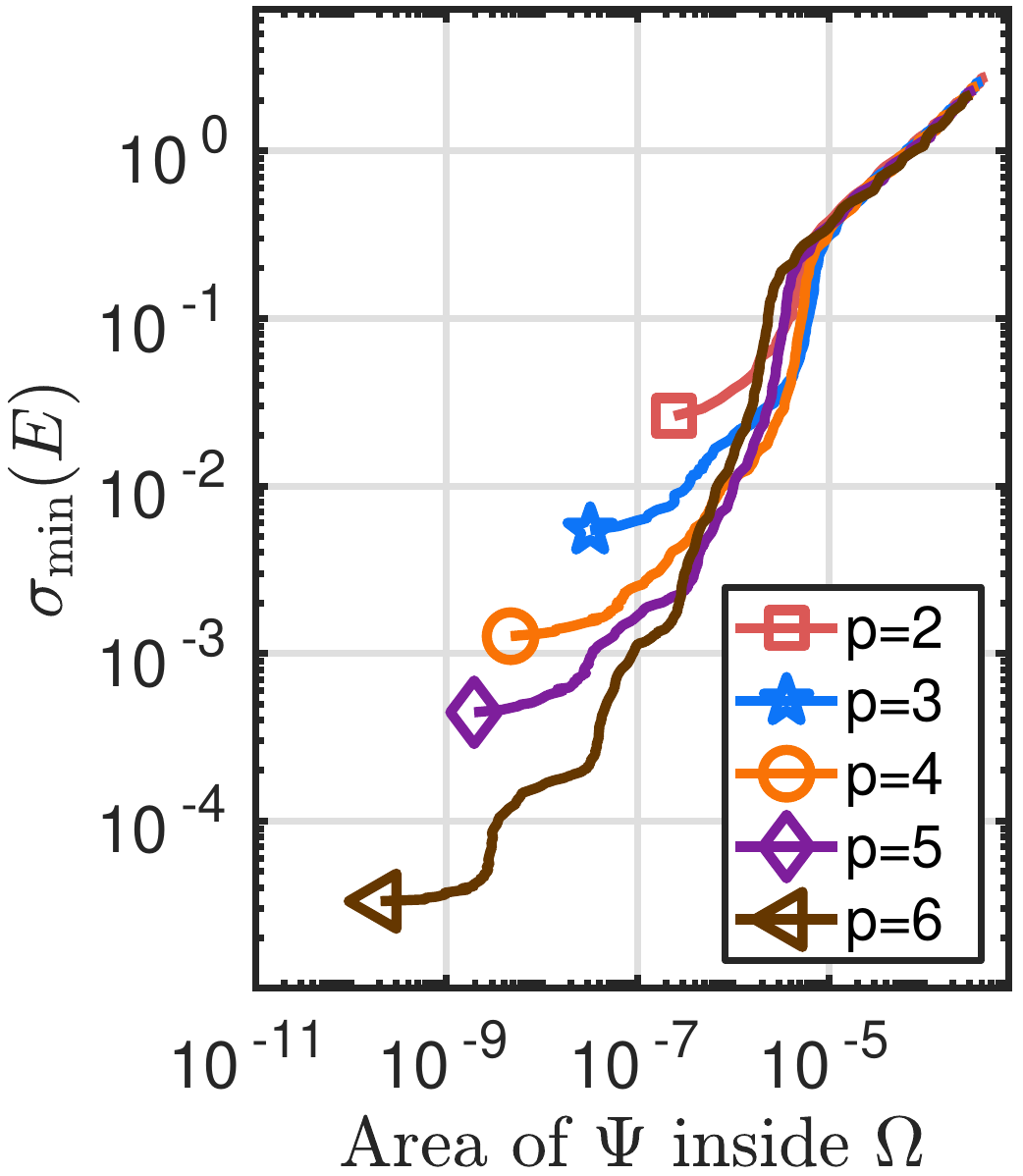}
\includegraphics[width=0.325\linewidth]{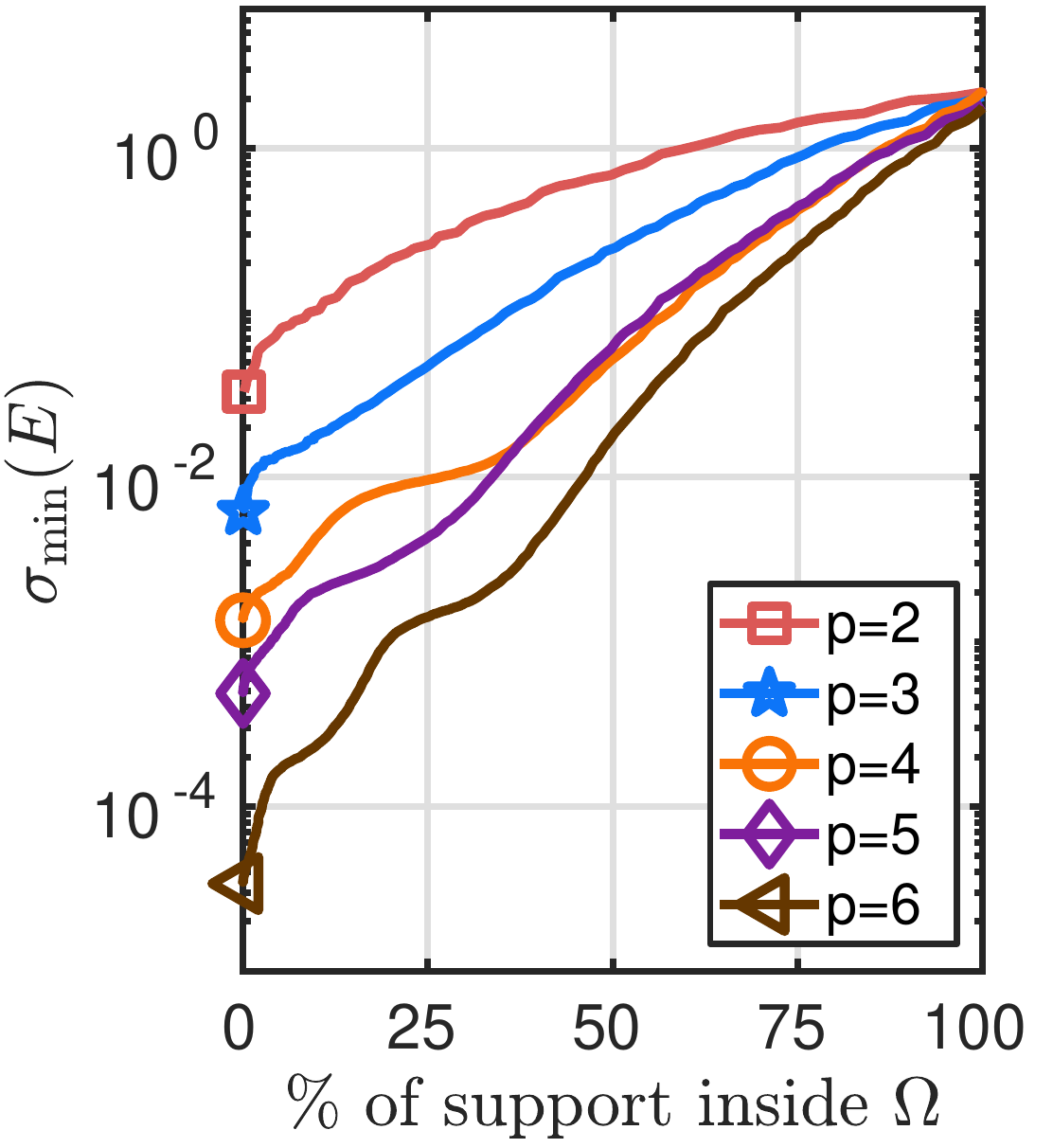}
\includegraphics[width=0.325\linewidth]{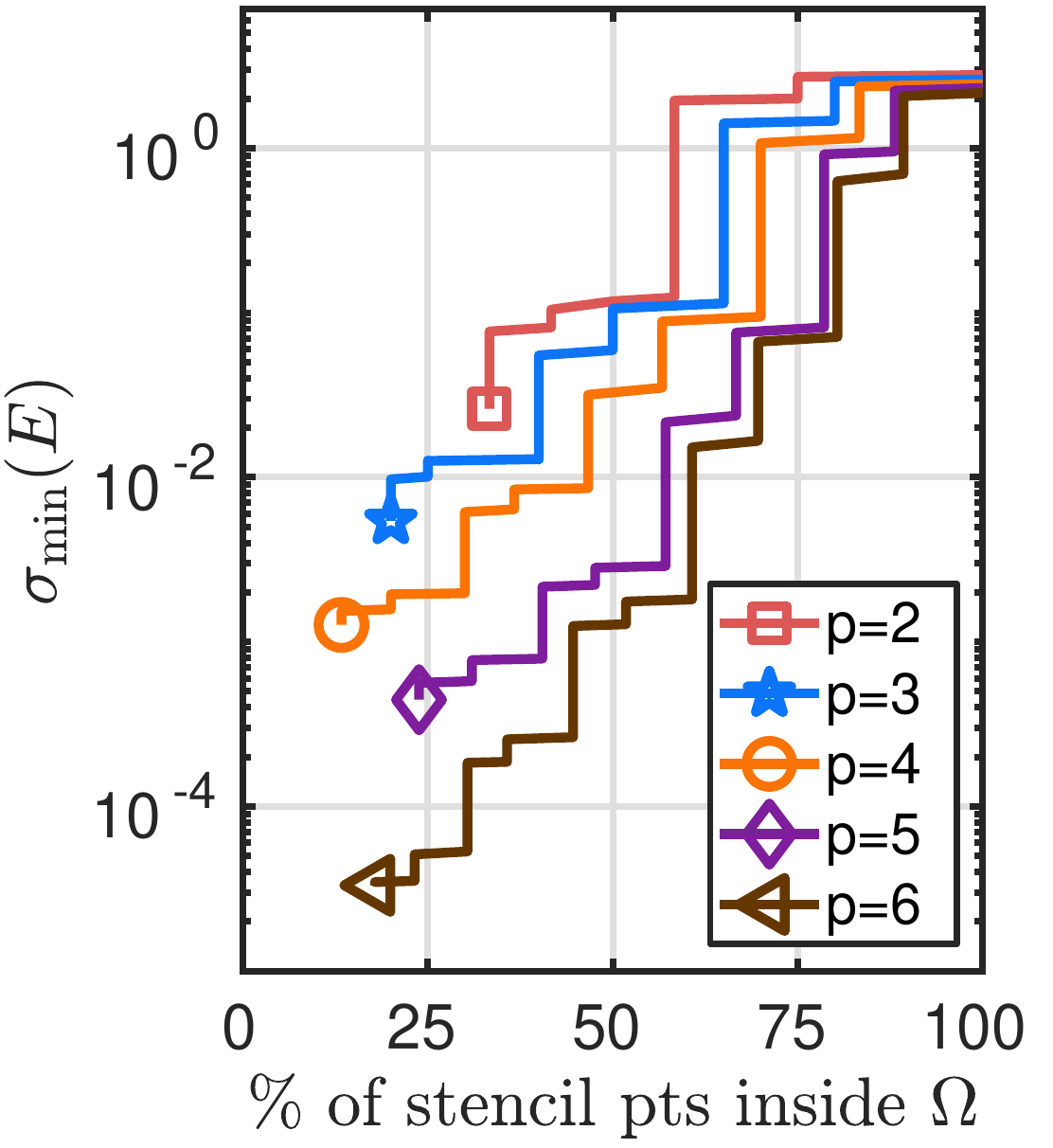}
\caption{Two-dimensional case: The smallest singular value as a function of (i) the approximate area of the outmost cardinal function $\Psi_{\text{outmost}}$ 
that penetrates inside $\Omega$,
(ii) the support of $\Psi_{\text{outmost}}$ inside $\Omega$, (iii) the percentage of stencil points inside $\Omega$, which are a part of the outmost stencil.}
\label{fig:method:cardinalTest_2d}
\end{figure}

Figure \ref{fig:method:cardinalTest_2d} shows analogous results to Figure \ref{fig:method:cardinalTest_1d} but for a
two-dimensional case, where the computational domain
is a square of size $[-1,1] \times [-1,1]$. Similar results can be observed as for the
one-dimensional case in terms of the area under $\Psi_{\mathrm{left}}$ and the percentage of its support inside $\Omega$, however,
the percentage of stencil points inside $\Omega$ is allowed to be smaller in the 2-dimensional case.
This indicates that the relation between
the support of a stencil and the compact support of a cardinal function is in this case tighter.

In our experience the smallest percentage of stencil points inside $\Omega$ is the criterion which is the
easiest to implement prior to computing RBF-FD differentiation matrices and it is therefore our choice for all further experiments.
We remove some of the initial interpolation points (placed in the box around $\Omega$) such that at least $50\%$ of points of each stencil are contained inside $\Omega$.
Once we decide on a stencil size $n$ and the initial $X$-points and $Y$-points are computed, the criterion can be used by invoking
one command in Matlab.

\begin{verbatim}
X = X(unique(knnsearch(X,Y,'k',ceil(0.5*n)), :);
\end{verbatim}
\section{Analysis of the approximation error under node refinement}
\label{section:erroranalysis}
In this section we develop an understanding of the behavior of the error $e$ between the true solution $u(Y)$ and the approximate solution $u_h(Y)$, both restricted to
the evaluation points.
\subsection{Preliminaries I: A norm for measuring the error}
Our choice of a vector norm that measures the error $e(Y) = u(Y) - u_h(Y)$ is given by:
\begin{equation}
    \label{eq:theory:norm}
    \|e(Y)\|^2_{\ell_2} = \frac{1}{M} \sum_{j=1}^M e(y_i)^2 = \frac{1}{M} \|e(Y)\|^2_2.
\end{equation}
The $\ell_2$ norm
is a good choice for the discrete least-squares problems since it is up to $\mathcal{O}(h_y)$ ($h_y$ is the average spacing between the $Y$-points) equivalent to 
the $L_2$ norm which is a natural
norm for the continuous least-squares problem \cite{tominec2020squares}. In this sense the stability properties can carry over to the discrete formulation, under some assumptions \cite{tominec2020squares}.

\refereeSecond{
An analogous norm for a matrix $A \in \mathbb{R}^{M \times N}$ is:
\begin{equation}
        \label{eq:theory:norm_matrix}
        \|A\|^2_{\ell_2} = \frac{1}{M}\, \sup_{x \neq 0} \frac{\|A x\|^2_2}{\|x\|^2_2} = \frac{1}{M} \|A\|^2_{2}.
\end{equation}
}

Furthermore we have the relation $\|e_h(Y)\|_{\ell_2} \leq \|e_h(Y)\|_{\infty}$ due to:
\begin{equation}
    \label{eq:theory:normRelation}
    \|e(Y)\|^2_{\ell_2} = \frac{1}{M} \sum_{j=1}^M e(y_i)^2 \leq \frac{1}{M} M \max_i e(y_i)^2 = \|e(Y)\|^2_\infty
\end{equation}
which is later
used to bound the consistency terms.

\subsection{Preliminaries II: Consistency estimates}
Evaluate a function $\hat u: \hat \Omega \to \mathbb{R}$ in the $X$-points and denote this data by $\hat u(X)$.
Then we can, for any $y\in \hat \Omega$, where $\hat \Omega$ contains $\Omega$, construct a finite dimensional
representation $u_{\text{\refereeSecond{ext}}}(y)$ that interpolates the data $\hat u(X)$:
\begin{equation}
    \label{eq:theory:u_findim}
    u_{\text{\refereeSecond{ext}}}(y) = \sum_{i=1}^N \hat u(x_i) \Psi_i(y),
\end{equation}
where $\Psi_i(y)$ are the RBF-FD cardinal functions defined in \eqref{eq:cardinalF}. The equivalent matrix-vector formulation is:
\begin{equation}
    \label{eq:theory:u_findim_matrix}
u_{\text{\refereeSecond{ext}}}(y) = E_h(y,X) \hat u(X).
\end{equation}
The interpolation error
is estimated by \cite{Bayona19,tominec2020squares}:
\begin{equation}
    \label{eq:theory:interpError}
\|\hat u(y) - u_{\text{\refereeSecond{ext}}}(y)\|_\infty \leq C_I h^{p+1} |\hat u|_{W_\infty^{p+1}},
\end{equation}
where and $|\hat u|_{W_\infty^{p+1}}$ is a Sobolev semi-norm defined as
$|\hat u|_{W_\infty^{p+1}} = \max_{y\in\hat\Omega} |\mathcal{D}^{p+1} \hat u(y)|$, where $\mathcal{D}^{p+1}$ is a partial derivative 
of degree $p+1$. \refereeSecond{The constant $C_I$ depends on $p$ and the quality of $X$ point set $c_q=\frac{h}{h_s}$, where $h$ is the fill distance 
in the point set $X$ and $h_s$ is the separation 
distance in the point set $X$. The two distances are defined as: 
$$h = \sup_{x\in \Omega}\, \min_{x_j \in X}\, \| x- x_j \|_2 \leq h_s = \frac{1}{2}\, \min_{j\neq k, x_j,x_k \in X} \|x_j - x_k\|_2.$$
}

An application of a PDE operator $D$ from \eqref{eq:model:Poisson_system} to \eqref{eq:theory:u_findim} gives:
$$D u_{\text{\refereeSecond{ext}}}(y) = \sum_{i=1}^N \hat u(x_i) D \Psi_i(y).$$
The (semi-discrete) matrix-vector equivalent is:
\begin{equation}
    \label{eq:theory:Du_findim_matrix}
D u_{\text{\refereeSecond{ext}}}(y) = D_h(y,X) \hat u(X).
\end{equation}
There are three parts involved in $D$, see \eqref{eq:model:Poisson_system}, the Laplacian, the normal derivative and the Dirichlet condition.
Each of them has its own consistency estimate depending on the order of the derivative. The overall
consistency \cite{tominec2020squares} is bounded by:
\begin{eqnarray}
    \label{eq:theory:diffError}
\|D \hat u(y) - Du_{\text{\refereeSecond{ext}}}(y)\|_\infty &\leq& \left(C_2 h^{p-1} + C_1 h^p + C_0 h^{p+1}\right) |\hat u|_{W_\infty^{p+1}} \nonumber \\
&\leq& C h^{p-1} |\hat u|_{W_\infty^{p+1}}.
\end{eqnarray}

\subsection{Preliminaries III: Extension of the true solution}
Generally speaking the error estimation for the unfitted RBF-FD requires a special treatment
since an intermediate numerical solution $u_h(X)$, see \eqref{eq:method:PDE_solution}, 
also lives in the exterior of the computational domain.
In order to be able to compare the solution with the true solution we have to define an extension of the true solution on some extended open domain.
Let $u=u(y)$ be a true unique solution of the PDE problem \eqref{eq:model:Poisson}, where $y\in \bar\Omega$ (the closure of $\Omega$) and
let $\hat\Omega = (\bar \Omega + \tilde{\Omega}) \subset \mathbb{R}^2$ be an extended domain. The extended smooth solution $\hat{u}=\hat u (\hat y)$,
where $\hat y\in (\bar \Omega + \tilde{\Omega})$ is then defined such that:
\begin{equation}
    \label{eq:theory:u_extended}
\hat{u}|_{\bar \Omega} = u.
\end{equation}
This definition makes it possible to bound the error in terms of the stability and consistency terms,
where the latter then depends on the size of the partial derivative of degree $p+1$.
There exists an extension lemma for smooth functions that we can use, which is stated below.
\begin{lemma}[\protect{\cite[Lemma 4.1]{Lee}}]
    Suppose M is a smooth
    manifold with or without boundary, $A \subseteq M$ is a closed subset, and $f:A\to\mathbb{R}^k$ is a
    smooth function. For any open subset $U$ containing $A$, there exists a smooth function
    $\tilde f: M \to \mathbb{R}^k$ such that $\tilde f|_A = f$ and $\text{supp } \tilde f \subseteq U$.
\end{lemma}
The lemma from above is directly applicable to extending $u: \bar\Omega \to \mathbb{R}$, since $\bar\Omega$ is
a closed domain and since our $u$ is a smooth function. We can therefore conclude that the extension defined in \eqref{eq:theory:u_extended} exists, and that
a discrete implication is the relation:
\begin{equation}
    \label{eq:theory:u_extended_discrete}
    \hat u(Y) = u(Y),
\end{equation}
where $Y$ is an evaluation point set which conforms to $\Omega$.
\subsection{The PDE error estimate}
Now we estimate the error $e = u_h(Y) - u(Y)$ between the numerical solution $u_h(Y)$ and the true solution $u(Y)$:
\begin{eqnarray}
    \label{eq:theory:estimate_tmp1}
\|u(Y) - u_h(Y)\|_{\ell_2} &=& \|u(Y) - u_h(Y) + u_{\text{\refereeSecond{ext}}}(Y) - u_{\text{\refereeSecond{ext}}}(Y)\|_{\ell_2} \nonumber \\
&\leq& \|u_h(Y) - u_{\text{\refereeSecond{ext}}}(Y)\|_{\ell_2} + \|u(Y) - u_{\text{\refereeSecond{ext}}}(Y)\|_{\ell_2}
\end{eqnarray}
where we first added and subtracted $u_{\text{\refereeSecond{ext}}}(Y)$, then used the triangle inequality to make a split into the PDE error
$\|u_h(Y) - u_{\text{\refereeSecond{ext}}}(Y)\|_{\ell_2}$ and the
interpolation error $\|u(Y) - u_{\text{\refereeSecond{ext}}}(Y)\|_{\ell_2}$. The latter is \refereeSecond{then measured and} trivially bounded by \eqref{eq:theory:interpError}.

The term $\|u_h(Y) - u_{\text{\refereeSecond{ext}}}(Y)\|_{\ell_2}$ remains to be estimated. We first use the definition of $u_h(Y)$ from \eqref{eq:evalMatrix} and
the definition of $u_{\text{\refereeSecond{ext}}}(Y)$ from \eqref{eq:theory:u_findim_matrix} and then multiply with $D_h^+ D_h = I$ where $D_h$ defined in
\eqref{eq:method:PDE_discrete} \refereeSecond{is assumed to have a full column rank}.
\begin{align*}
\|u_h(Y) - u_{\text{\refereeSecond{ext}}}(Y)\|_{\ell_2} &= \|E_h u_h(X) - E_h \hat{u}(X)\|_{\ell_2} & \\
&= \|E_h\, D_h^+D_h\, \left( u_h(X) - \hat{u}(X) \right)\|_{\ell_2} \\
&= \|E_h\, D_h^+\, \left( D_h u_h(X) - D_h \hat{u}(X) \right)\|_{\ell_2}.
\end{align*}
Using the relation $D_h u_h(X) = F(Y) + r_h(Y)$ from \eqref{eq:method:residual} and the fact that
$D_h^+ r_h(Y) = 0$ from \eqref{eq:method:residual_ortho} we then obtain:
\begin{eqnarray*}
\|u_h(Y) - u_{\text{\refereeSecond{ext}}}(Y)\|_{\ell_2} &=& \|E_h\, D_h^+\, \left( F(Y) + r_h(Y) - D_h \hat{u}(X) \right)\|_{\ell_2} \\
&\refereeSecond{=}& \|E_h\, D_h^+ \left( F(Y)  - D_h \hat{u}(X) \right)\|_{\ell_2}.
\end{eqnarray*}
After that we use that $F(Y) = Du(Y)$ by \eqref{section:modelproblem_system}, $u(Y) = \hat u(Y)$ by \eqref{eq:theory:u_extended},
$\|E_h\,D_h^+\|_{\ell_2} = \frac{1}{\sqrt{M}} \|E_h\, D_h^+\|_2$ by the \refereeSecond{matrix norm relation from \eqref{eq:theory:norm_matrix}}
and use the \refereeSecond{submultiplicative property of the 2-norm to bound the matrix norm $\|E_h\, D_h^+\|_{2}$, and arrive to}:
\begin{align}
    \label{eq:theory:estimate_tmp2}
\|u_h(Y) - u_{\text{\refereeSecond{ext}}}(Y)\|_{\ell_2} &\leq \frac{1}{\sqrt{M}} \|E_h\|_{2}\,\|D_h^+\|_{2}\, \|D\hat{u}(Y) - D_h \hat{u}(X)\|_{\ell_2}
\end{align}

It now remains to insert the estimate \eqref{eq:theory:estimate_tmp2} into \eqref{eq:theory:estimate_tmp1} and then also combine 
this with consistency estimates \eqref{eq:theory:interpError}, \eqref{eq:theory:diffError} to arrive at the final error estimate:
\begin{eqnarray}
    \label{eq:theory:finalestimate}
    \|u(Y) - u_h(Y)\|_{\ell_2} &\leq& \frac{1}{\sqrt{M}} \|E_h\|_{2}\,\|D_h^+\|_{2}\, C_D\, h^{p-1}\, \max|D^{p+1}\hat{u}| + C_E\, h^{p+1}\, \max|D^{p+1}\hat{u}|. \nonumber \\
\end{eqnarray}
The term $\frac{1}{\sqrt{M}} \|E_h\|_{2}\,\|D_h^+\|_{2}$ is what we call the stability norm, which should remain constant
so that the error overall decays with at least order $p-1$. We numerically test that in the following section.

\section{Detailed numerical experiments on a 2D butterfly domain}
\label{section:experiments_butterfly}
In this section we perform computational experiments to further explore the numerical properties of the unfitted RBF-FD method when
solving \eqref{eq:model:Poisson}
and compare them to the classical RBF-FD method in the least-squares (RBF-FD-LS) and collocation settings (RBF-FD-C).
The involved parameters are (the internodal distance of the $X$-points),
$p$ (the polynomial degree used to form the interpolant over a stencil) and $q$ (the oversampling parameter).
More precisely, we compute $h$ as the average distance between all pairs of the neighboring interpolation points:
\begin{equation*}
    h = \frac{1}{N} \sum_{i\refereeSecond{,j}=1}^N \min_{x_{\refereeSecond{j}}\neq x_i} \|x_{\refereeSecond{j}}-x_i\|_2.
\end{equation*}
The relation between the stencil size $n$ and the polynomial degree $p$ is\cite{BFFB17}:
$$n=2\binom{p+d}{d}.$$
Throughout the section we compute the relative error $\|e\|$ as:
\begin{equation*}
    \|e\| = \frac{\|u_h(Y) - u(Y)\|_{2}}{\|u(Y)\|_{2}},
\end{equation*}
where $u_h(Y)$ and $u(Y)$ are the numerical and exact solutions sampled in the $Y$-points.

All of the computations are performed in Matlab on a laptop with an Intel i7-7500U processor and 16 Gb of RAM.

\subsection{Domain $\Omega$}
The boundary of a computational domain has a butterfly-like shape (see Figure \ref{fig:intro:butterfly})
and is prescribed using a polar function:
$$r(\theta) = \frac{1}{4} \left(2 + \sin\left(2t\right) - 0.01\cos\left(5t-\frac{\pi}{2}\right) + 0.63\sin\left(6t - 0.1\right)\right),$$
where $r$ is a radial coordinate and $\theta \in [0, 2\pi]$ the angle.
\begin{figure}[h!]
    \centering
    \begin{tabular}{cc}
        $u_1$                                                                                     & $u_2$ \\
        \includegraphics[width=0.46\linewidth]{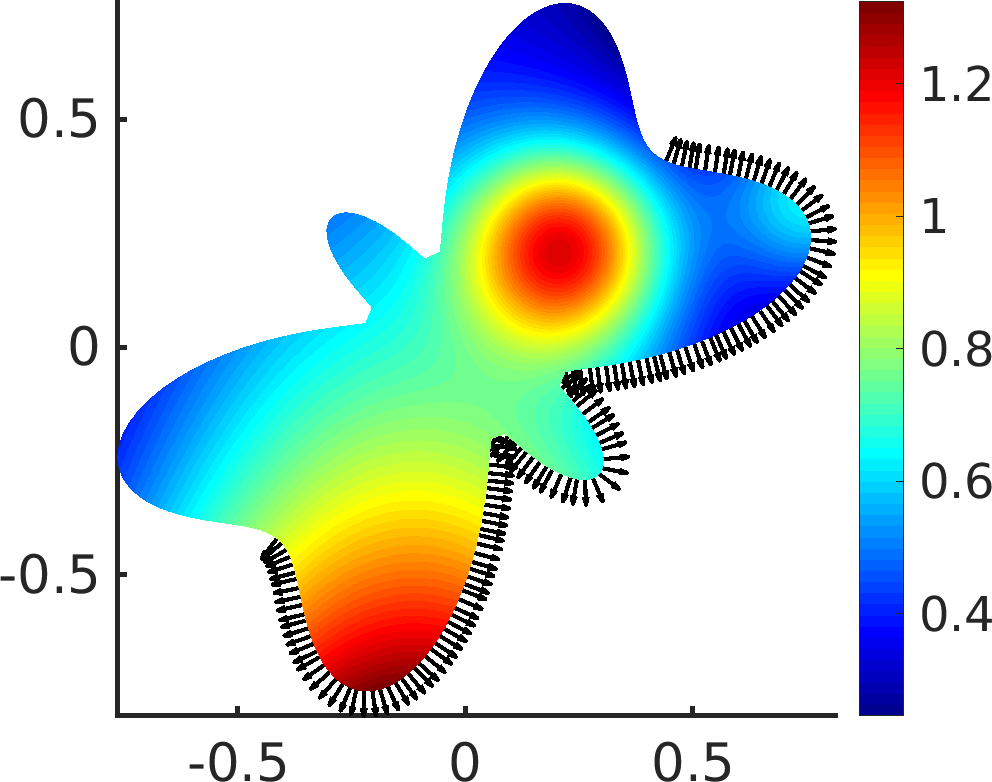} &
        \includegraphics[width=0.46\linewidth]{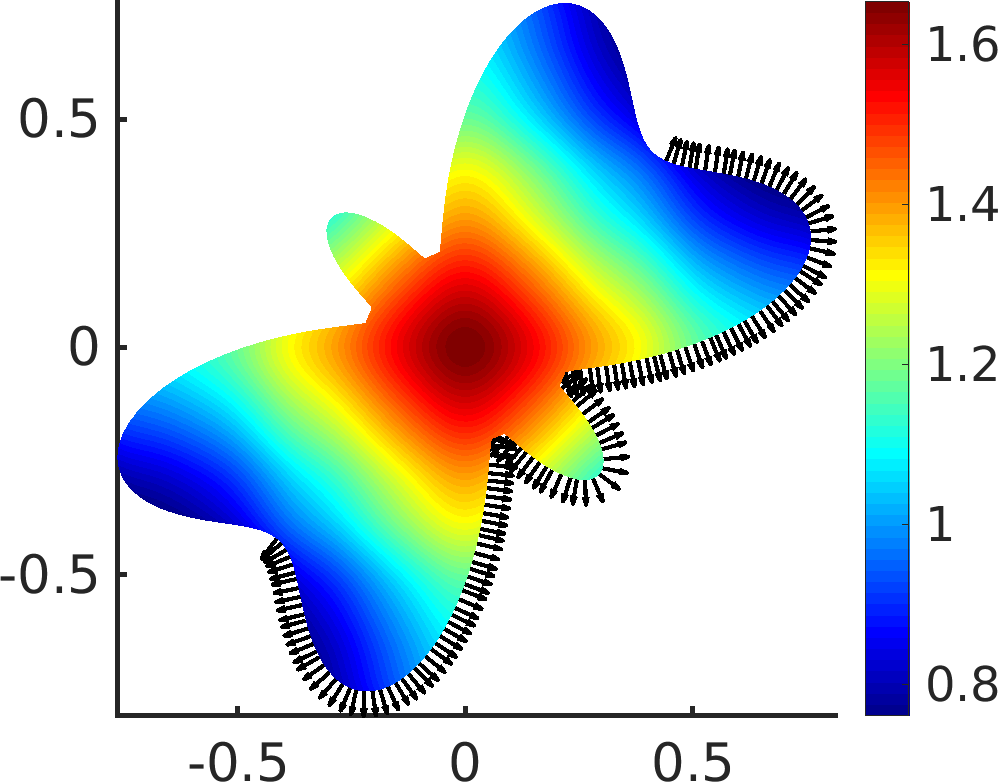}    \\
    \end{tabular}
    \caption{Solution functions on a butterfly domain. Franke function over a butterfly ($u_1$), and the \refereeSecond{Truncated non-analytic} function over a butterfly ($u_2$).  The black outward normals indicate
        the locations of the Neumann condition. The Dirichlet condition is enforced at locations where there are no normals displayed over the boundary.}
    \label{fig:experiments:solutionFunctions}
\end{figure}

\refereeSecond{
    \subsection{Point sets}
    \label{sec:experiments_butterfly:pointsets}
    When using the unfitted RBF-FD-LS method, the point sets $X$ and $Y$ are computed according to the description in Section \ref{section:method}, where for constructing
    $Y$ we use the oversampling parameter $q=5$.
    For a more detailed study of the
    oversampling parameter and its impact on the approximation error and the stability norm, we refer the reader to \cite{tominec2020squares}.
    When using the RBF-FD-LS method the point set $X$ is computed using the DistMesh algorithm \cite{distmesh}, while the point set $Y$ is computed in the same way as for the unfitted
    RBF-FD-LS method. For the RBF-FD-C method, the point set $X$ is also computed using the DistMesh algorithm, note that in the collocation case we have $Y=X$.
}

\subsection{Solution functions}
We pick two solution functions to compute the right-hand-sides of \eqref{eq:model:Poisson}. Those are:
\begin{eqnarray}
    \label{eq:experiments:franke}
    u_1(x,y) &=& \frac{3}{4} e^{-\frac{1}{4}((9x-2)^2 + (9y-2)^2)} + \frac{3}{4} e^{-(\frac{1}{49}(9x+1)^2 + \frac{1}{10}(9y+1)^2)} \\
    &...& \refereeSecond{+ \frac{1}{2} e^{-\frac{1}{4}((9x-7)^2 + (9y-3)^2)} - \frac{1}{5} e^{-((9x-4)^2 + (9y-7)^2)}}, \nonumber \\
    \label{eq:experiments:nonanalytic}
    u_2(x,y) &=& \sum_{k=0}^5 e^{-\sqrt{2^k}}\left( \cos(2^k x) + \cos(2^k y)\right),
\end{eqnarray}
where $u_1$ is the Franke function, a commonly used infinitely smooth test function for benchmarking multivariate approximations.
Function $u_2$ is a truncated series of an infinitely smooth function that is at the same time not analytic: we refer to $u_2$ as the \refereeSecond{Truncated non-analytic} function.
Both functions over the butterfly domain are displayed in Figure \ref{fig:experiments:solutionFunctions}.

\subsection{Convergence under node refinement with Dirichlet condition}
\begin{figure}[h!]
    \centering
    \begin{tabular}{ccc}
        \multicolumn{3}{c}{\textbf{Franke: Dirichlet on $\partial \Omega$}}                                                                                       \\
        \hspace{0.7cm} $p=2$                                                                                        & \hspace{0.7cm} $p=4$ & \hspace{0.7cm} $p=6$ \\
        \includegraphics[width=0.3\linewidth]{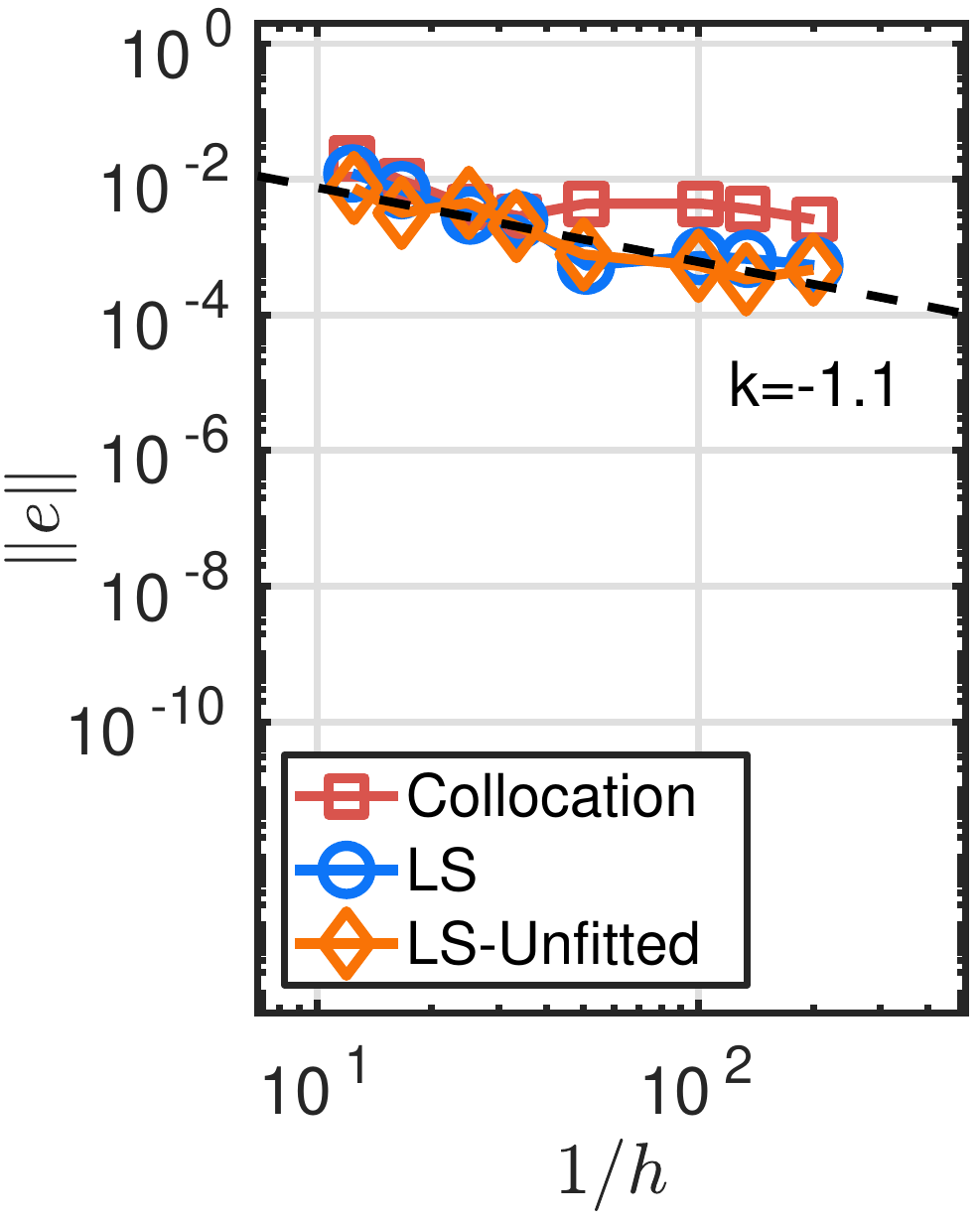} &
        \includegraphics[width=0.3\linewidth]{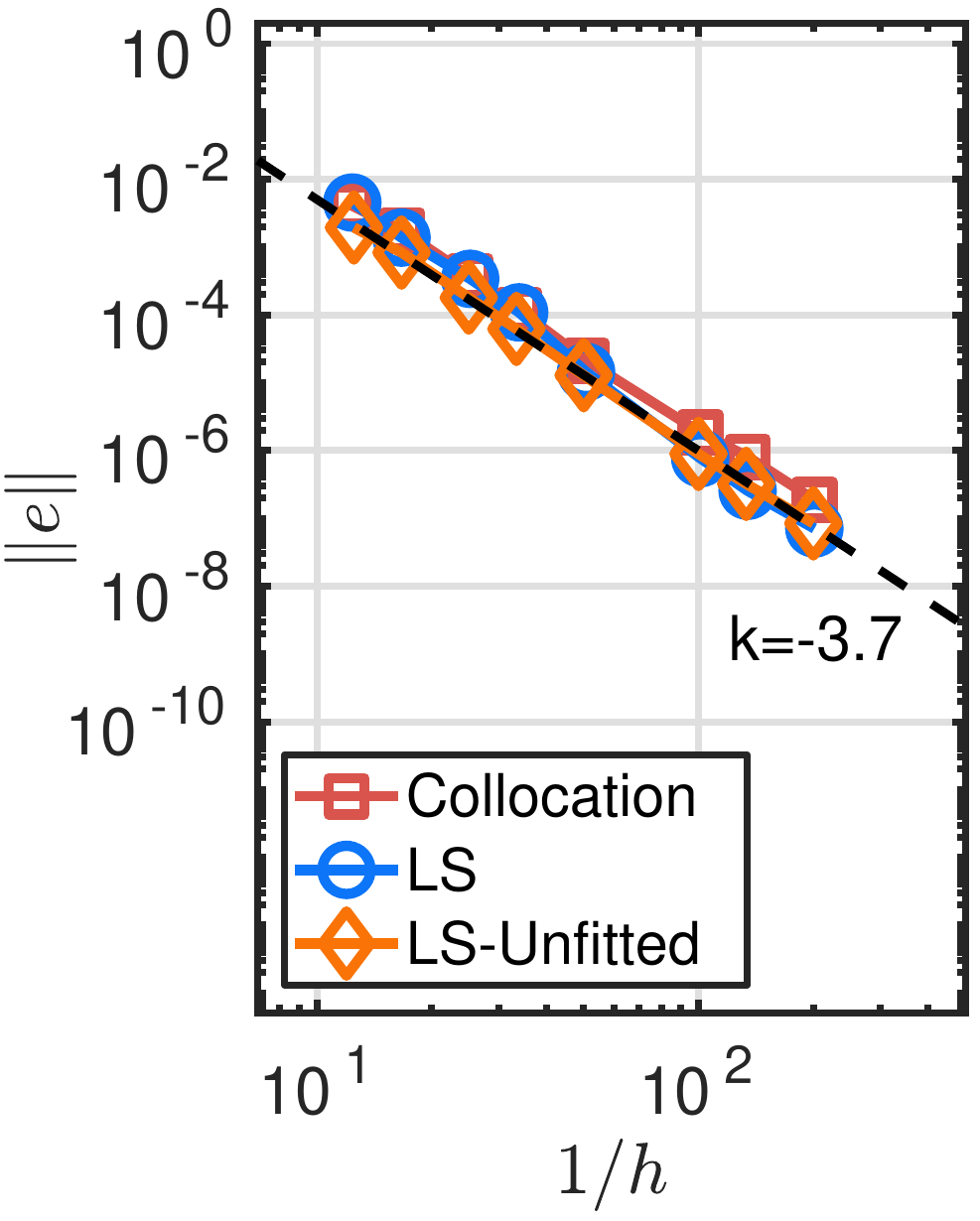} &
        \includegraphics[width=0.3\linewidth]{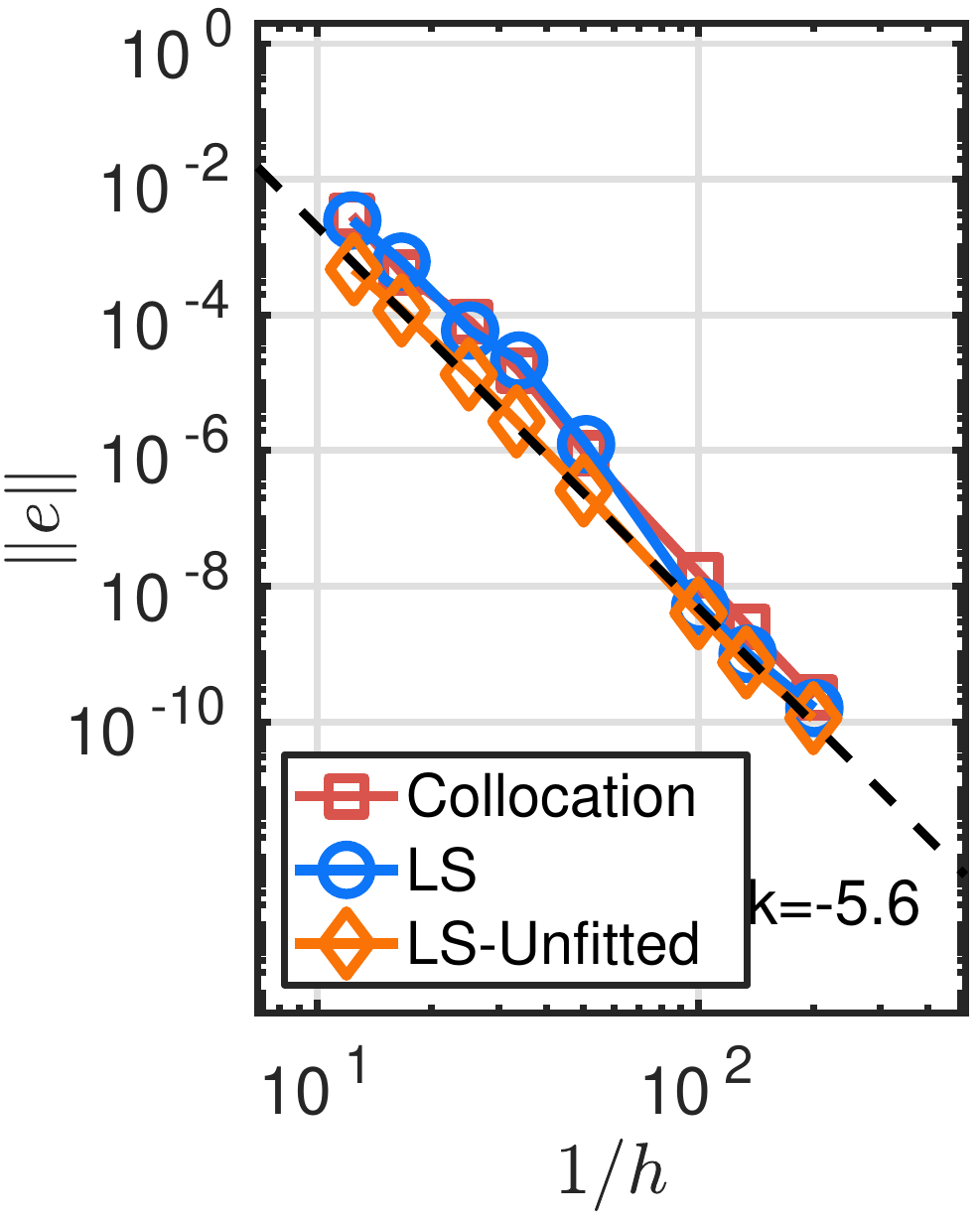}
    \end{tabular}
    \caption{Error as a function of the inverse internodal distance $1/h$ for different polynomial degrees $p$.
        In this case the Franke function is used to manufacture the right-hand-sides of the Poisson equation with the Dirichlet boundary condition.}
    \label{fig:experiments:butterfly:href:franke_dir}
\end{figure}

\begin{figure}[h!]
    \centering
    \begin{tabular}{ccc}
        \multicolumn{3}{c}{\textbf{\refereeSecond{Truncated non-analytic}: Dirichlet on $\partial \Omega$}}                                                                                      \\
        \hspace{0.7cm} $p=2$                                                                                             & \hspace{0.7cm} $p=4$ & \hspace{0.7cm} $p=6$ \\
        \includegraphics[width=0.3\linewidth]{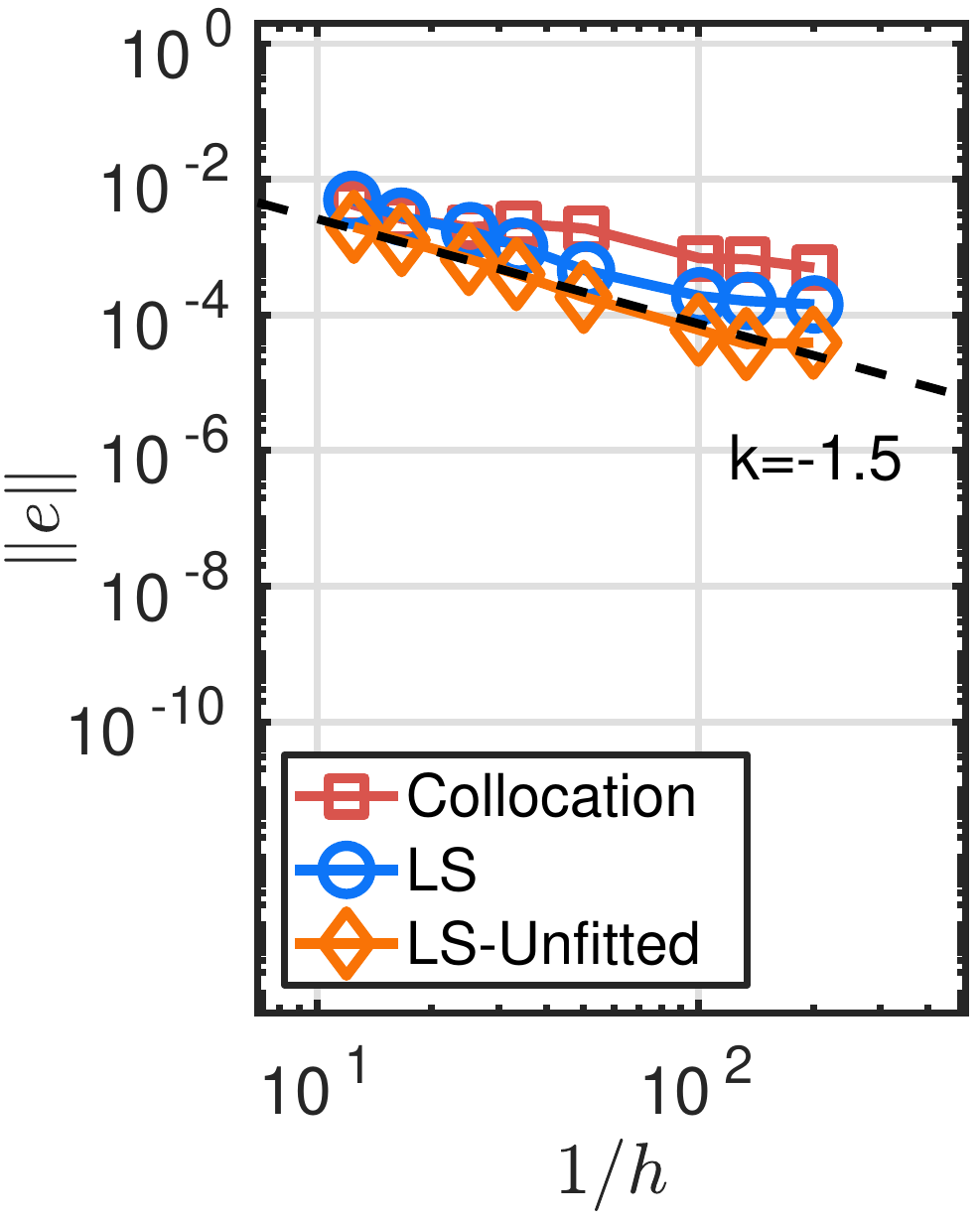} &
        \includegraphics[width=0.3\linewidth]{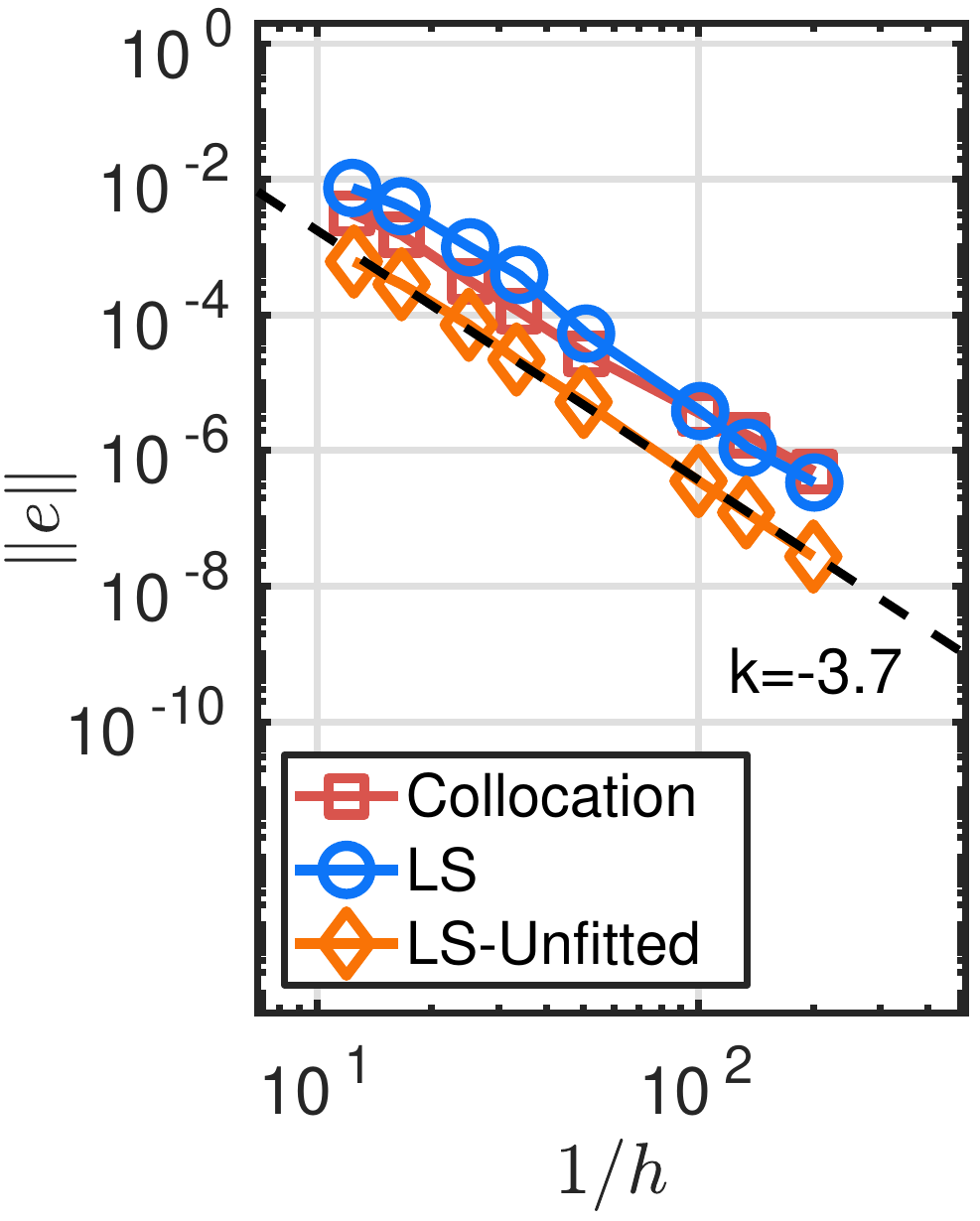} &
        \includegraphics[width=0.3\linewidth]{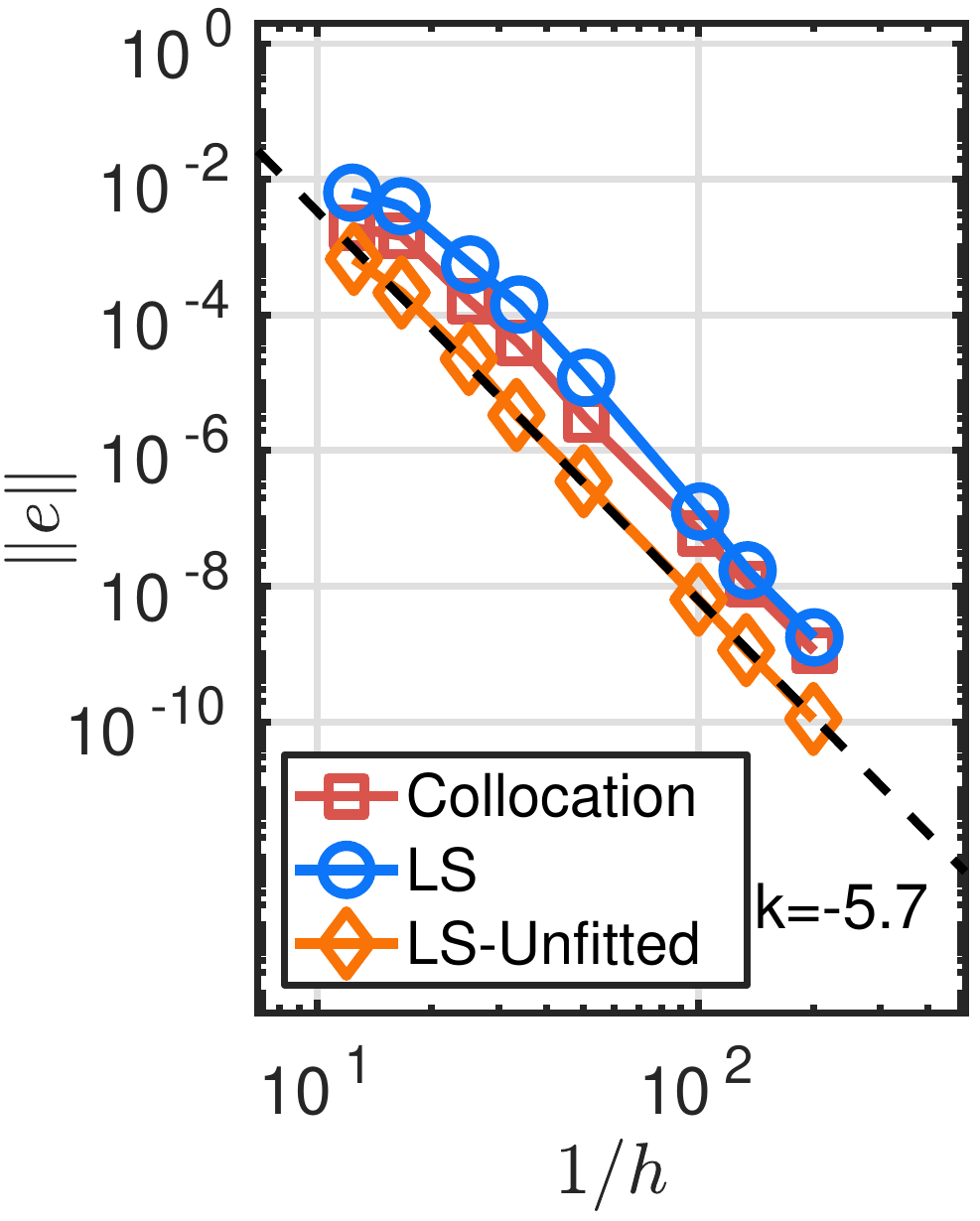}
    \end{tabular}
    \caption{Error as a function of the inverse internodal distance $1/h$ for different polynomial degrees $p$.
        In this case the \refereeSecond{Truncated non-analytic} function is used to manufacture the right-hand-sides of the Poisson equation with the Dirichlet boundary condition.}
    \label{fig:experiments:butterfly:href:nonanalytic_dir}
\end{figure}

First we use the Dirichlet boundary condition on all of $\partial\Omega$ and
compute the error as a function of $1/h$ for three choices of polynomial degrees: $p=2$, $p=4$ and $p=6$.
The results for the Franke and \refereeSecond{Truncated non-analytic} functions are given in Figure \ref{fig:experiments:butterfly:href:franke_dir} and
Figure \ref{fig:experiments:butterfly:href:nonanalytic_dir}, respectively.
We observe that for the Franke function,
the error behavior is similar for all three methods, while for the \refereeSecond{Truncated non-analytic} function the error is smaller for the unfitted RBF-FD method when
the stencil sizes are larger ($p=4$ and $p=6$). We note that the Dirichlet condition is enforced exactly in
RBF-FD-C and RBF-FD-LS, but weakly in the unfitted RBF-FD-LS method. This could lead a reader to intuitively presume that the overall error
could behave in favor of RBF-FD-C and RBF-FD-LS.

The smaller error in the unfitted variant can be attributed to the smaller skeweness of
the stencils which are placed in the interior, but still touch the boundary of $\Omega$. \refereeSecond{Supporting discussion is available in Appendix \ref{sec:appendix:skeweness}.}

\subsection{Convergence under node refinement with Dirichlet and Neumann conditions}
\label{section:experiments_butterfly:mixedBCs}
\begin{figure}[h!]
    \centering
    \begin{tabular}{ccc}
        \multicolumn{3}{c}{\textbf{Franke: Dirichlet on $\partial \Omega_0$, Neumann on $\partial \Omega_1$}}                                       \\
        \hspace{0.7cm} $p=2$                                                                          & \hspace{0.7cm} $p=4$ & \hspace{0.7cm} $p=6$ \\
        \includegraphics[width=0.3\linewidth]{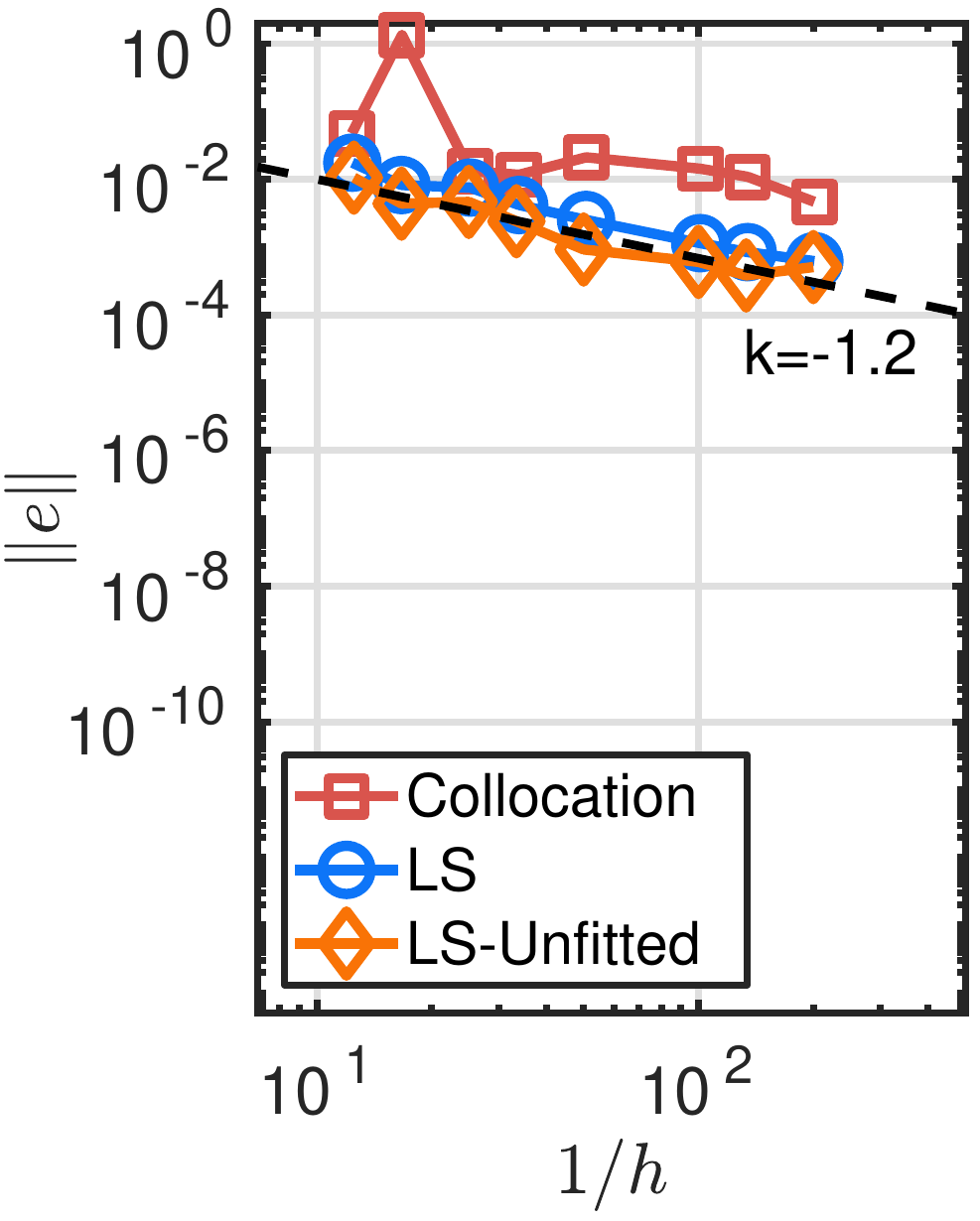} &
        \includegraphics[width=0.3\linewidth]{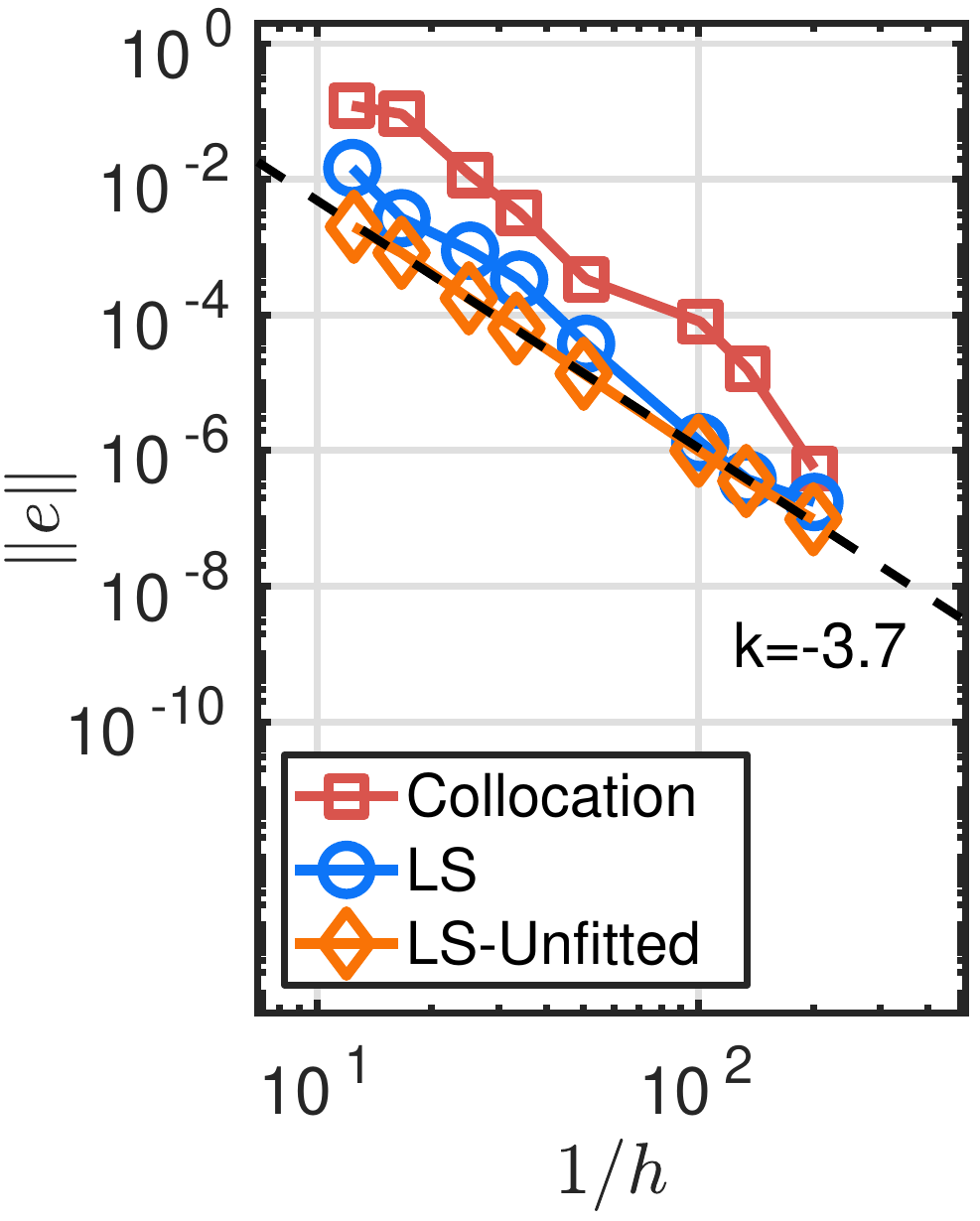} &
        \includegraphics[width=0.3\linewidth]{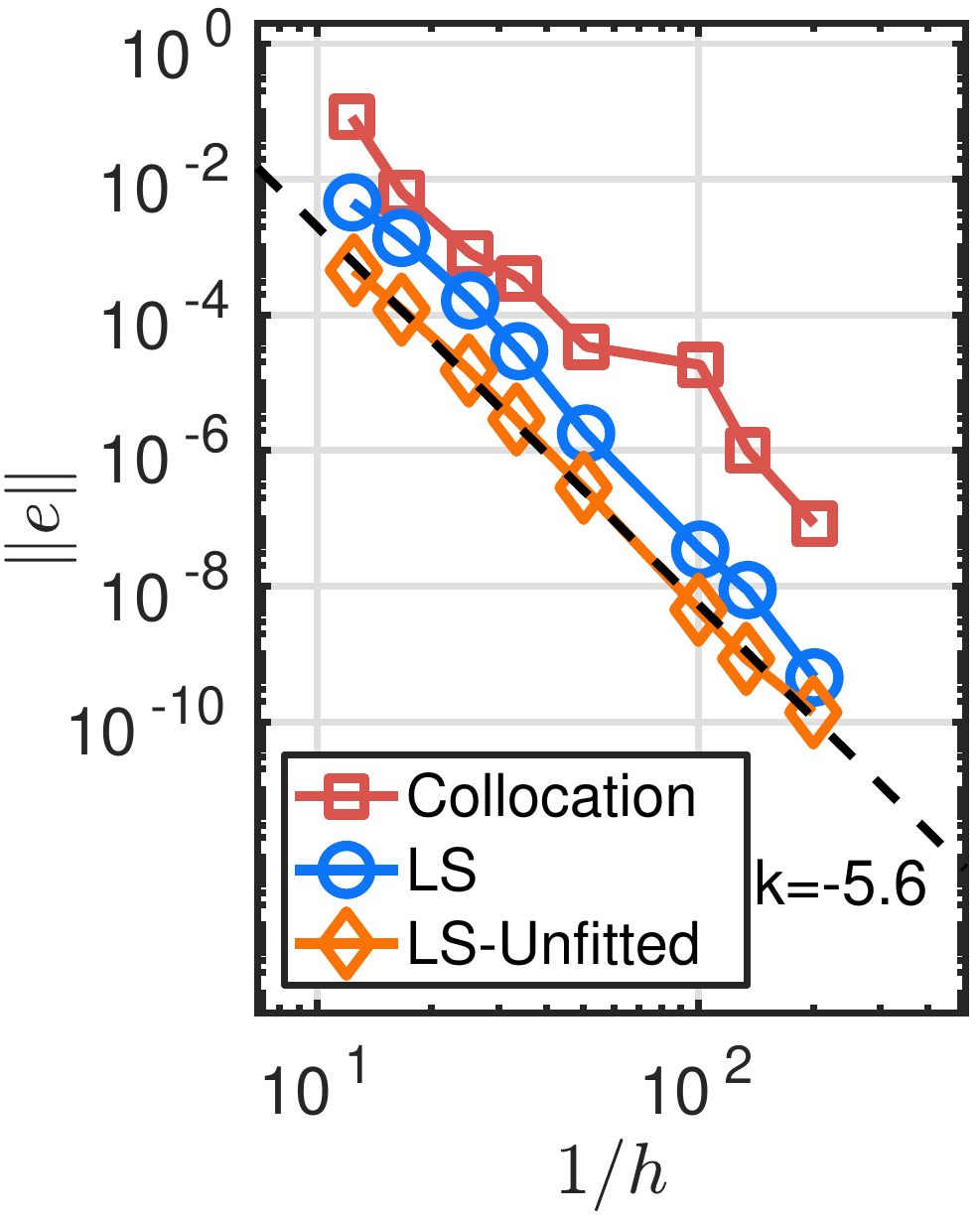}
    \end{tabular}
    \caption{Error as a function of the internodal distance $h$ for different polynomial degrees $p$.
        In this case the Franke function is used to manufacture the right-hand-sides of the Poisson equation with mixed boundary conditions.}
    \label{fig:experiments:butterfly:href:franke_dirneu}
\end{figure}

\begin{figure}[h!]
    \centering
    \begin{tabular}{ccc}
        \multicolumn{3}{c}{\textbf{\refereeSecond{Truncated non-analytic}: Dirichlet on $\partial \Omega_0$, Neumann on $\partial \Omega_1$}}                                      \\
        \hspace{0.7cm} $p=2$                                                                               & \hspace{0.7cm} $p=4$ & \hspace{0.7cm} $p=6$ \\
        \includegraphics[width=0.3\linewidth]{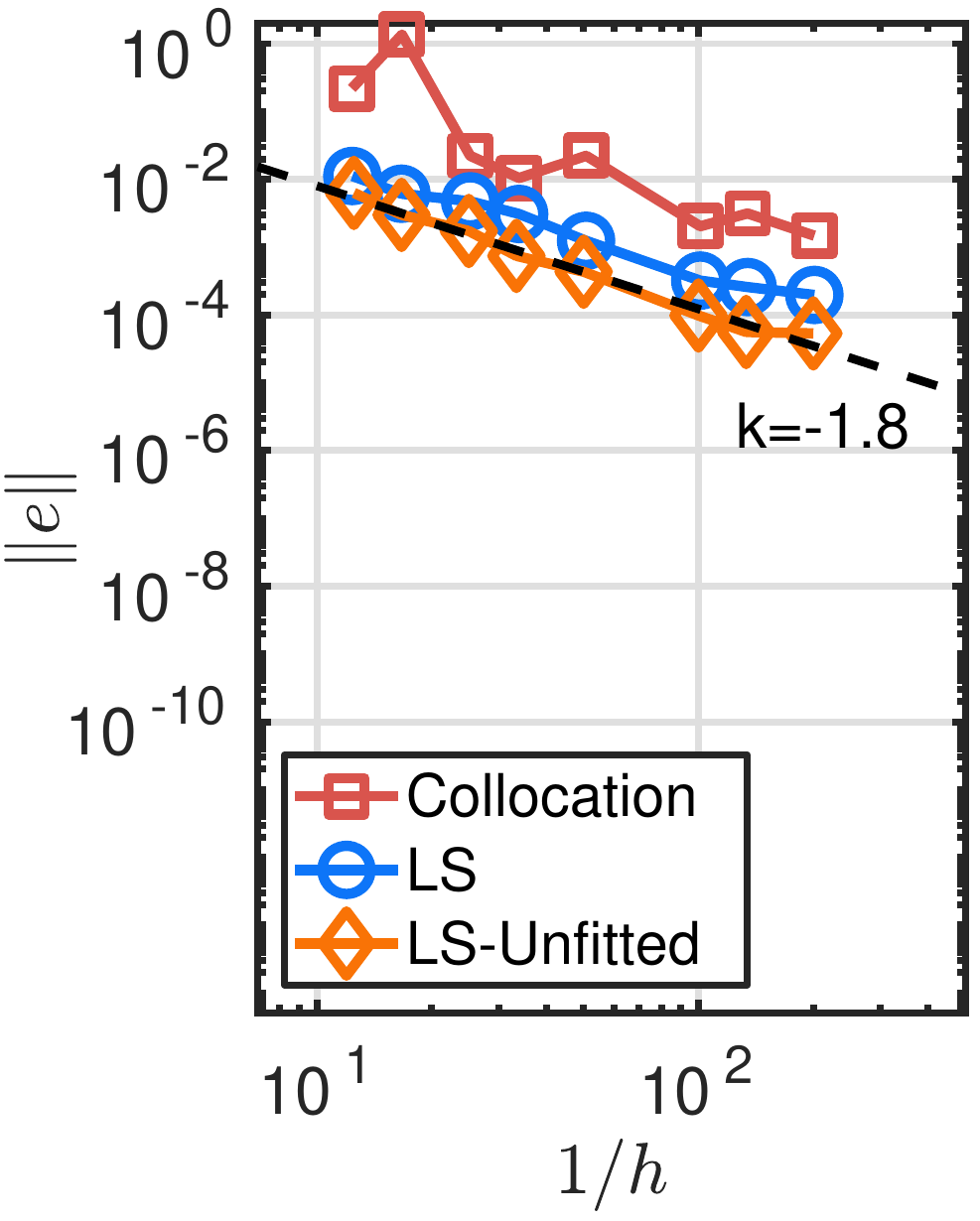} &
        \includegraphics[width=0.3\linewidth]{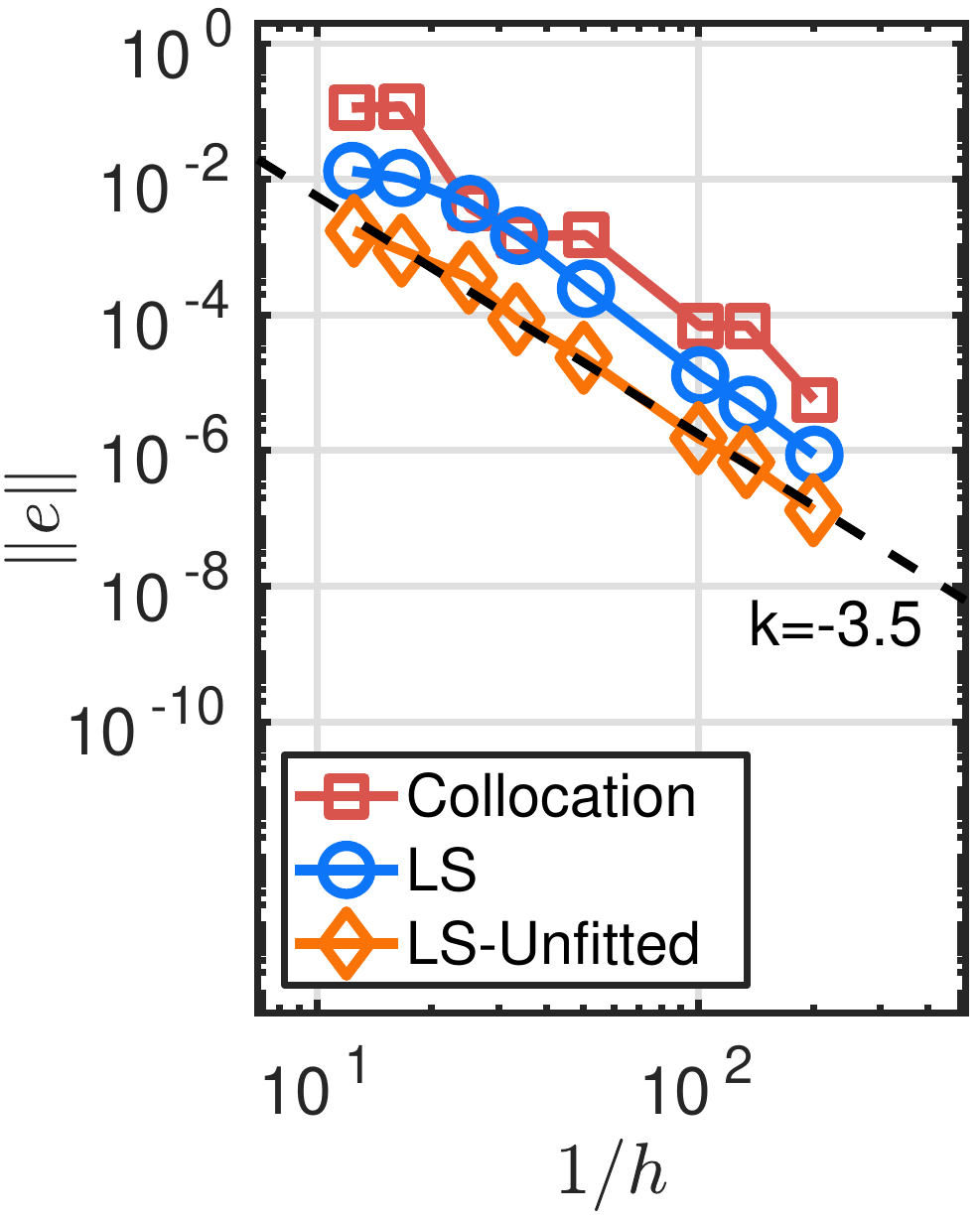} &
        \includegraphics[width=0.3\linewidth]{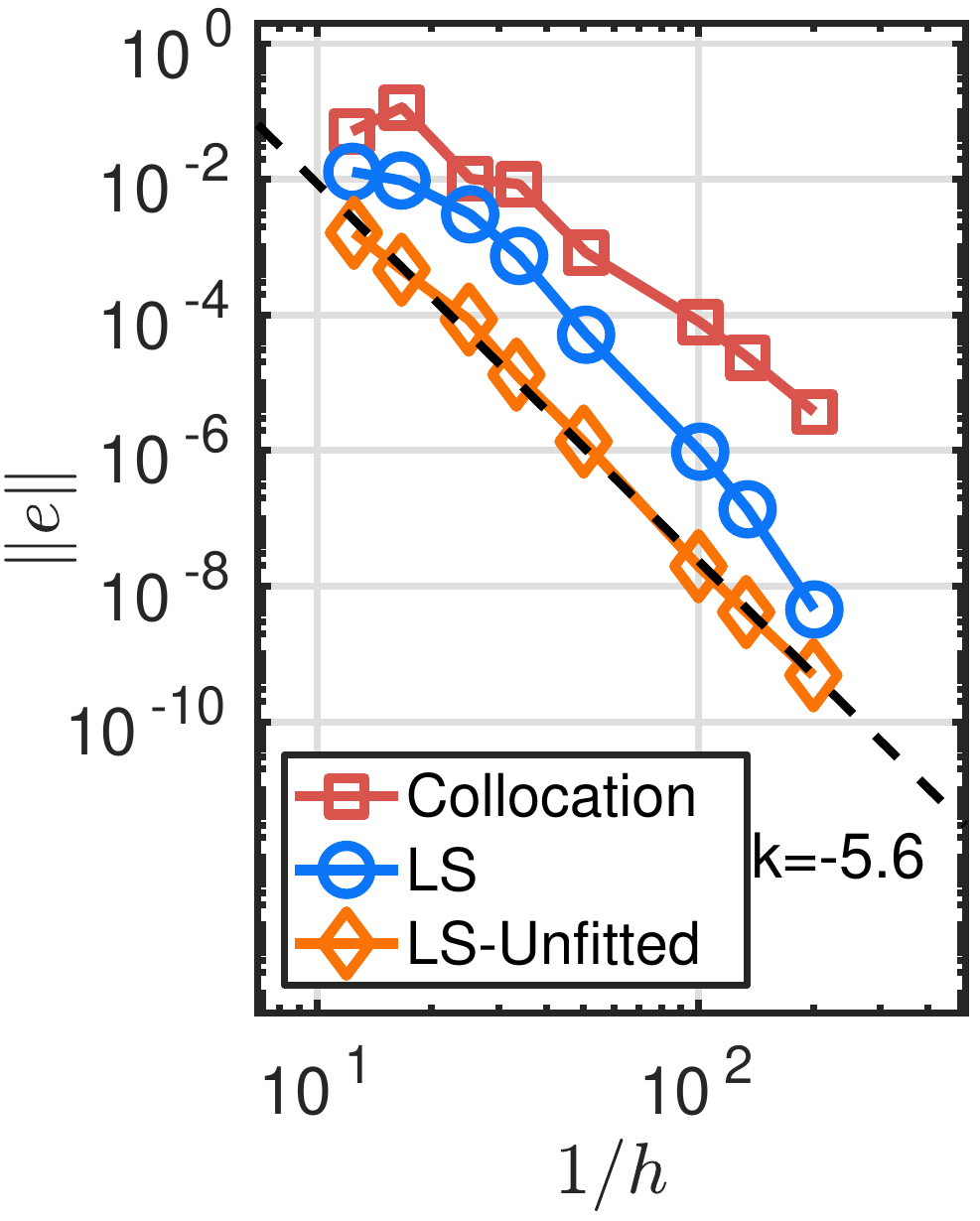}
    \end{tabular}
    \caption{Error as a function of the internodal distance $h$ for different polynomial degrees $p$.
        In this case the \refereeSecond{Truncated non-analytic} function is used to manufacture the right-hand-sides of the Poisson equation with mixed boundary conditions.}
    \label{fig:experiments:butterfly:href:nonanalytic_dirneu}
\end{figure}

Here we follow the formulation from \eqref{eq:model:Poisson} in the sense that the Dirichlet and Neumann boundary conditions
are used on two disjoint parts of the domain. The error is again computed as a function of $1/h$ for $p=2$, $p=4$ and $p=6$.
The results for the Franke and the \refereeSecond{Truncated non-analytic} functions are given in Figure \ref{fig:experiments:butterfly:href:franke_dirneu} and
Figure \ref{fig:experiments:butterfly:href:nonanalytic_dirneu} respectively.

We observe that the collocation setting
has larger errors than the least-squares counterparts, which is expected \cite{tominec2020squares}.
For all $p$, the unfitted RBF-FD-LS method has a slightly smaller error compared to RBF-FD-LS in the case of the Franke function.
The difference between the errors is more pronounced in the \refereeSecond{Truncated non-analytic} case, especially when $p=4$ and $p=6$, where the error of the unfitted
RBF-FD-LS is $10$-times smaller throughout the refinement.

We know that the Franke function does not oscillate around the boundary and that the
\refereeSecond{Truncated non-analytic} function is highly oscillatory on a fine scale over the whole $\Omega$.
In the first case the less skewed stencils are then not expected to have a significant impact, while in the second case the
less skewed stencils are expected to significantly contribute towards a smaller PDE error. This can be accounted to
smaller Lebesgue constants which enable a better approximation error \refereeSecond{see Appendix \ref{sec:appendix:skeweness} for a more detailed explanation.}

\subsection{Error as a function of runtime \refereeSecond{with Dirichlet and Neumann conditions}}
In the subsections above we confirmed that the approximation error is smaller for the unfitted RBF-FD-LS method compared with
RBF-FD-LS and RBF-FD-C\refereeSecond{, for all the cases that were considered.
The computational time for both RBF-FD-LS methods is expected to be slightly larger than for RBF-FD-C. In the RBF-FD-LS case, 
the system \eqref{eq:method:PDE_discrete} is rectangular and requires 
a sparse QR decomposition as the most expensive intermediate step, so that the solution is obtained. In the RBF-FD-C case, 
the system \eqref{eq:method:PDE_discrete} is square and requires a sparse LU decomposition. The sparse QR decomposition is 
slightly more expensive compared to the sparse LU decomposition.
In addition, the computational cost of the unfitted RBF-FD-LS method is expected to be slightly larger than for RBF-FD-LS and RBF-FD-C due to the 
additional degrees of freedom that extend over the boundary of $\Omega$ and give rise to an increased number of columns
in the matrices $D_h$ and $E_h$ from \eqref{eq:method:PDE_solution}.}
The question that we address in this subsection is whether the error for the unfitted RBF-FD-LS method is small enough to compensate for the larger computational cost. 
\refereeSecond{For this reason we measure the code runtime vs. the obtained accuracy. 
The runtime is a sum of the following execution times.
\begin{itemize}
\item The closest neighbor search for obtaining the stencil neighbor points around every $x_i \in X$.
\item Forming and inverting a set of local interpolation matrices \eqref{eq:M}.
\item Forming the evaluation and differentiation weights used in \eqref{eq:method:stencil_cardinal}.
\item Assembly of the PDE matrix used in \eqref{eq:method:PDE_discrete}.
\item Solving the linear system \eqref{eq:method:PDE_discrete} using \emph{mldivide()} in Matlab.
\end{itemize}
}

\begin{figure}[h!]
    \centering
    \begin{tabular}{ccc}
        \hspace{0.7cm} $p=2$                                                                                   & \hspace{0.7cm} $p=4$ & \hspace{0.7cm} $p=6$ \\
        \includegraphics[width=0.3\linewidth]{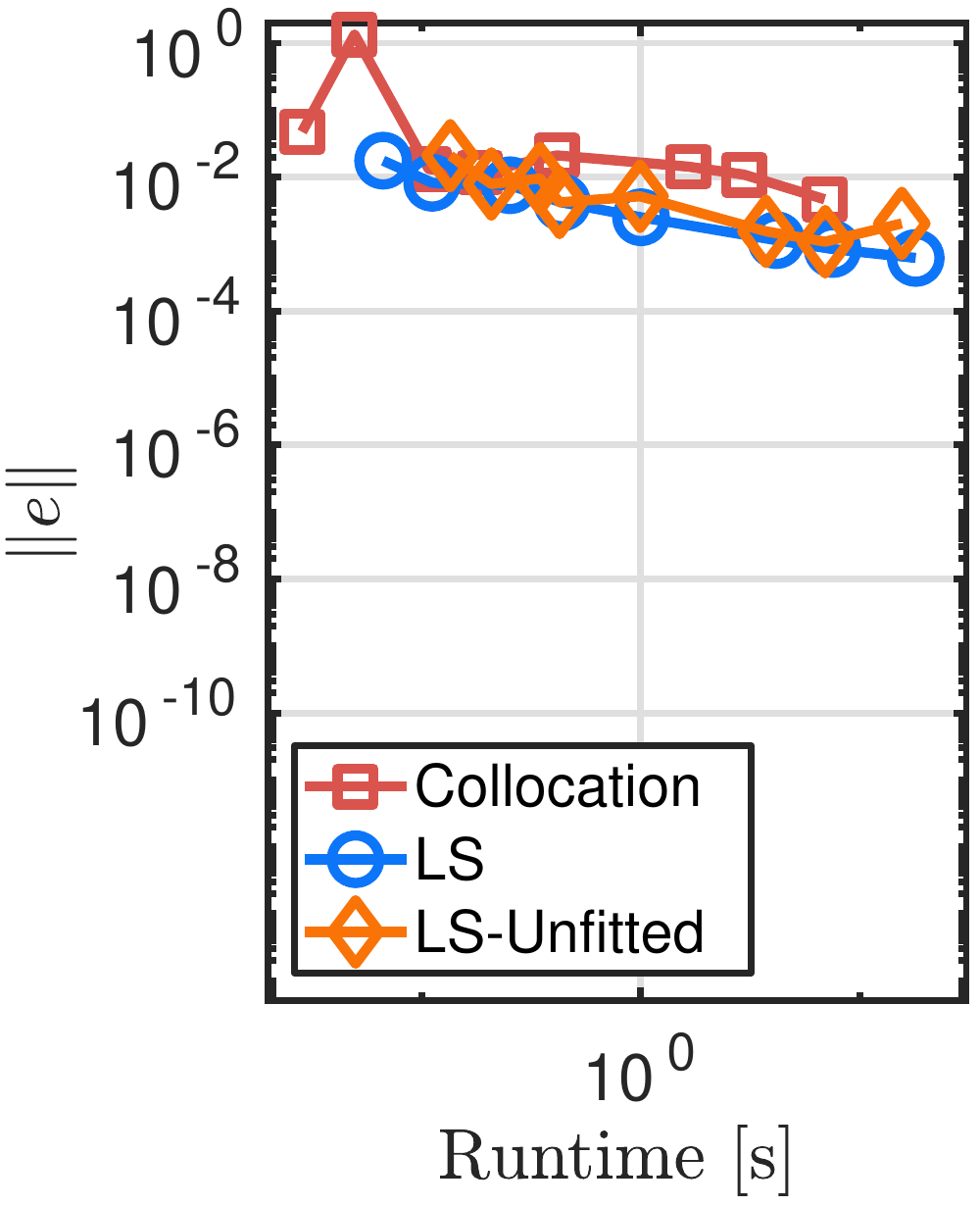} &
        \includegraphics[width=0.3\linewidth]{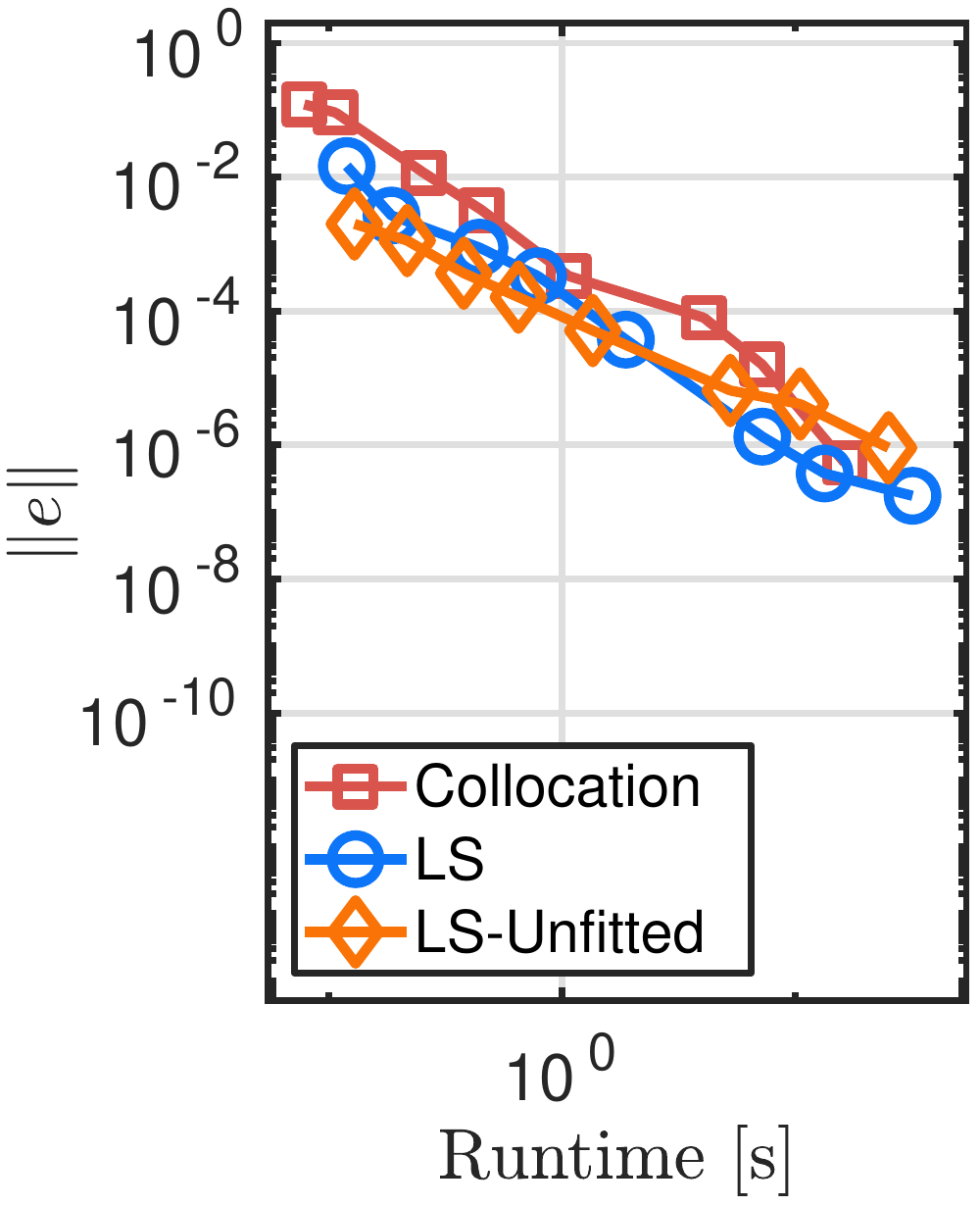} &
        \includegraphics[width=0.3\linewidth]{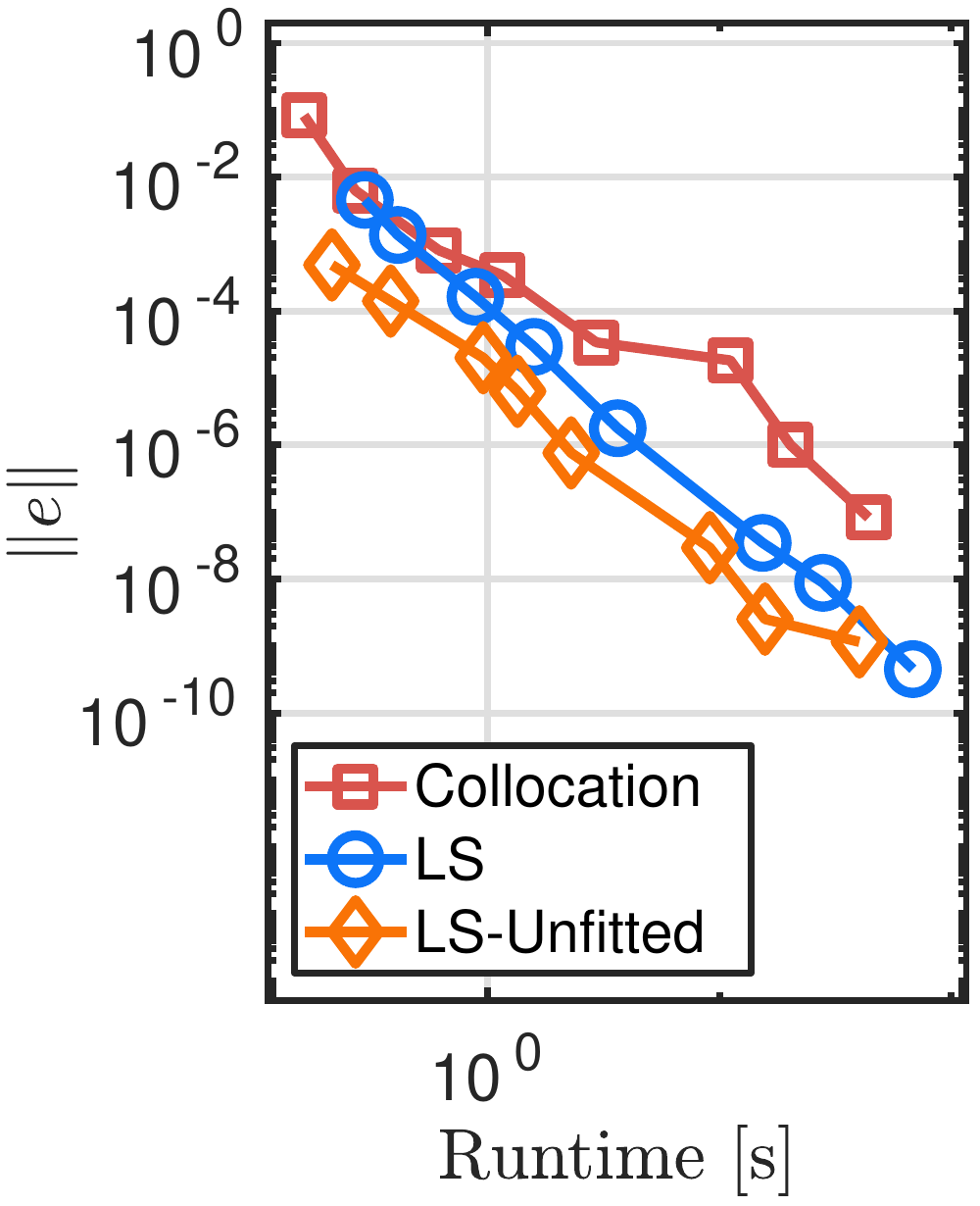}
    \end{tabular}
    \caption{Error versus runtime for three polynomial degrees $p$ when the Franke function is used as the exact solution.}
    \label{fig:experiments:butterfly:href:franke_errVsRuntime}
\end{figure}

\begin{figure}[h!]
    \centering
    \begin{tabular}{ccc}
        \hspace{0.7cm} $p=2$                                                                                        & \hspace{0.7cm} $p=4$ & \hspace{0.7cm} $p=6$ \\
        \includegraphics[width=0.3\linewidth]{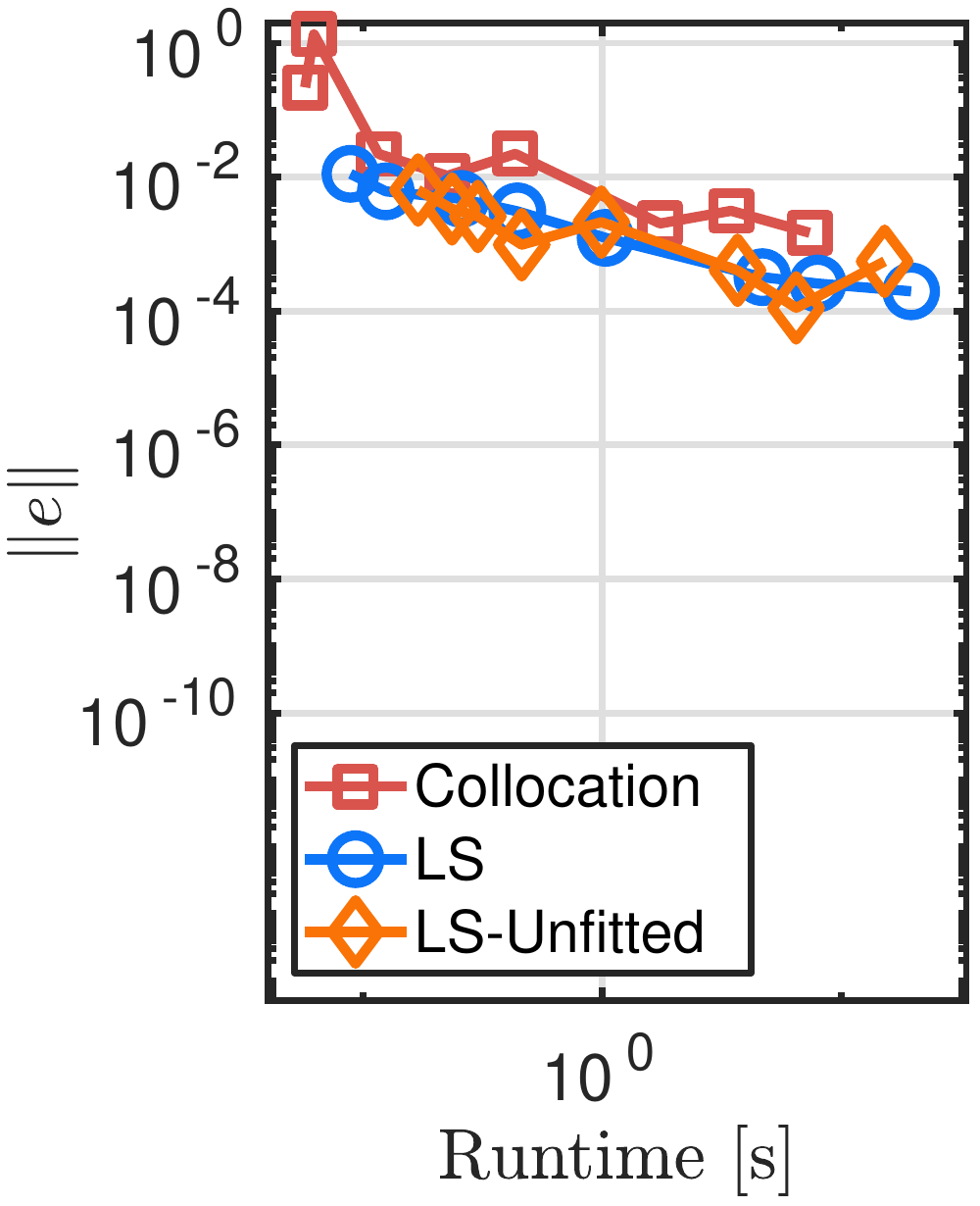} &
        \includegraphics[width=0.3\linewidth]{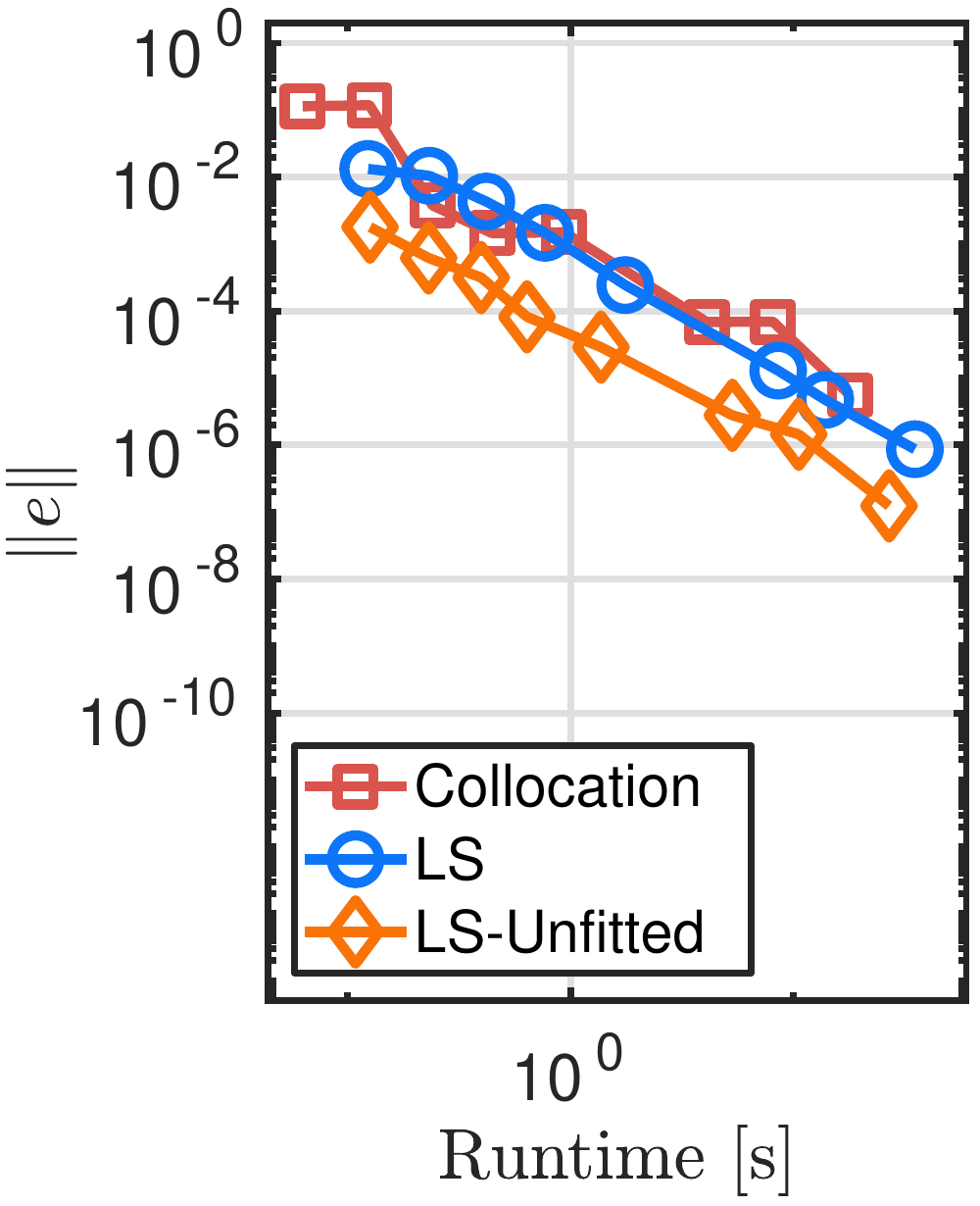} &
        \includegraphics[width=0.3\linewidth]{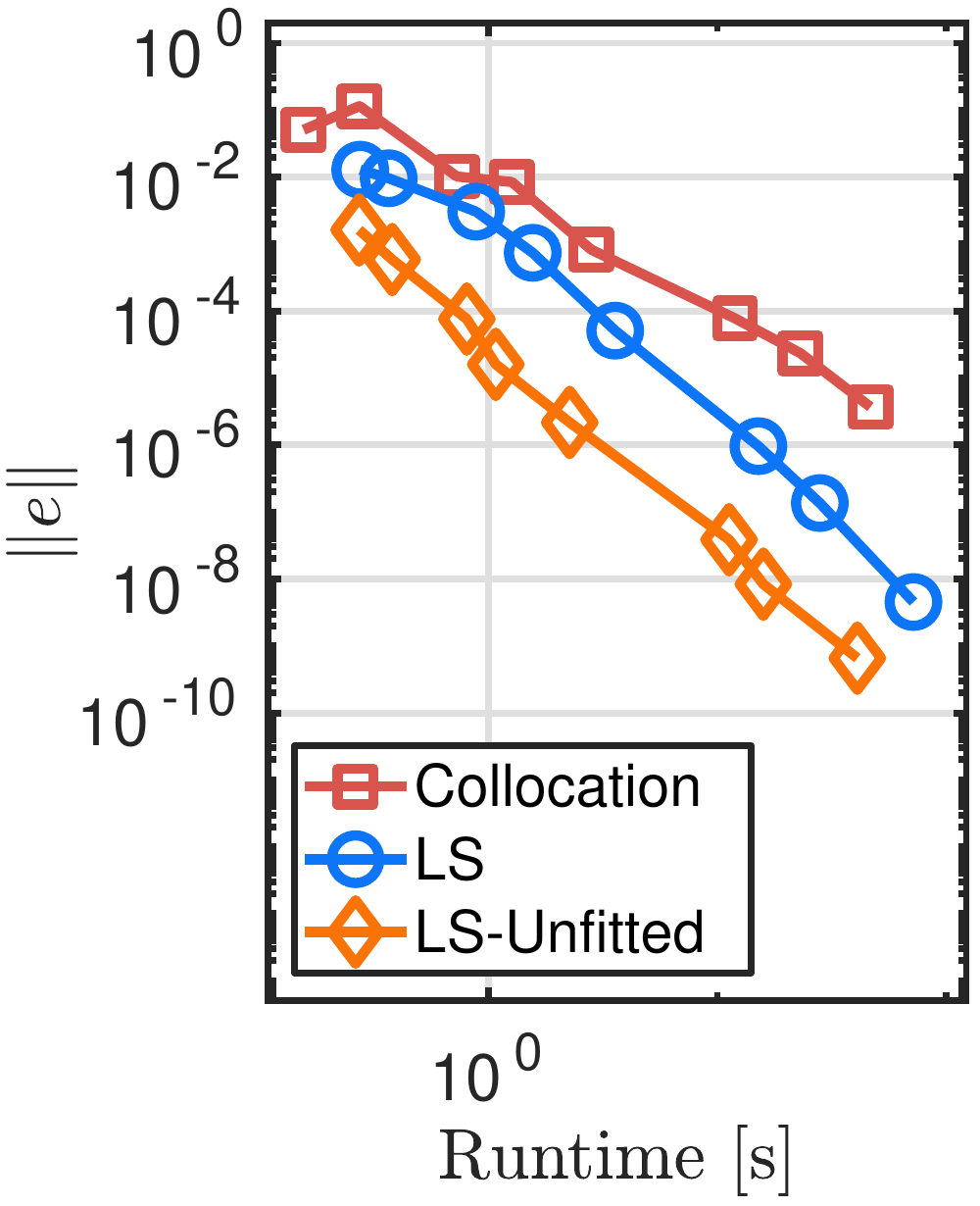}
    \end{tabular}
    \caption{Error versus runtime for three polynomial degrees $p$ when the \refereeSecond{Truncated non-analytic} function is used as the exact solution.}
    \label{fig:experiments:butterfly:href:nonanalytic_errVsRuntime}
\end{figure}

In Figure \ref{fig:experiments:butterfly:href:franke_errVsRuntime}
we can observe that in the
Franke case error vs. runtime is comparable for all three settings when $p=2$ and $p=4$. When $p=6$ the ratio
is in favor of the unfitted RBF-FD method. In Figure
\ref{fig:experiments:butterfly:href:nonanalytic_errVsRuntime} we see that in the \refereeSecond{Truncated non-analytic} case the ratio is in favor
of the unfitted RBF-FD method for $p=4$ and $p=6$, while it is tied between the three settings for $p=2$.

The conclusion is that the unfitted RBF-FD method does not lag behind in error vs. runtime compared to the two fitted settings, despite
a larger amount of degrees of freedom.
Moreover, in cases when stencils are large and the approximated solution oscillatory around the boundaries (the \refereeSecond{Truncated non-analytic} case) it performs significantly better.

\subsection{Stability norms and condition numbers under node refinement}
Here we measure the condition number of the PDE matrix $D_h$ given in
\eqref{eq:method:PDE_discrete}. 
In addition we also measure the stability norm $||E_h||_2\,||D_h^+||_2$ which can be understood as the well-posedness constant that multiplies
the consistency term in the error estimate \eqref{eq:theory:finalestimate}. 
\refereeSecond{In both cases, we impose Neumann and Dirichlet boundary conditions in the same way as in 
Section \ref{section:experiments_butterfly:mixedBCs}.}
The condition number of (for example) a matrix $D_h$ is computed by:
$$\kappa\left(D_h\right) = ||D_h||_2\,||D_h^+||_2 = \frac{\sigma_{\max}\left(D_h\right)}{\sigma_{\min}\left(D_h\right)},$$
and the stability norm involved in the error estimate \eqref{eq:theory:finalestimate} by:
$$\frac{1}{\sqrt{M}}||E_h||_2\, ||D_h^+||_2 = \frac{1}{\sqrt{M}} \frac{\sigma_{\max}\left(E_h\right)}{\sigma_{\min}\left(D_h \right)}.$$
\refereeSecond{Note that when using the two relations stated above in the RBF-FD-C case, the evaluation and interpolation point sets are the same ($Y = X$), which implies $M=N$. 
Furthermore, in the RBF-FD-C case $E_h$ is a square identity matrix.}
\begin{figure}[h!]
    \centering
    \begin{tabular}{ccc}
        \hspace{0.7cm} $p=2$                                                                           & \hspace{0.7cm} $p=4$ & \hspace{0.7cm} $p=6$ \\
        \includegraphics[width=0.3\linewidth]{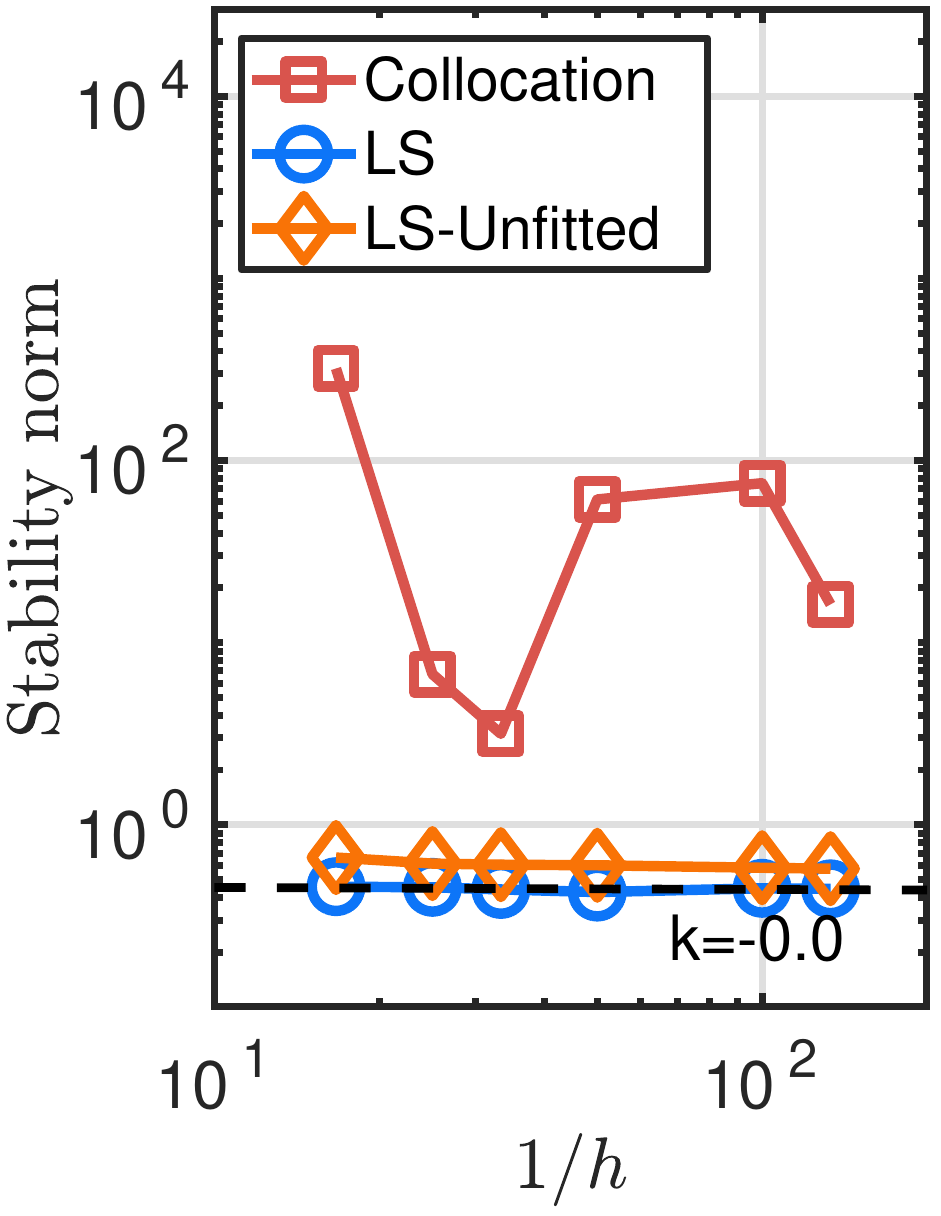} &
        \includegraphics[width=0.3\linewidth]{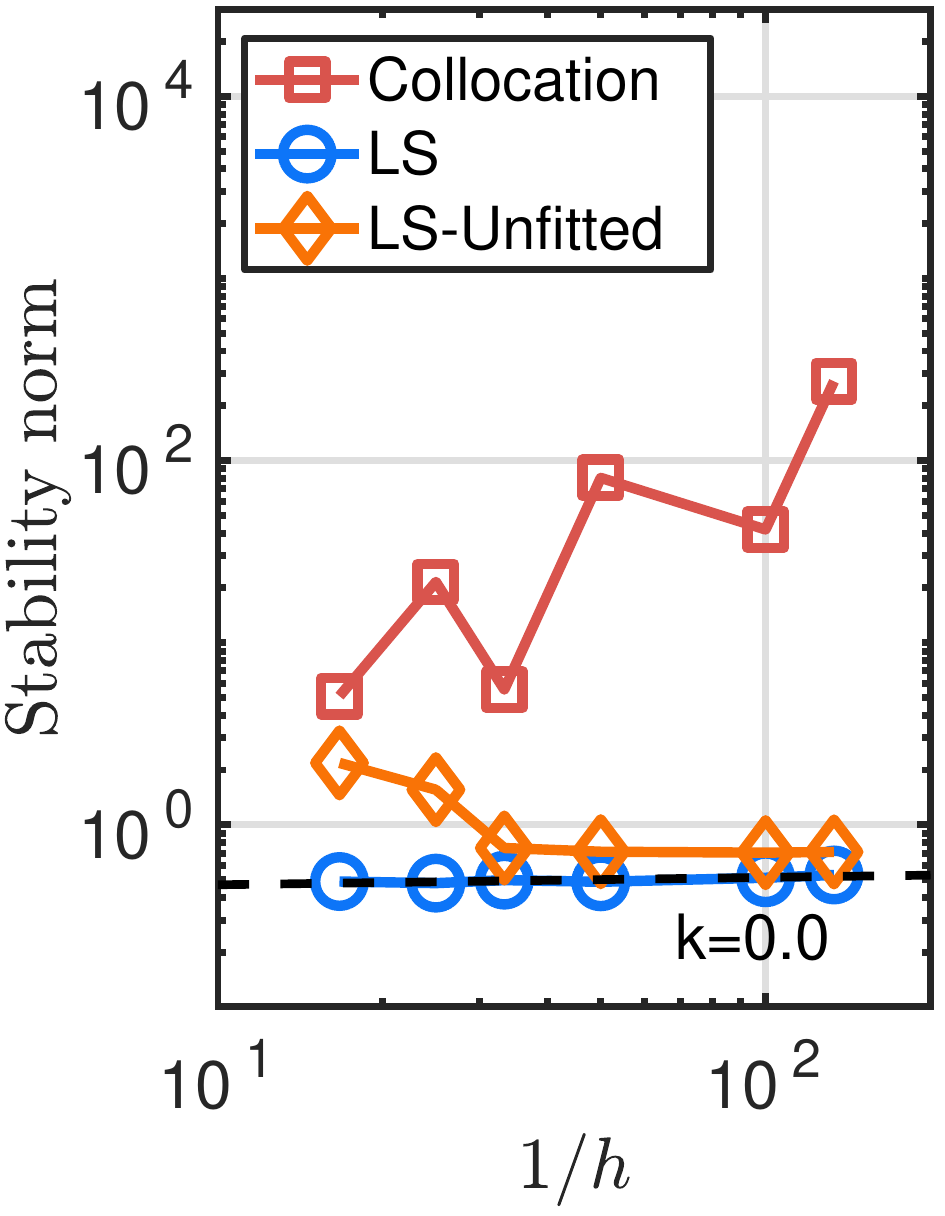} &
        \includegraphics[width=0.3\linewidth]{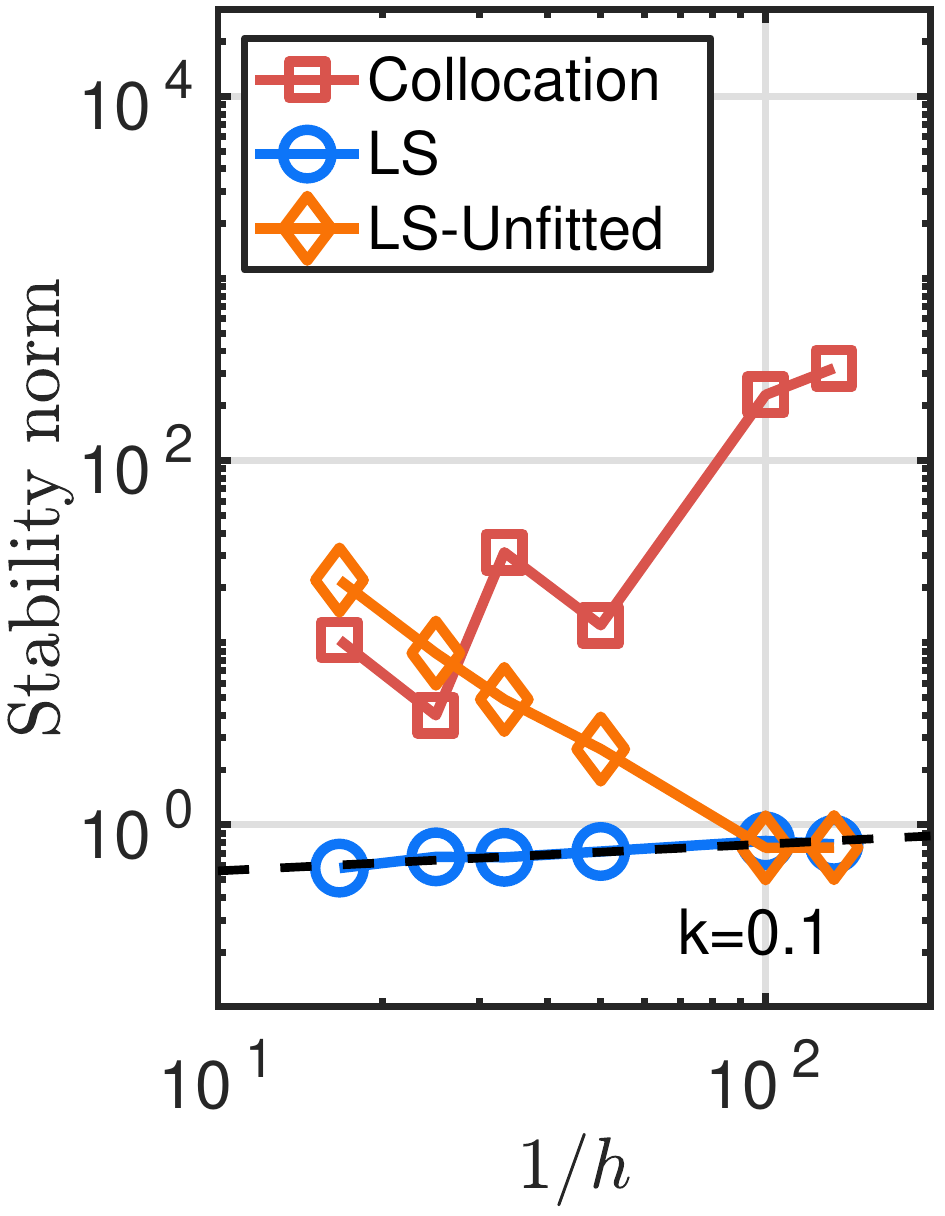}
    \end{tabular}
    \caption{Stability norm as a function of the inverse internodal distance $1/h$ for three different polynomial degrees $p$.}
    \label{fig:experiments:butterfly:href:stabilitynorm}
\end{figure}

In Figure \ref{fig:experiments:butterfly:href:stabilitynorm} we see that
the stability norm of the unfitted RBF-FD-LS method is constant for $p=2$. For $p=4$ and $p=6$ it is partially decaying
until it levels out at a value slightly larger than the stability norm of fitted RBF-FD-LS.
This effect is in line with the behavior of the exterior interpolation points, which move closer to the boundary as $\frac{1}{h}$ gets larger.
In this way the node layout of the unfitted RBF-FD-LS method is getting increasingly similar to the node layout of
the \refereeSecond{fitted} RBF-FD-LS method.

\begin{figure}[h!]
    \centering
    \begin{tabular}{ccc}
        \hspace{0.7cm} $p=2$                                                                          & \hspace{0.7cm} $p=4$ & \hspace{0.7cm} $p=6$ \\
        \includegraphics[width=0.3\linewidth]{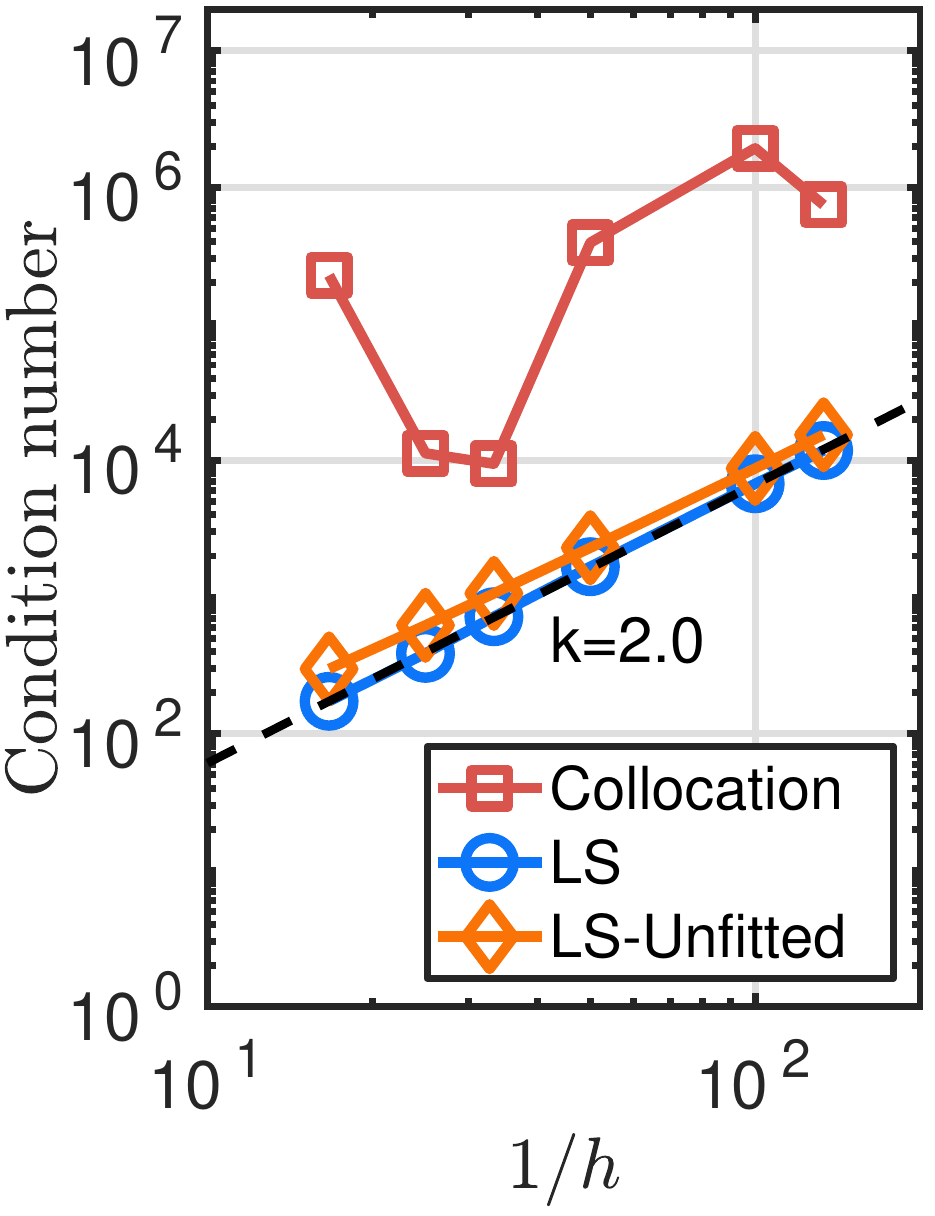} &
        \includegraphics[width=0.3\linewidth]{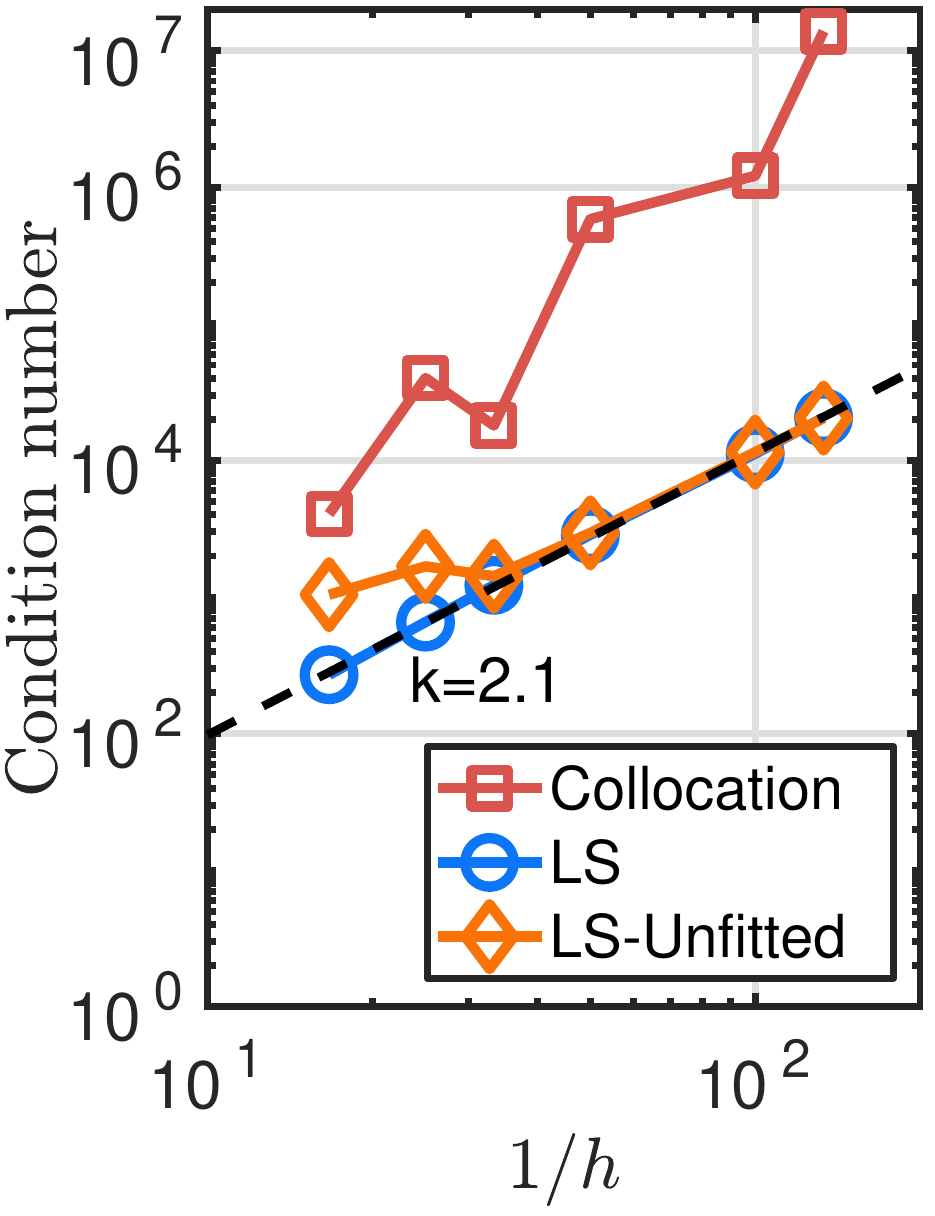} &
        \includegraphics[width=0.3\linewidth]{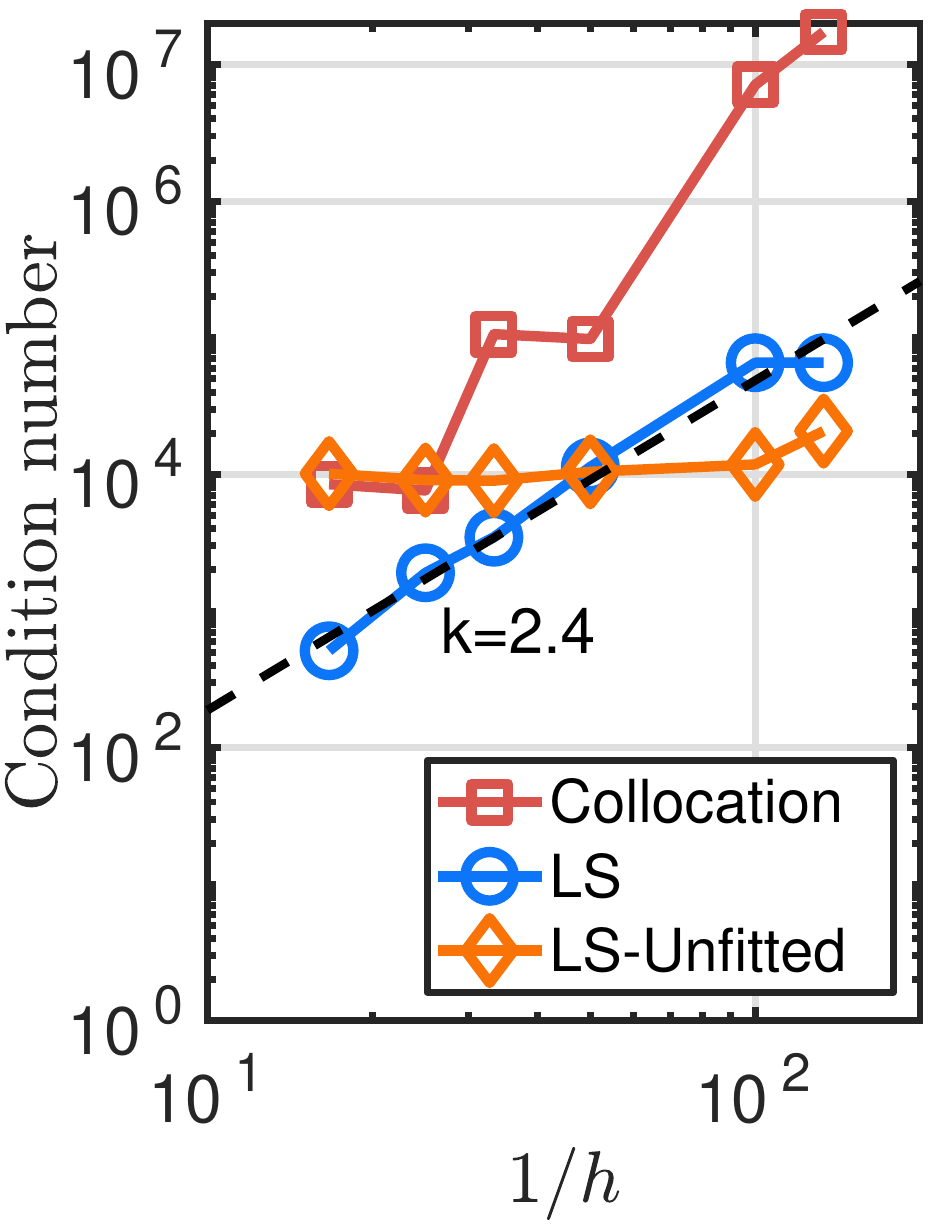}
    \end{tabular}
    \caption{Conditioning of the PDE matrix $D_h$ as a function of the inverse internodal distance $1/h$ for three different polynomial degrees $p$.}
    \label{fig:experiments:butterfly:href:conditioning}
\end{figure}

The condition numbers specific to the collocation, the least-squares and the unfitted least-squares setups are referred to
by $\kappa_{\text{C}}$, $\kappa_{\text{LS}}$ and $\kappa_{\text{U-LS}}$
respectively. In Figure \ref{fig:experiments:butterfly:href:conditioning} we see that when
$p=2$ the condition numbers $\kappa_{\text{LS}}$ and $\kappa_{\text{U-LS}}$ behave as $h^{-2}$
while $\kappa_{\text{C}}$ does not follow any specific pattern (already observed in \cite{tominec2020squares}). Growth with $h^{-2}$ is expected,
since it is known that the conditioning of any discretization involving the Laplacian operator
scales with at least that rate.
In the $p=4$ and $p=6$ cases we see a different behavior of $\kappa_{\text{U-LS}}$ compared with $\kappa_{\text{LS}}$: here
$\kappa_{\text{U-LS}}$ is at first large and remains constant until it coincides with $\kappa_{\text{LS}}$ and starts
growing with approximately $h^{-2}$.
This effect is conceptually very similar to the behavior of the stability norm of the unfitted RBF-FD-LS method which approaches the
stability norm of RBF-FD-LS as $h\to 0$.

\subsection{Approximation properties as the polynomial degree is increased}
\begin{figure}[h!]
    \centering
    \begin{tabular}{cc}
        \hspace{0.7cm} $h=0.05$                                                                     & \hspace{0.7cm} $h=0.015$ \\
        \includegraphics[width=0.31\linewidth]{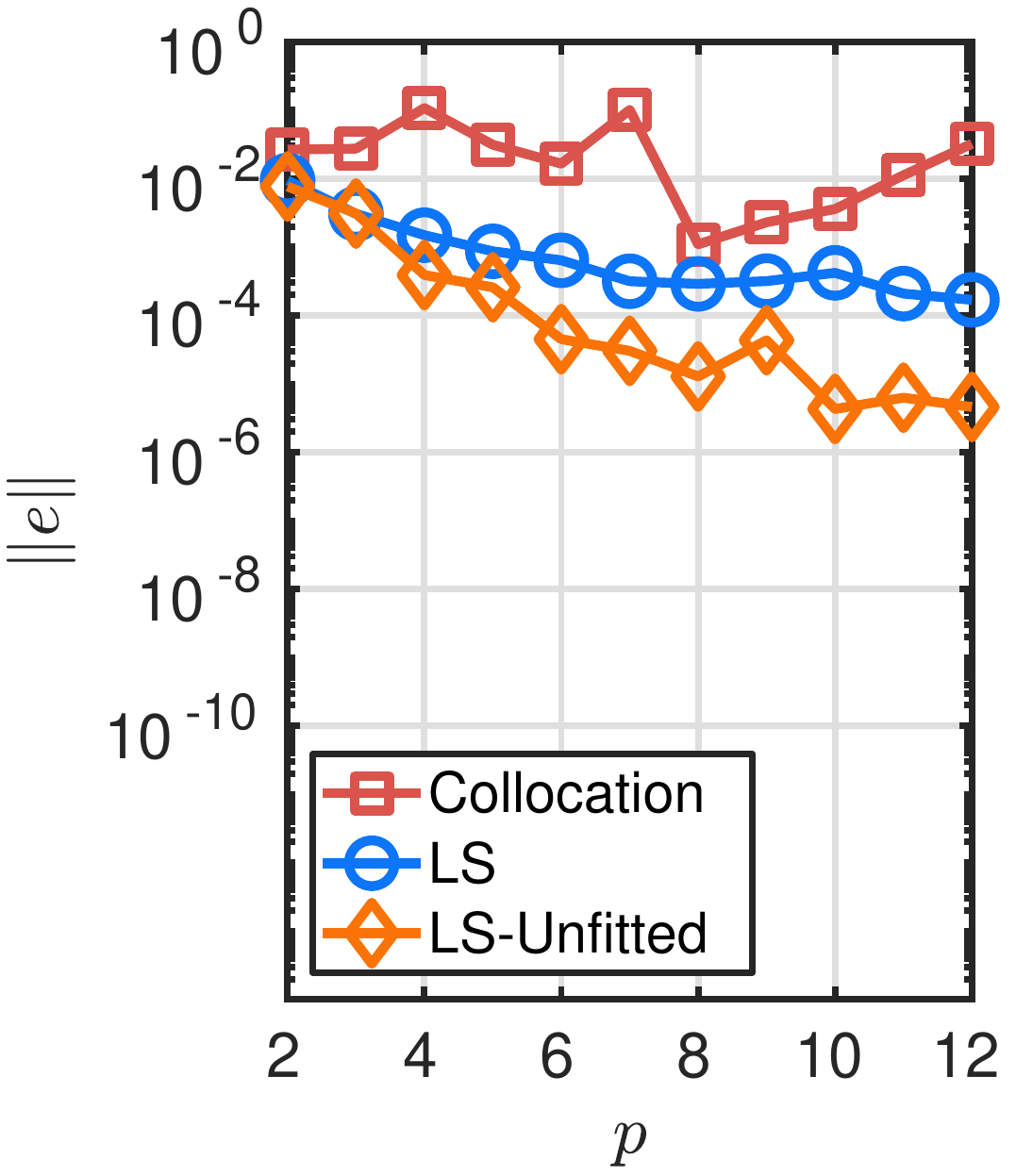} &
        \includegraphics[width=0.31\linewidth]{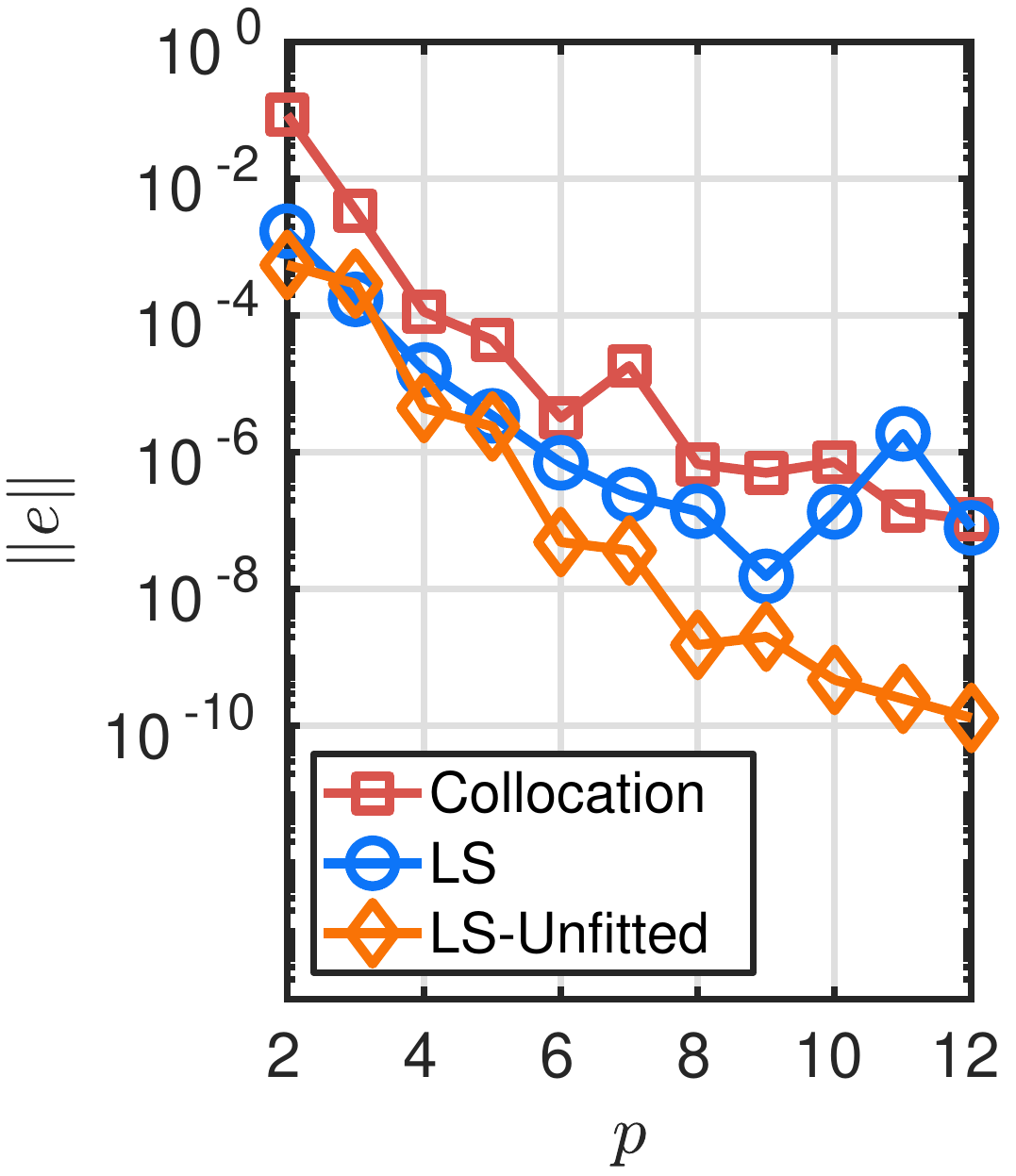}
    \end{tabular}
    \caption{Error as a function of the polynomial degree $p$ for an under-resolved case ($h=0.05$) and a well-resolved case ($h=0.015$) when
        the Franke function is chosen as the exact solution.}
    \label{fig:experiments:butterfly:pref:franke_error}
\end{figure}

\begin{figure}[h!]
    \centering
    \begin{tabular}{cc}
        \hspace{0.7cm} $h=0.05$                                                                          & \hspace{0.7cm} $h=0.015$ \\
        \includegraphics[width=0.31\linewidth]{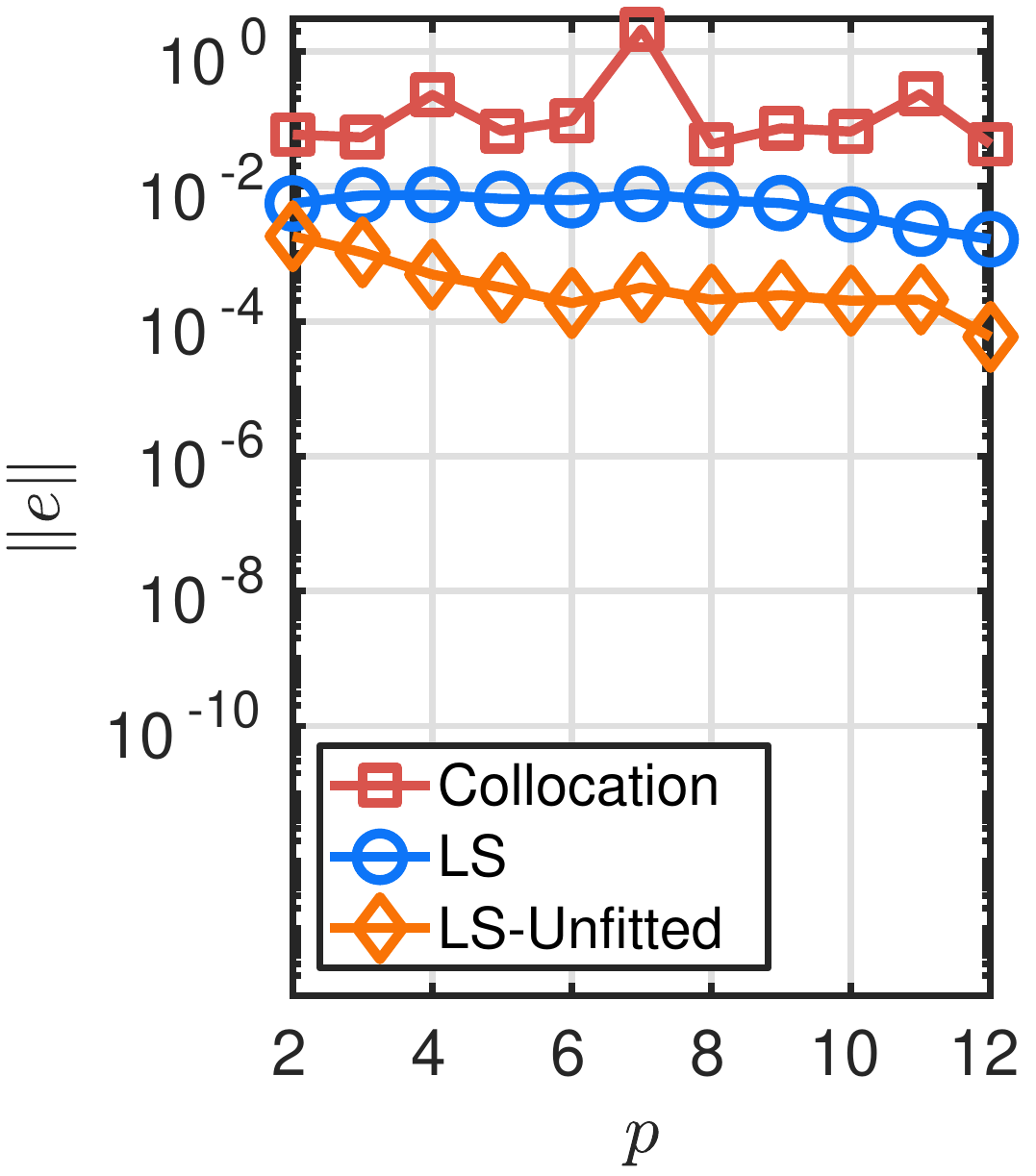} &
        \includegraphics[width=0.31\linewidth]{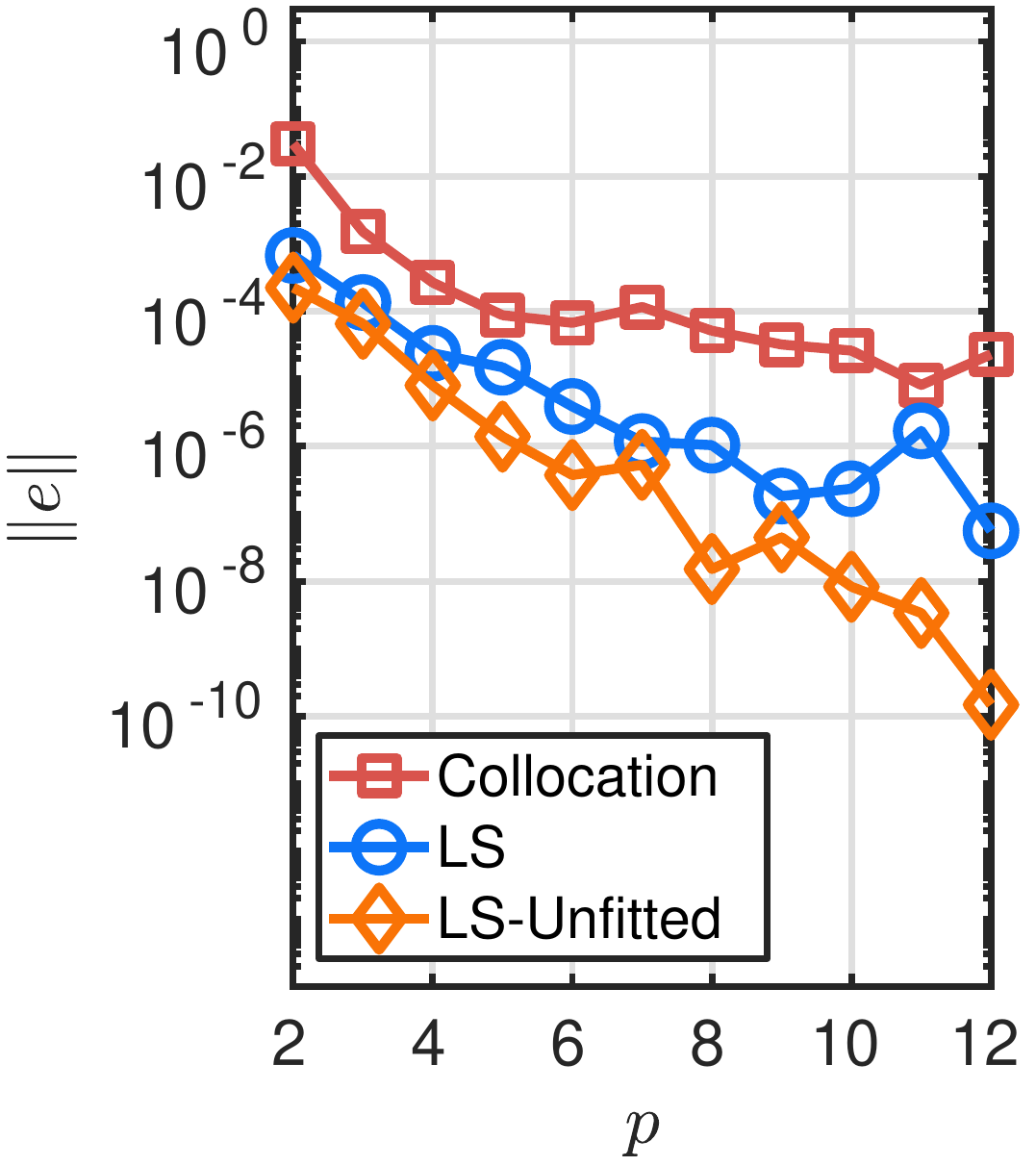}
    \end{tabular}
    \caption{Error as a function of the polynomial degree $p$ for an under-resolved case ($h=0.05$) and a well-resolved case ($h=0.015$) when
        the \refereeSecond{Truncated non-analytic} function is chosen as the exact solution.}
    \label{fig:experiments:butterfly:pref:nonanalytic_error}
\end{figure}

This test gives an insight into the approximation error when the stencil size is increased together with the polynomial degree, while
the internodal distance $h$ is fixed. We choose to work with two fixed values of $h$, namely $h=0.05$ and $0.015$, where the former corresponds to
a case where the approximated solution is unresolved and the latter when the approximated solution is well resolved. \refereeSecond{We use Neumann and Dirichlet boundary conditions in the same way as in 
Section \ref{section:experiments_butterfly:mixedBCs}.}
Results for the Franke and \refereeSecond{Truncated non-analytic} functions
in the role of exact solutions are given in Figure \ref{fig:experiments:butterfly:pref:franke_error} and Figure \ref{fig:experiments:butterfly:pref:nonanalytic_error}
respectively. In both figures the unfitted RBF-FD-LS method is superior in error to RBF-FD-LS and RBF-FD-C, especially when $p$ is large.
This is expected: when the stencil size is increasing, the stencils on the boundary get increasingly more skewed in both fitted setups,
while in the unfitted RBF-FD-LS method the stencils remain fairly unskewed. \refereeSecond{For a discussion about the behavior of the interpolation error when a stencil is skewed see Appendix \ref{sec:appendix:skeweness}.}

\section{Experiments on a 2D drilled sprocket domain with interior boundary conditions}
\label{section:experiments_sprocket}
In this section we consider the drilled 24-tooth sprocket from Figure \ref{fig:intro:sprocket} as our  computational domain.
We solve \eqref{eq:model:Poisson}, but in addition to the exterior mixed boundary conditions
we also introduce interior Dirichlet and Neumann
boundary conditions and compare the unfitted RBF-FD-LS method,
RBF-FD-LS and RBF-FD-C with a focus on \refereeSecond{the approximation error under node refinement measured in different norms, and} the spatial distribution of the error.
\refereeSecond{The point sets for all methods are chosen in the same way as in Section \ref{sec:experiments_butterfly:pointsets}.}

\subsection{Domain $\Omega$}
The points of the domain were acquired from a simple monochrome sprocket drawing.
By using tools from mathematical morphology we extracted thin borders of the object and applied the Harris feature detector to produce a set of points of interest (eg. corners) adequately representing the object shape.
Each connected component was then parametrized separately via linear arc-length and resampled equidistantly (using the one-dimensional RBF-FD method for interpolation).

\subsection{Solution function}
\begin{figure}[h!]
    \centering
    \includegraphics[width=0.46\linewidth]{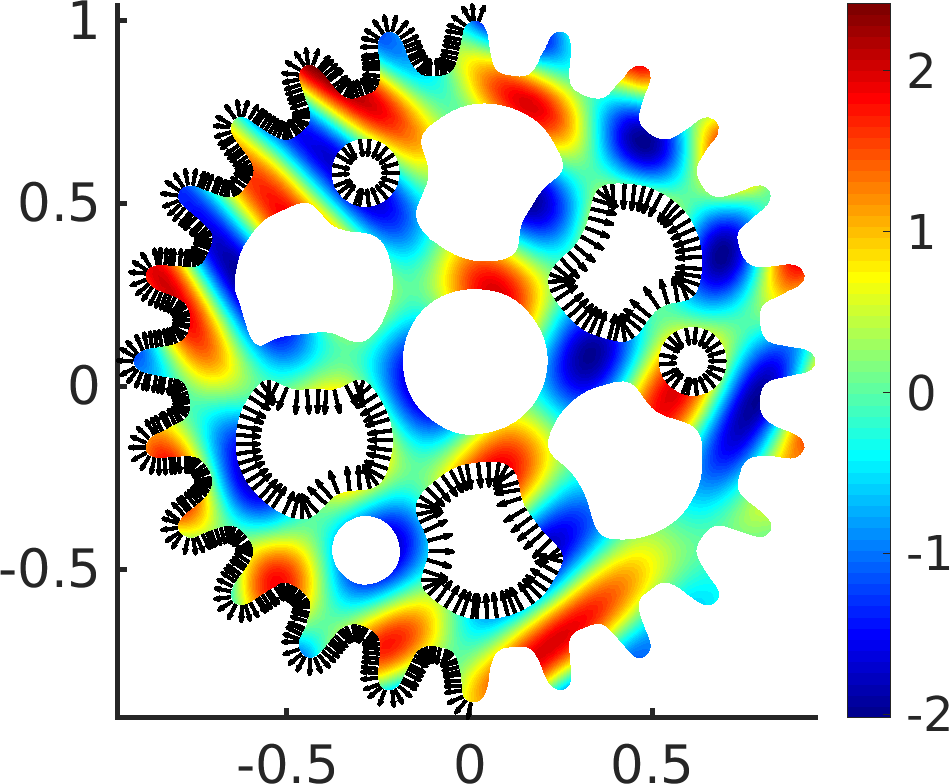}
    \caption{Solution function $u_3$ on a sprocket. The black outward normals indicate the locations of the Neumann condition. The locations over the boundary which do not
        include an arrow indicate the Dirichlet condition.}
    \label{fig:experiments:solutionFunction2}
\end{figure}

Inspired by \cite{BFFB17} the solution function used to
compute the right-hand-sides of \eqref{eq:model:Poisson} is set to:
$$u_3 = \sin\left(3\pi y^2 + 4.5\pi x \right) - \cos\left(4\pi y - 3\pi x^2 \right).$$
Its visual representation is given in Figure \ref{fig:experiments:solutionFunction2}.
The setup of boundary conditions
that we use in the experiments can be seen in Figure \ref{fig:experiments:solutionFunction2}, where the black normals
over the boundaries indicate the locations of the Neumann boundary conditions and the Dirichlet boundary conditions are employed where there is no marker.

\subsection{Convergence under node refinement}
We numerically verify that the solution in all three setups (RBF-FD-C, RBF-FD-LS and unfitted RBF-FD-LS) converges for polynomial degrees $p=2$, $p=4$, $p=6$
as $h\to 0$.
\refereeSecond{We compute the relative errors in the following three norms:
    \begin{equation}
        \|e\|_1 = \frac{\|u_h(Y) - u(Y)\|_1}{\|u(Y)\|_1},\quad \|e\|_2 = \frac{\|u_h(Y) - u(Y)\|_2}{\|u(Y)\|_2},\quad \|e\|_\infty = \frac{\|u_h(Y) - u(Y)\|_\infty}{\|u(Y)\|_\infty},
    \end{equation}
    where $u_h(Y)$ and $u(Y)$ are the approximate solution and the exact solution consecutively, both sampled in $Y$ points.}

\refereeSecond{The oversampling parameter is set to $q=5$. For a more detailed study of the
    oversampling parameter and its impact on the approximation error and the stability norm, we refer the reader to \cite{tominec2020squares}.}

\refereeSecond{
    The convergence results measured in $1$-norm, $2$-norm and $\infty$-norm are presented in Figure \ref{fig:experiments:sprocket:href_error_l1},
    \ref{fig:experiments:sprocket:href_error_l2} and Figure \ref{fig:experiments:sprocket:href_error_linf} consecutively. Here
    we see that for the considered problem, the unfitted RBF-FD-LS method is the most accurate among all three setups, in all three different norms.
    \begin{figure}[h!]
        \centering
        \begin{tabular}{ccc}
            \hspace{0.7cm} $p=2$                                                                                    & \hspace{0.7cm} $p=4$ & \hspace{0.7cm} $p=6$ \\
            \includegraphics[width=0.3\linewidth]{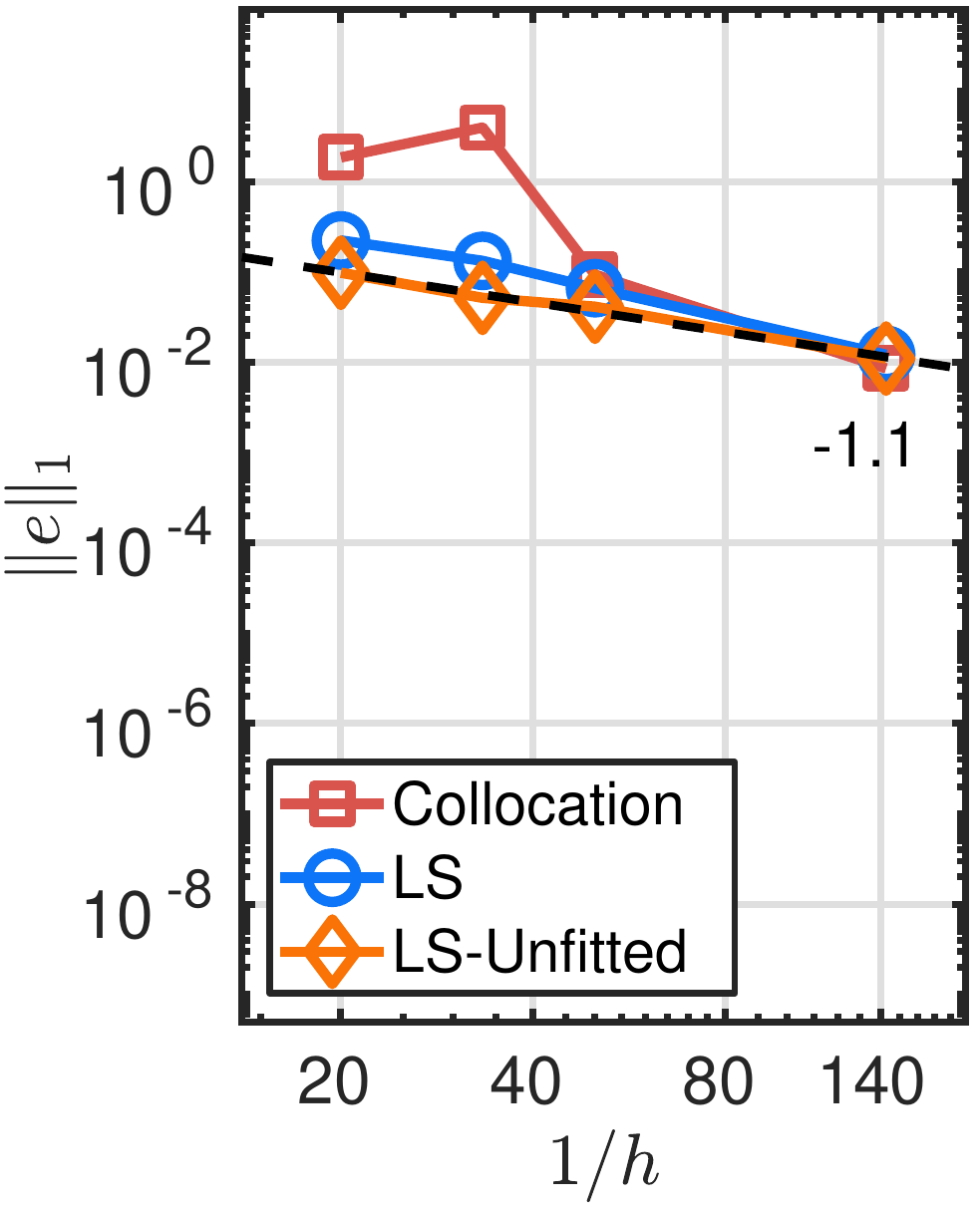} &
            \includegraphics[width=0.3\linewidth]{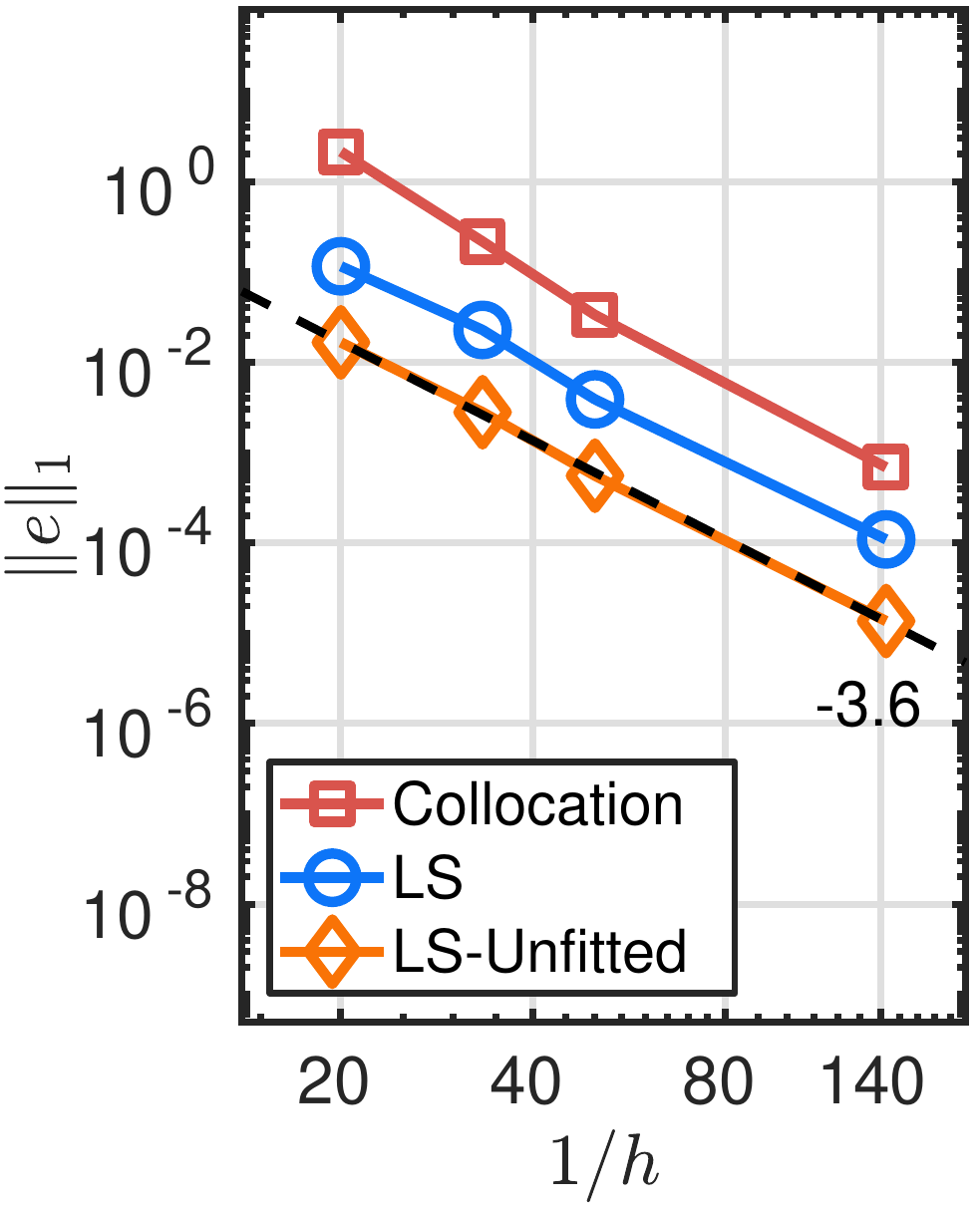} &
            \includegraphics[width=0.3\linewidth]{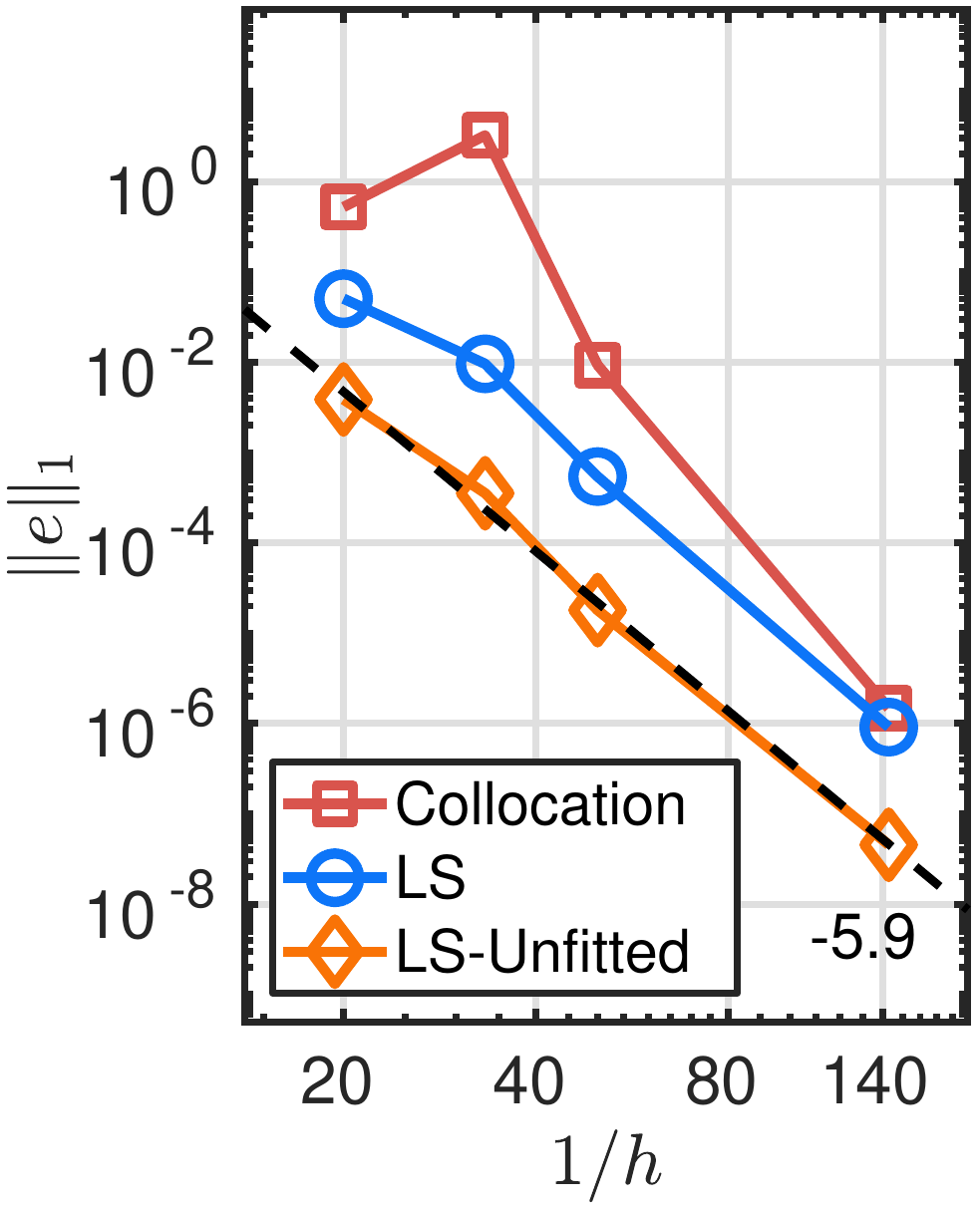}
        \end{tabular}
        \caption{Sprocket case. Error in $1$-norm as a function of the inverse internodal distance $1/h$ for three different polynomial degrees $p$.}
        \label{fig:experiments:sprocket:href_error_l1}
    \end{figure}
    \begin{figure}[h!]
        \centering
        \begin{tabular}{ccc}
            \hspace{0.7cm} $p=2$                                                                                    & \hspace{0.7cm} $p=4$ & \hspace{0.7cm} $p=6$ \\
            \includegraphics[width=0.3\linewidth]{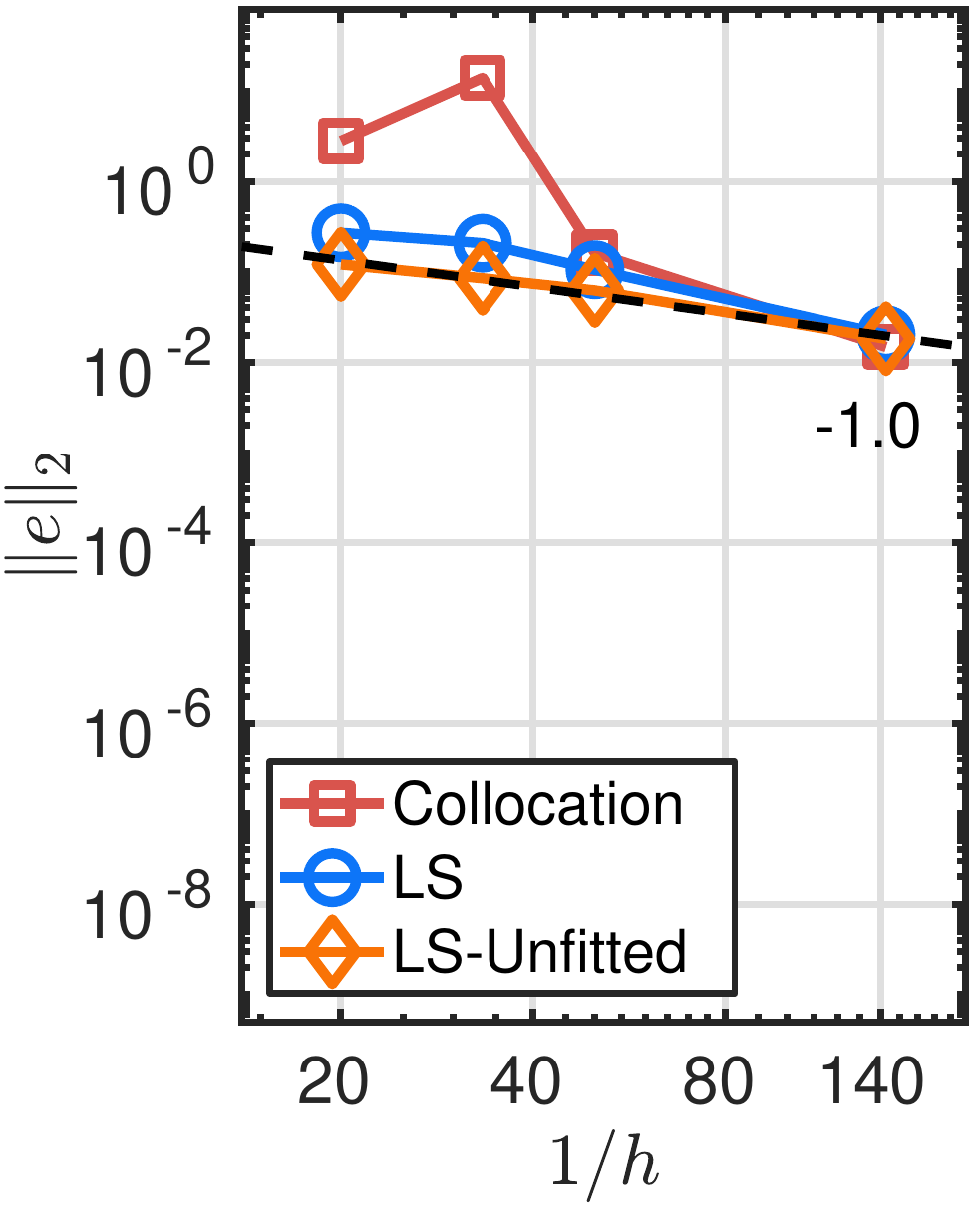} &
            \includegraphics[width=0.3\linewidth]{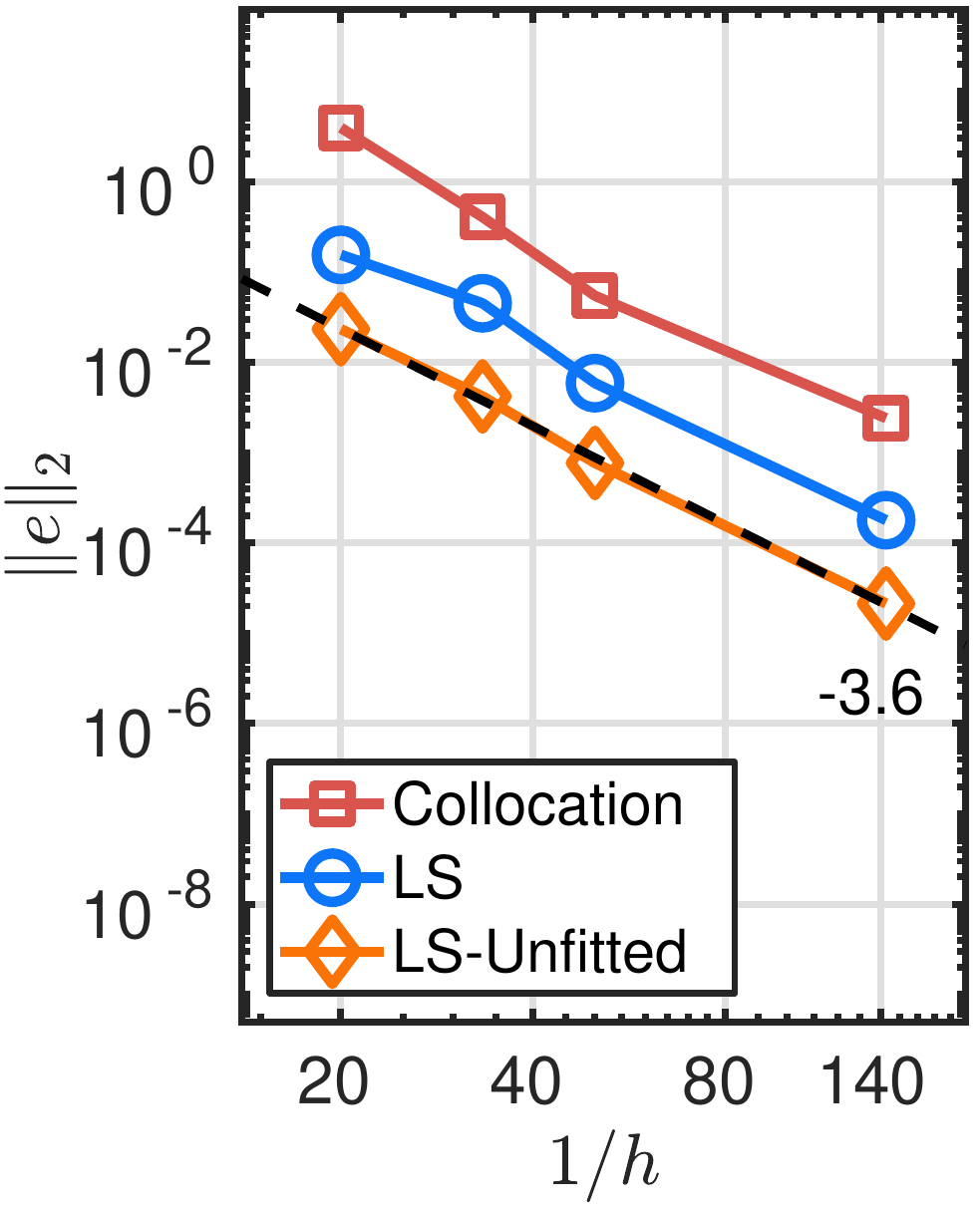} &
            \includegraphics[width=0.3\linewidth]{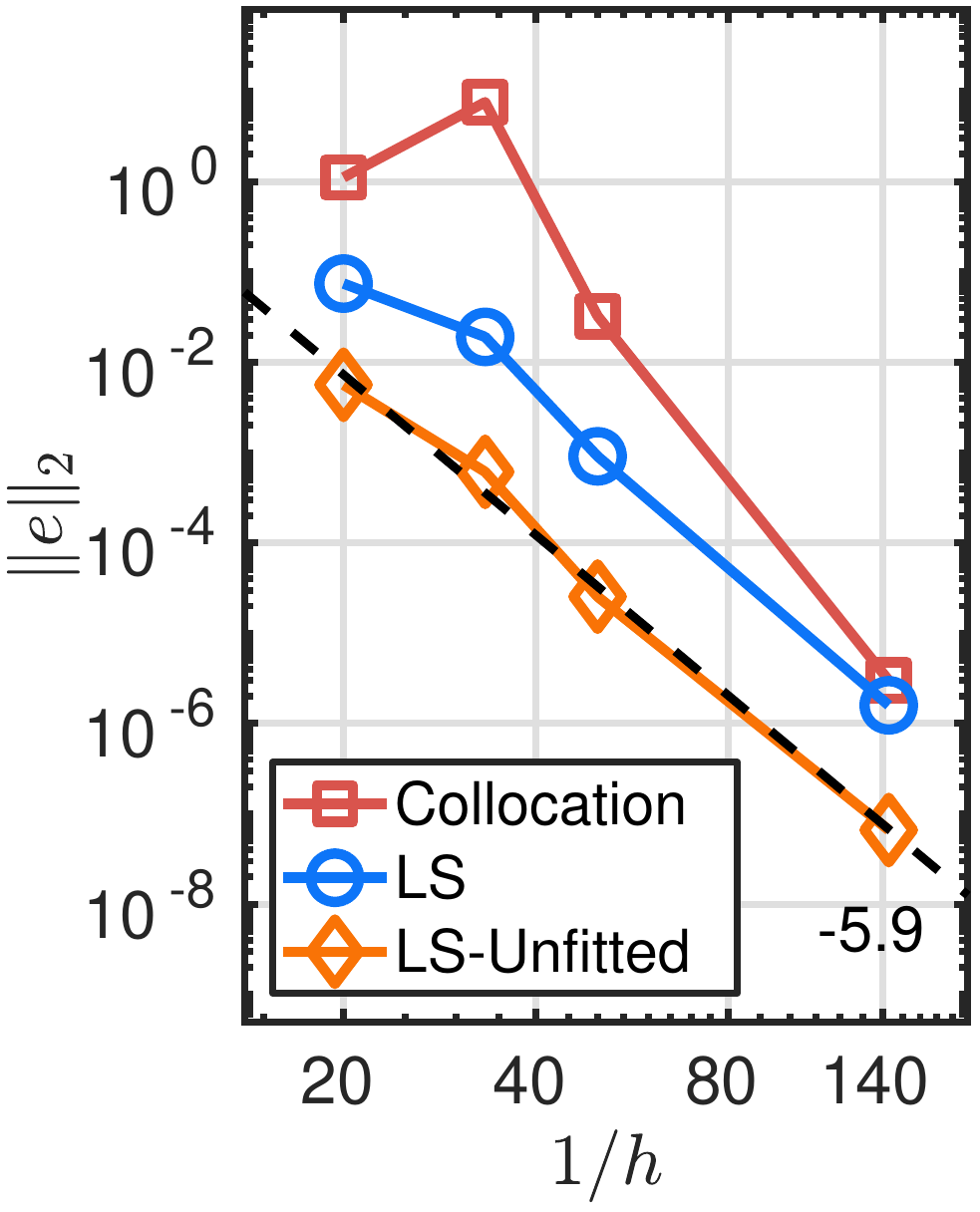}
        \end{tabular}
        \caption{Sprocket case. Error in $2$-norm as a function of the inverse internodal distance $1/h$ for three different polynomial degrees $p$.}
        \label{fig:experiments:sprocket:href_error_l2}
    \end{figure}
    \begin{figure}[h!]
        \centering
        \begin{tabular}{ccc}
            \hspace{0.7cm} $p=2$                                                                                      & \hspace{0.7cm} $p=4$ & \hspace{0.7cm} $p=6$ \\
            \includegraphics[width=0.3\linewidth]{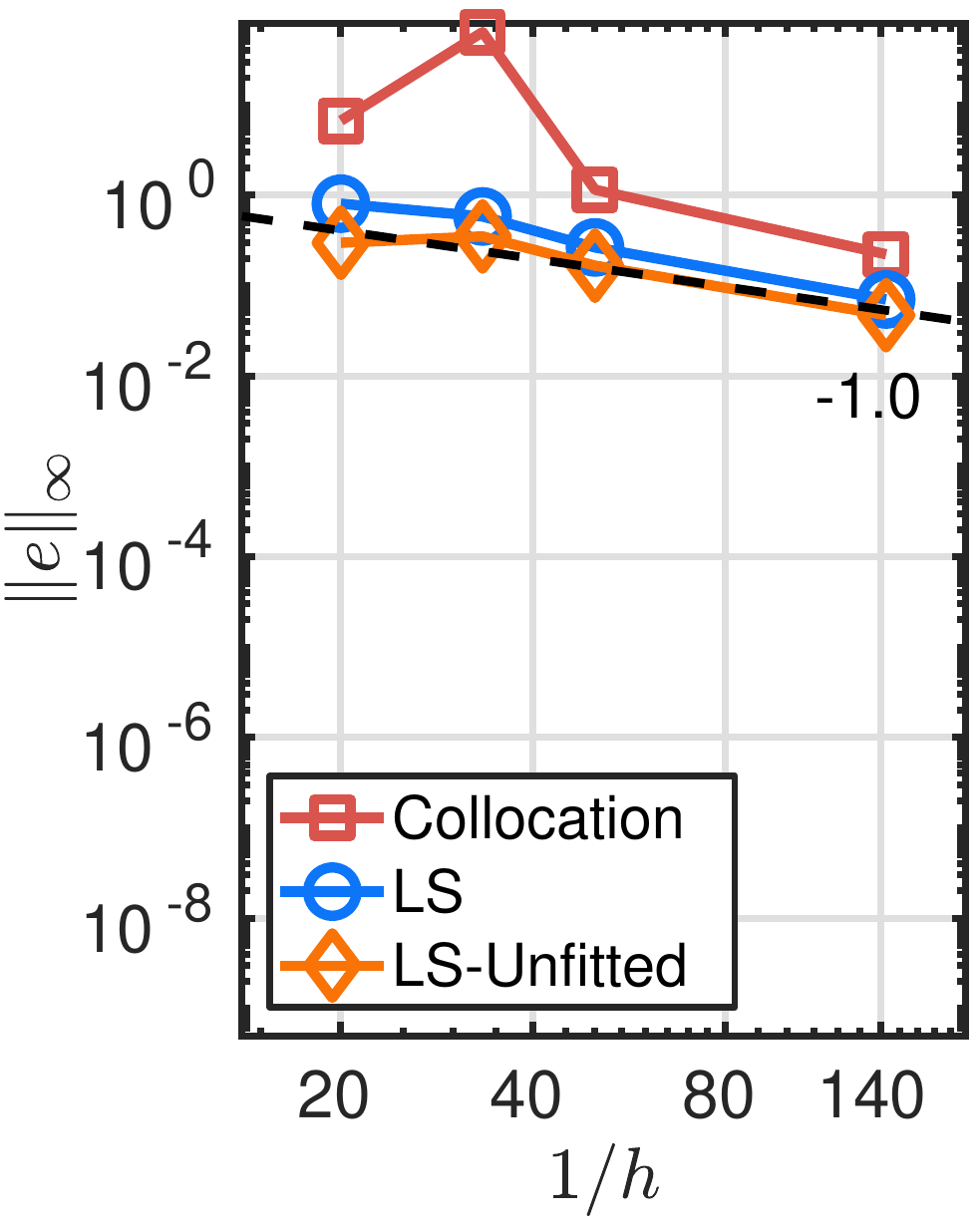} &
            \includegraphics[width=0.3\linewidth]{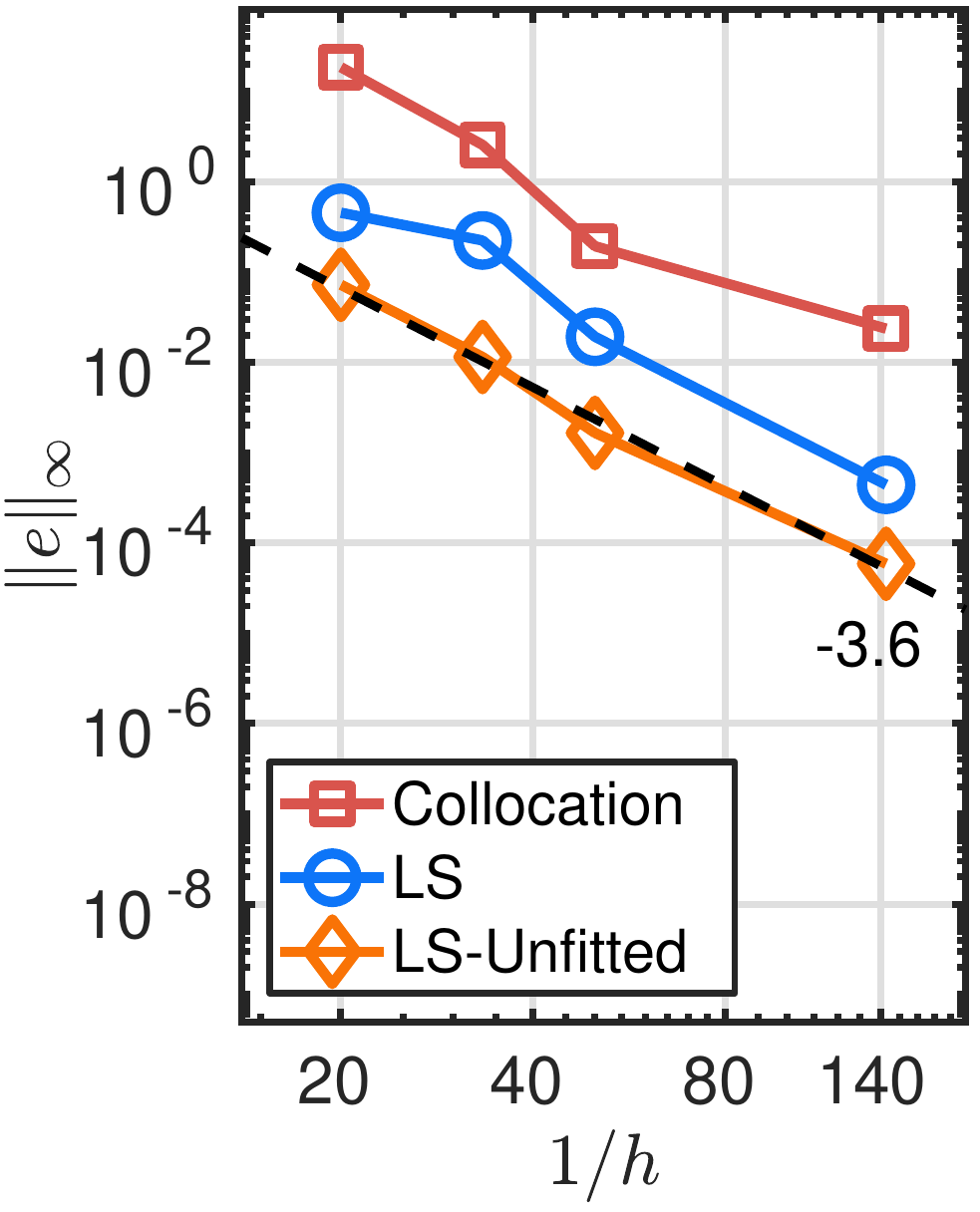} &
            \includegraphics[width=0.3\linewidth]{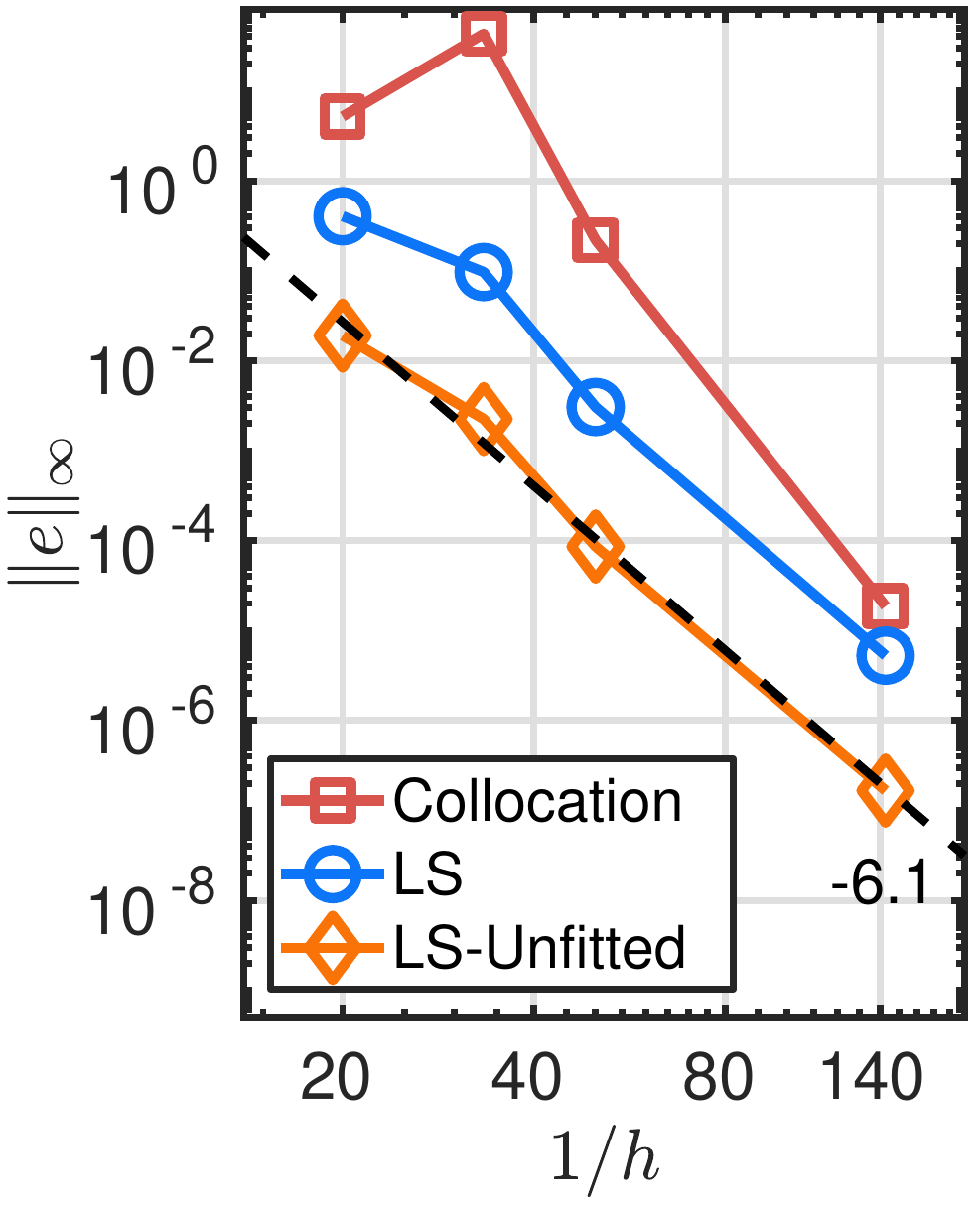}
        \end{tabular}
        \caption{Sprocket case. Error in $\infty$-norm as a function of the inverse internodal distance $1/h$ for three different polynomial degrees $p$.}
        \label{fig:experiments:sprocket:href_error_linf}
    \end{figure}
}
\subsection{Spatial distribution of the error}
\begin{figure}[h!]
    \centering
    \begin{tabular}{ccc}
        \textbf{RBF-FD-LS}                                                                                               & \textbf{Unfitted RBF-FD-LS} & \textbf{RBF-FD-C} \vspace{0.07cm} \\
        \multicolumn{3}{c}{$p=2$}\vspace{0.3cm}                                                                                                                                            \\
        \includegraphics[width=0.295\linewidth]{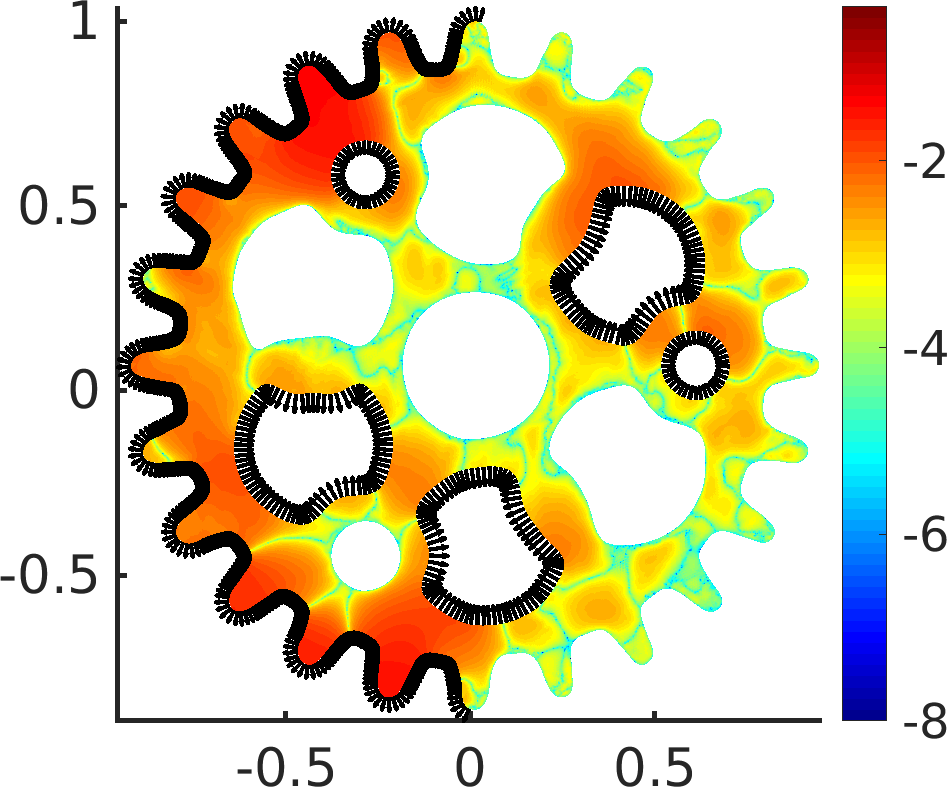}          &
        \includegraphics[width=0.295\linewidth]{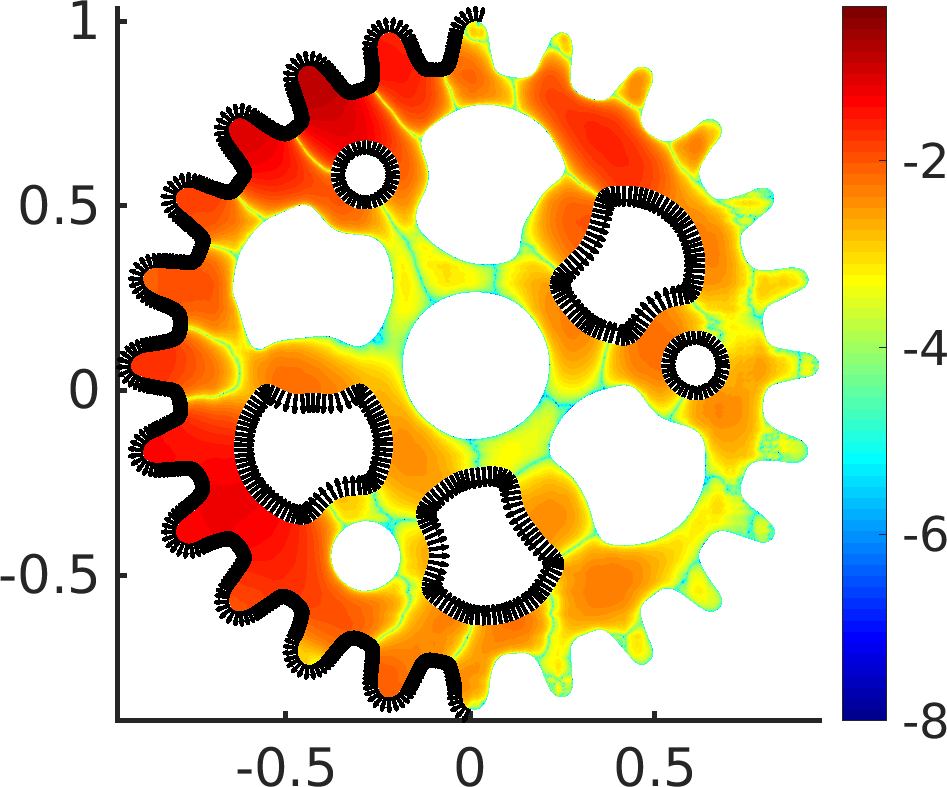} &
        \includegraphics[width=0.295\linewidth]{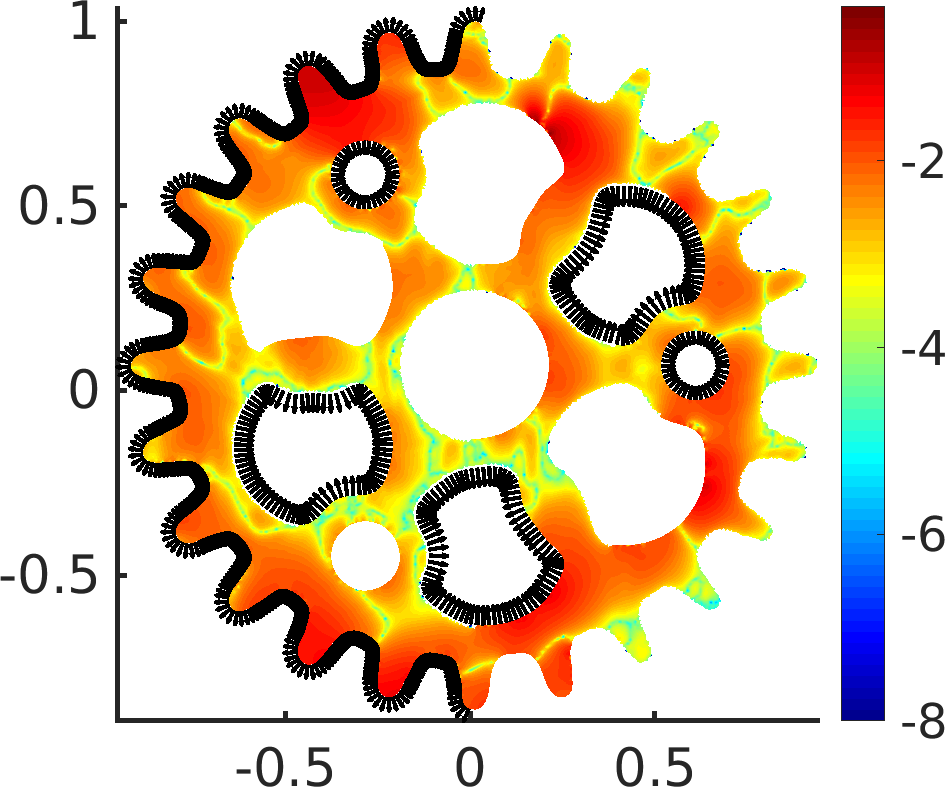} \vspace{0.07cm}                                                             \\

        \multicolumn{3}{c}{$p=4$}\vspace{0.3cm}                                                                                                                                            \\
        \includegraphics[width=0.295\linewidth]{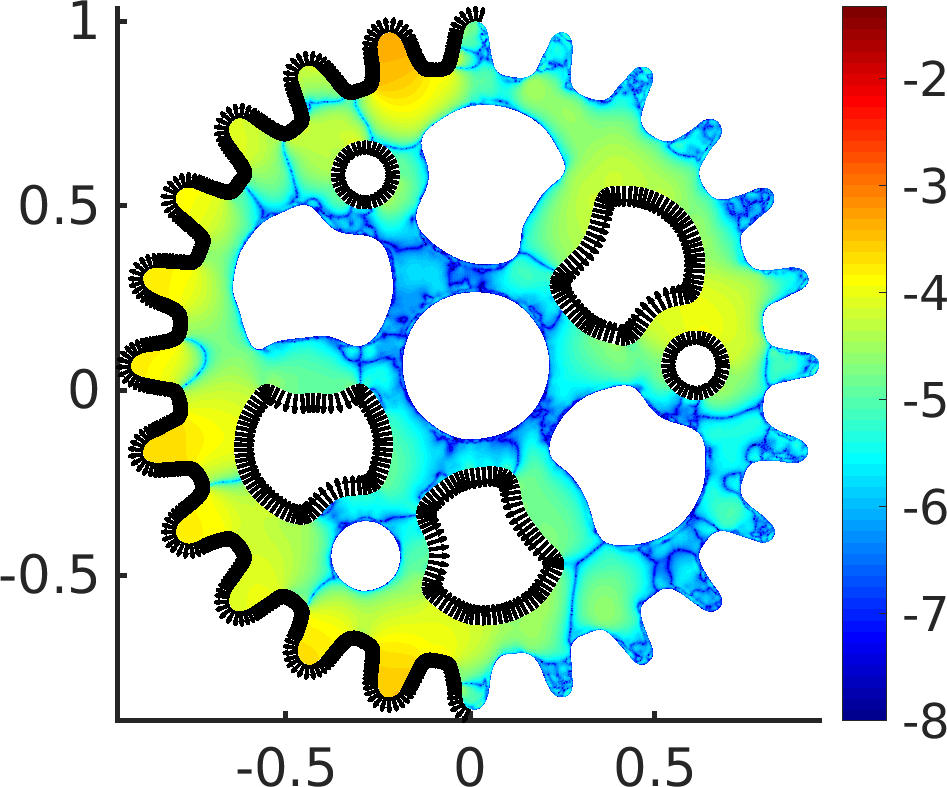}          &
        \includegraphics[width=0.295\linewidth]{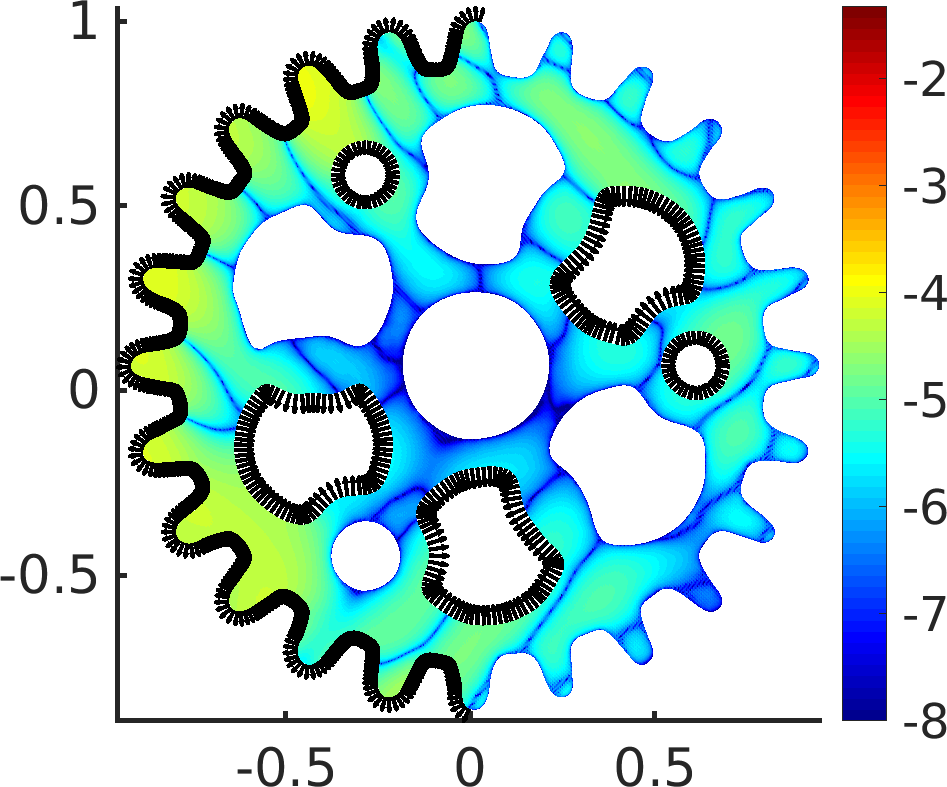} &
        \includegraphics[width=0.295\linewidth]{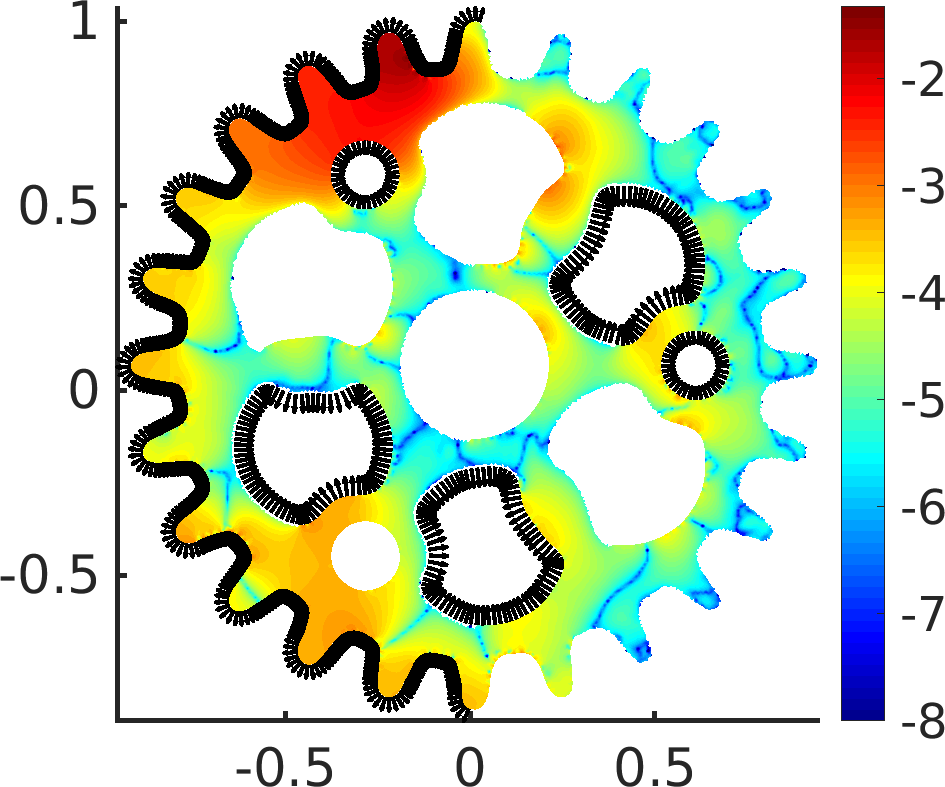} \vspace{0.07cm}                                                             \\

        \multicolumn{3}{c}{$p=6$}\vspace{0.3cm}                                                                                                                                            \\
        \includegraphics[width=0.295\linewidth]{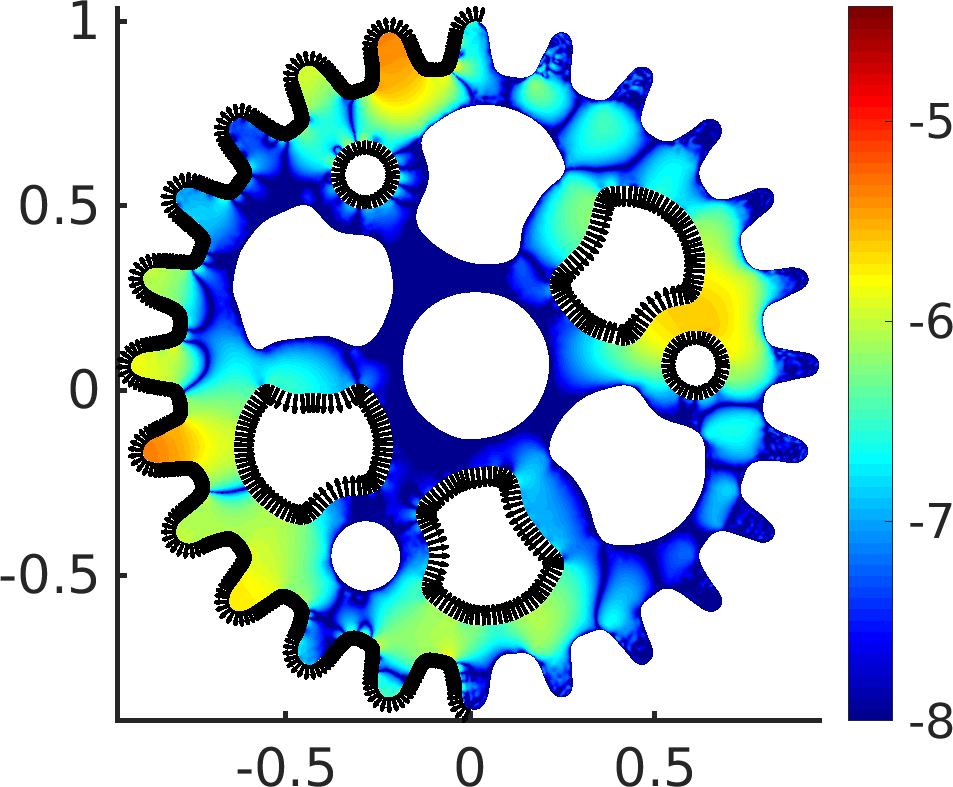}          &
        \includegraphics[width=0.295\linewidth]{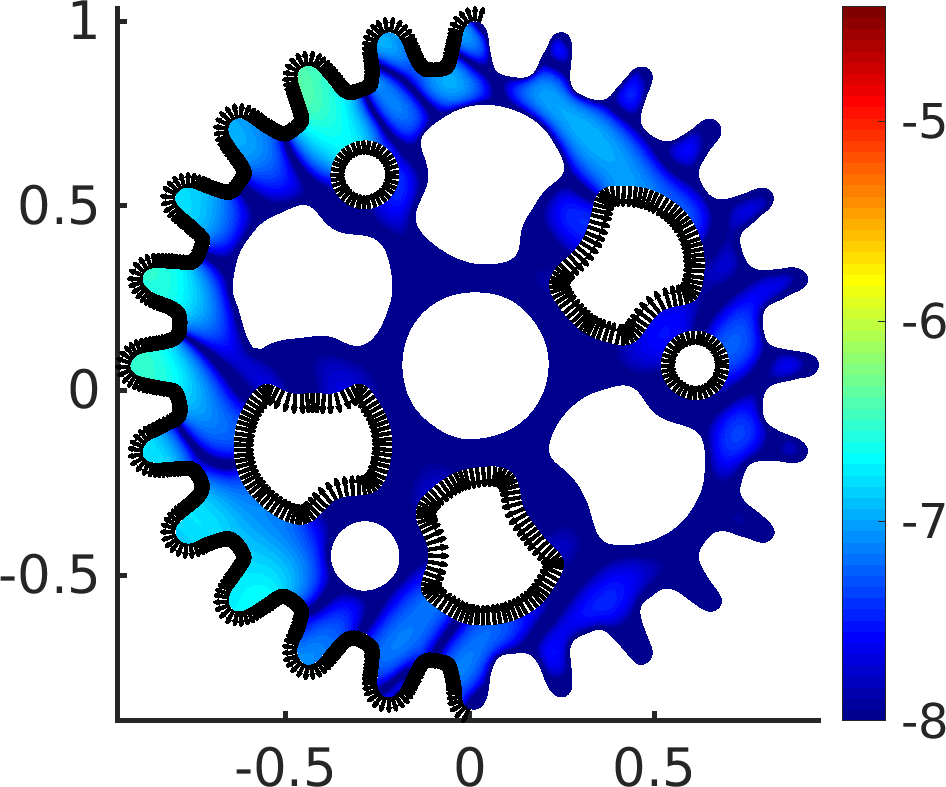} &
        \includegraphics[width=0.295\linewidth]{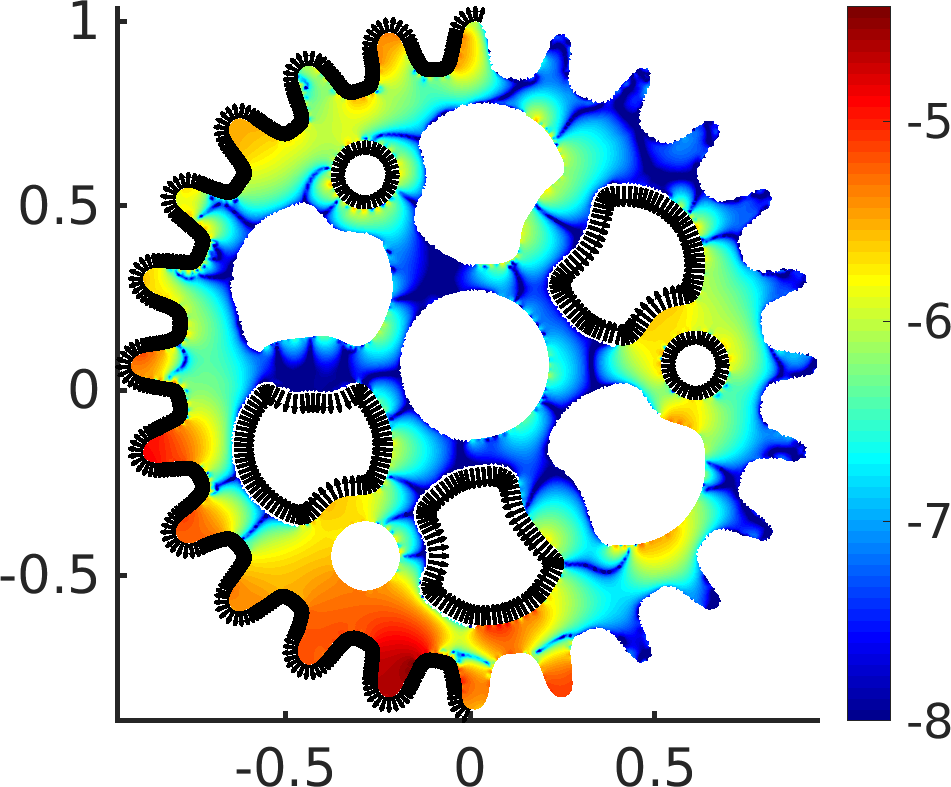}
    \end{tabular}
    \caption{A comparison of the error distribution (logarithmic scale) over a drilled 24-tooth sprocket .
        The outward normals indicate the parts of $\partial\Omega$ where the Neumann condition is imposed. The internodal distance
        is set to $h=0.007$, the oversampling parameter to $q=5$ and the polynomial degree to $p=2$ (first row), $p=4$ (second row),
        and $p=6$ (third row).}
    \label{fig:experiments:sprocket:abs_error_spatial}
\end{figure}

Next, we examine the spatial distribution of the error. In Figure \ref{fig:experiments:sprocket:abs_error_spatial}
we display the spatial distribution of error in logarithmic scale (the negative integers $k$ in the colorbar are the exponents $10^k$) for
$p=2$, $p=4$, $p=6$. When the stencils are small (the case $p=2$) the error distributions for all three settings look fairly similar.
A slight increase in the stencil size (the case $p=4$) causes RBF-FD-C to accumulate more error around the boundaries compared to the
other two setups.
Increasing the stencil size even further (the case $p=6$) reveals that RBF-FD-C again collects
more error around the boundaries than the other
two setups, but also that RBF-FD-LS has larger errors compared with the unfitted RBF-FD-LS method. This is ascribed to the effect of using the extended interpolation points
which make the stencils at the domain boundary less skewed (\refereeSecond{see Appendix \ref{sec:appendix:skeweness}.})

\section{Experiments on a 3D diaphragm domain}
\label{section:experiments_diaphragm}
\begin{figure}[h!]
    \centering
    \begin{tabular}{lrc}
        \includegraphics[width=0.42\linewidth]{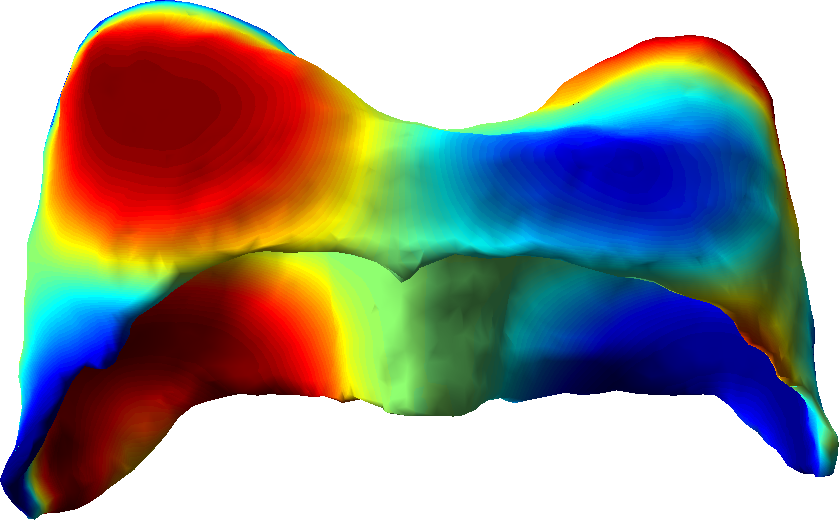} &
        \includegraphics[width=0.3\linewidth]{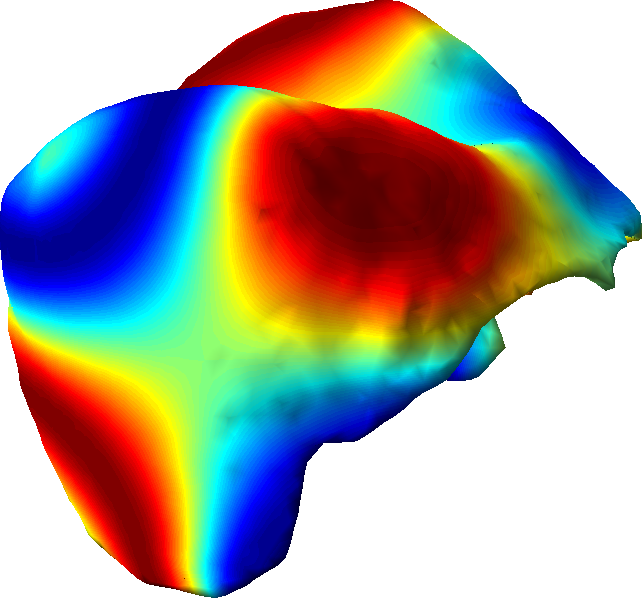}  &
        \multirow{2}{*}{\includegraphics[width=0.065\linewidth]{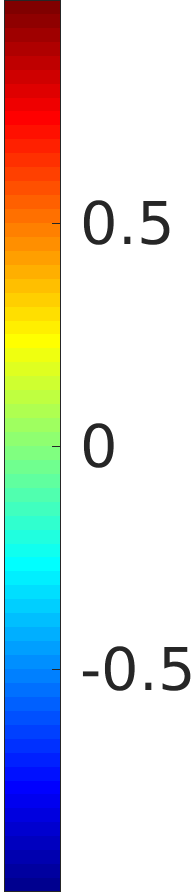}} \\
        \includegraphics[width=0.33\linewidth]{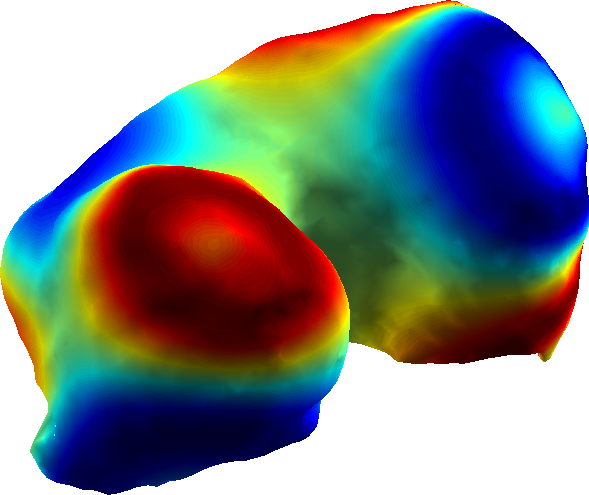} &
        \includegraphics[width=0.33\linewidth]{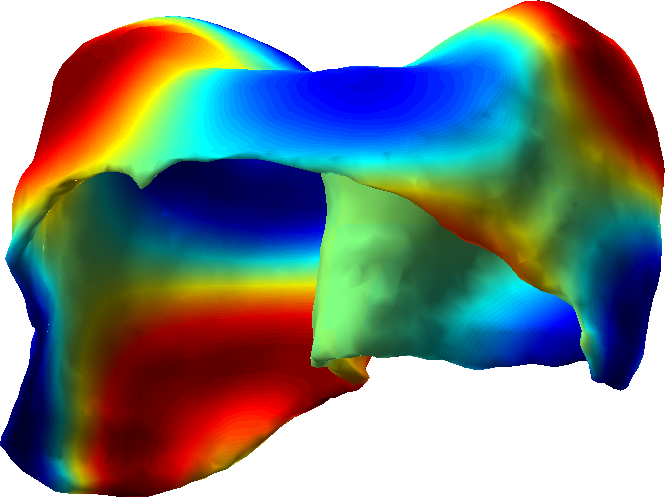} &
    \end{tabular}

    \caption{The solution function $u_4=sin(6\pi x y z)$ over a diaphragm viewed from four different angles.}
    \label{fig:experiments:diaphragm:solutionFunction}
\end{figure}

To show that the previously introduced results also generalize to a realistic scenario in three dimensions we now consider a
thoracic diaphragm of
a human being in
the role of a computational domain $\Omega$. The diaphragm is
displayed in Figure \ref{fig:experiments:diaphragm:solutionFunction} from three different angles.
The diaphragm is a thin and non-convex geometry, of which the thickness is approximately 100-times smaller than the largest circumference
over its surface. Interpolation points that conform to the diaphragm are thus hard to obtain, which is a good
motivation to use the unfitted RBF-FD-LS method for computing a solution to a PDE. In this section we solve the same problem as in
\eqref{eq:model:Poisson}, but this time in three dimensions and only using the unfitted RBF-FD-LS method.
\subsection{Solution function}
The solution function used to compute the right-hand-sides of a 3D equivalent of \eqref{eq:model:Poisson} is given by:
$$u_4(x,y,z) = \sin (6\pi x y z).$$
The function $u_4$ over the diaphragm is drawn in Figure \ref{fig:experiments:diaphragm:solutionFunction}.
The Dirichlet boundary condition is enforced over the whole top surface of the diaphragm,
and the Neumann boundary condition is enforced over the whole bottom surface of the diaphragm.
\subsection{Point sets}
\begin{figure}[h!]
    \centering
    \includegraphics[width=0.42\linewidth]{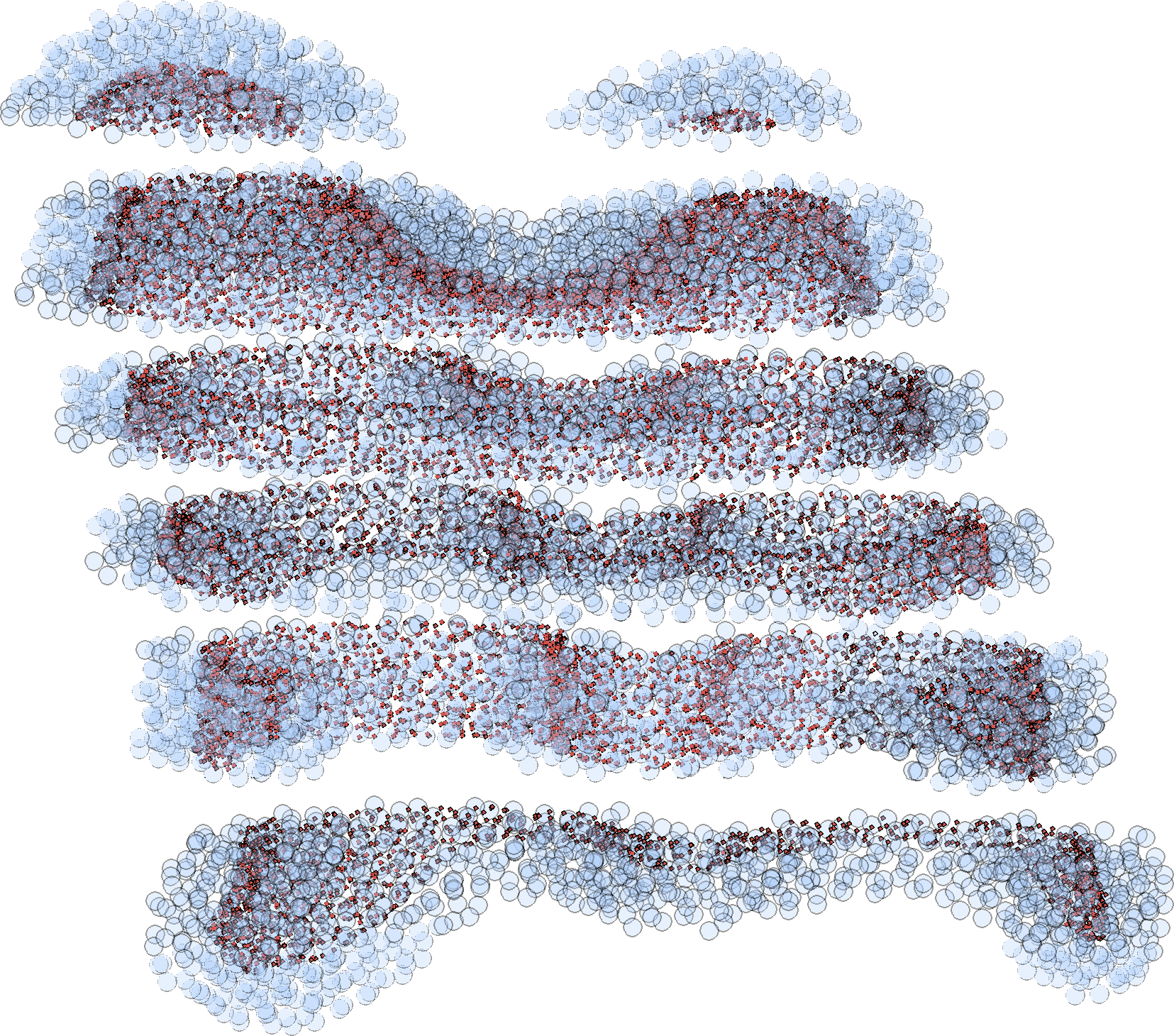}
    \includegraphics[width=0.42\linewidth]{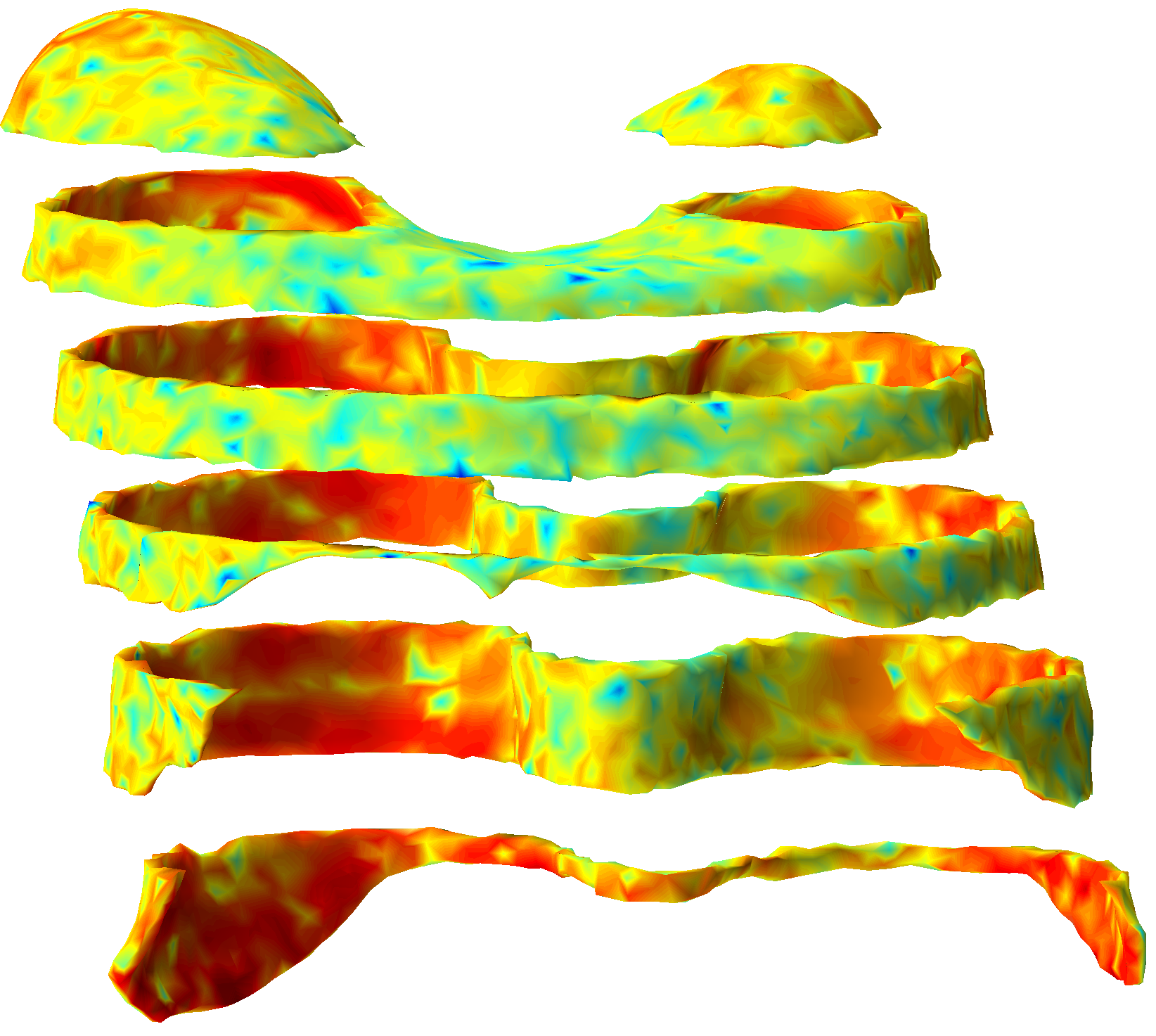}
    \includegraphics[width=0.04\linewidth]{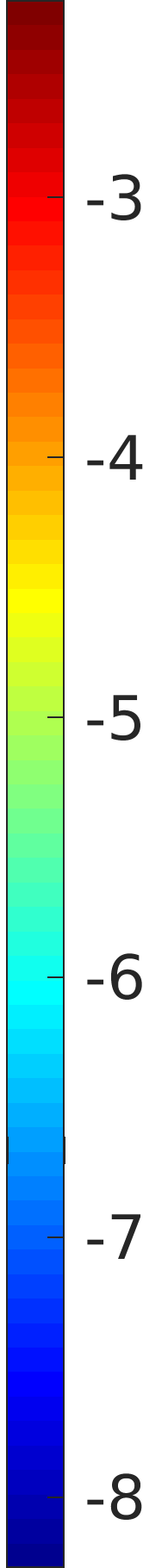}
    \caption{Both pictures show a diaphragm which is sliced for visualization purposes only. Left: interpolation points (blue markers) are placed over the  evaluation points (smaller red markers). The evaluation points discretize the diaphragm.
        Right: The spatial distribution of the magnitude of the relative error in the logarithmic scale when $h=0.066$, $p=5$ and $q=10$.}
    \label{fig:experiments:diaphragm:absError}
\end{figure}

The interpolation and evaluation points are -- conceptually speaking --
constructed in the same way as in the two-dimensional cases.
We start by using a surface point-cloud of the diaphragm \cite{Larsson20} in place of the boundary
evaluation points, which are placed in a three-dimensional box that contains interpolation points computed using an algorithm from
\cite{Kiera19}. Then we enforce \refereeSecond{$q=10$} evaluation
points (again computed by an algorithm from \cite{Kiera19}) around each interpolation point and after that remove those evaluation
points which are placed outside of the diaphragm.
At last we remove the interpolation points according to the criterion given in Section \ref{sec:method:linearindep}. An instance of
the resulting two point sets can be observed in the left image of Figure \ref{fig:experiments:diaphragm:absError}.
\subsection{Convergence under node refinement}
\begin{figure}[h!]
    \centering
    \includegraphics[width=0.37\linewidth]{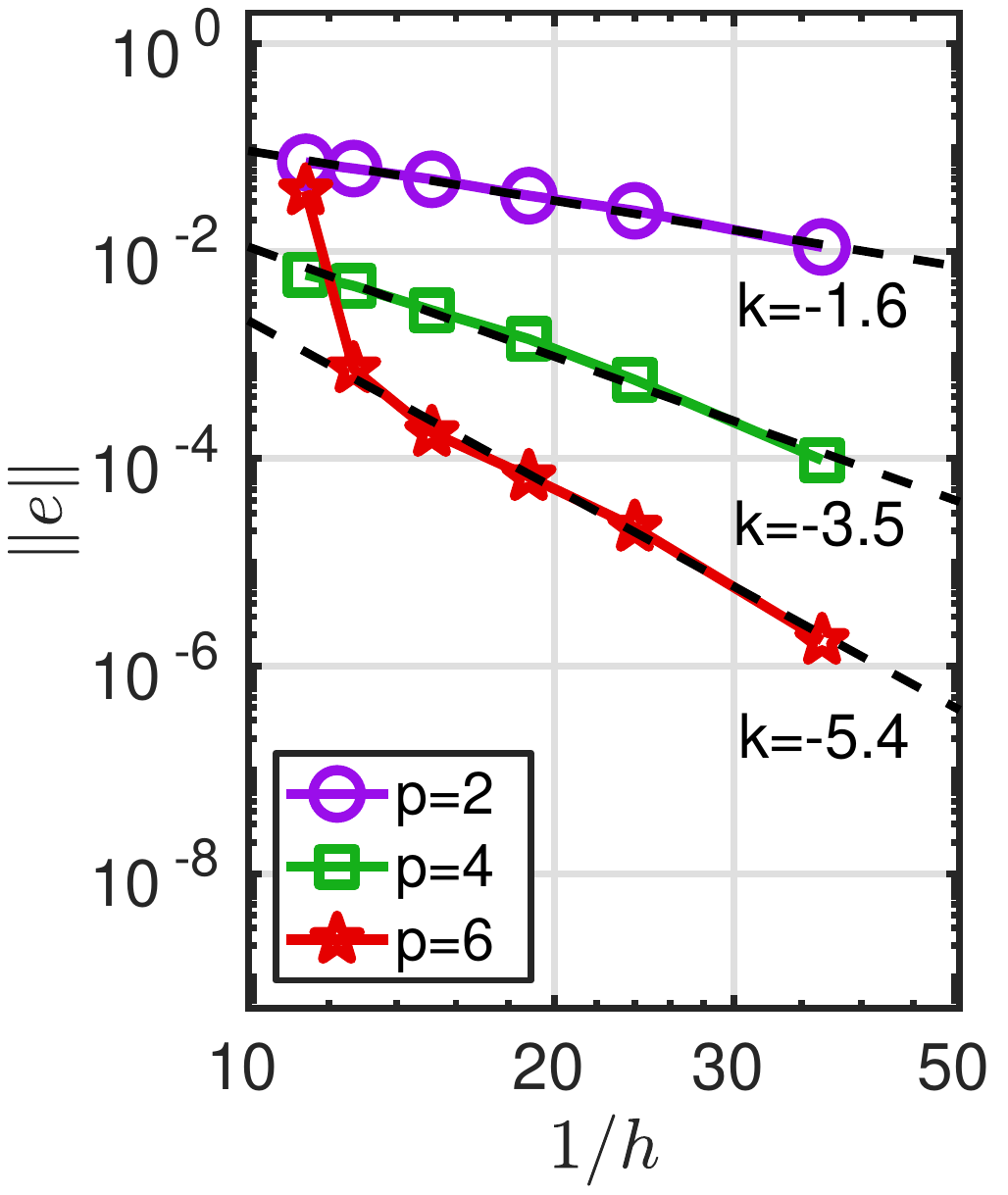}
    \includegraphics[width=0.37\linewidth]{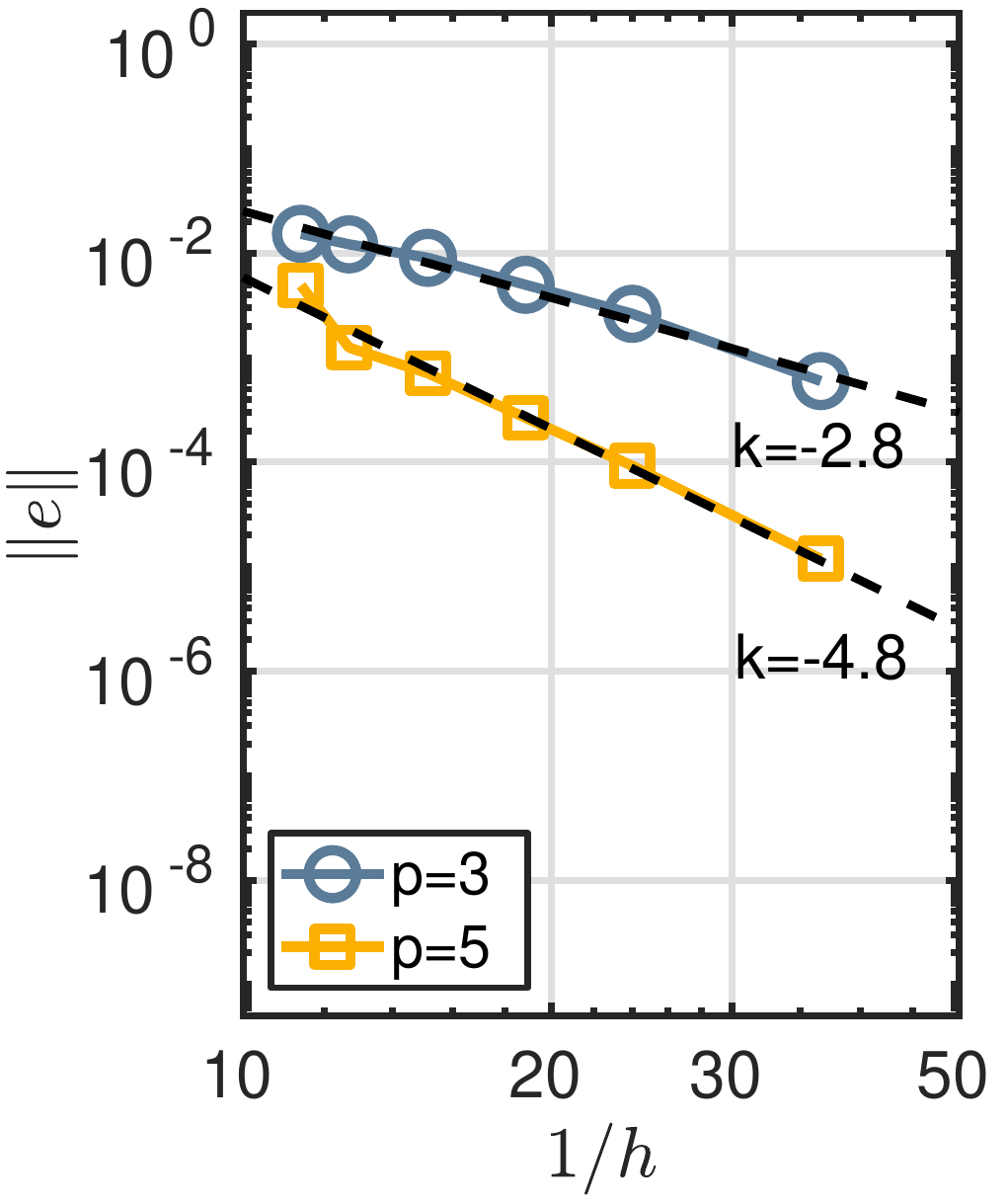}
    \caption{Error as a function of the inverse internodal distance $1/h$ when the oversampling parameter is set to $q=10$,
        for polynomial degrees $p=2$, $p=4$, $p=6$ (left image) and polynomial degrees $p=3$, $p=5$ (right image).}
    \label{fig:experiments:diaphragm:err_hrefinement}
\end{figure}

We test the convergence of the solution under node refinement for different polynomial degrees used to form the stencil
approximations.
The oversampling parameter is always fixed to $q=10$.
The spatial distribution of the magnitude of the relative error in logarithmic scale when $h=0.066$ and $p=5$ is given in Figure \ref{fig:experiments:diaphragm:absError}.
The convergence results can be observed in Figure \ref{fig:experiments:diaphragm:err_hrefinement}. We can see that the numerical solution converges
for all polynomial degrees $p$ with at least $\mathcal O(h^{p-1})$. This is an expected result according to the error
estimate \eqref{eq:theory:finalestimate} and also according to the numerical experiments previously made for the two-dimensional
cases in Section \ref{section:experiments_butterfly} and Section \ref{section:experiments_sprocket}.

\section{Final remarks}
\label{section:finalremarks}
In this paper we presented the unfitted version of the RBF-FD method in the least-squares setting and in this way 
simplified handling of complex 2D and 3D geometries.

We 
developed a criterion that numerically establishes the linear independence of the cardinal functions when the interpolation points 
are placed in the exterior of $\Omega$. 

Next, we numerically verified on a two-dimensional butterfly domain that the presented method 
is stable and that the 
stability properties closely follow the properties of the fitted RBF-FD-LS method. 
Moreover, the experiments confirmed that the error under node refinement decays as $\mathcal O(h^{p-1})$, which was also 
indicated in our theoretical work. The error was in the majority of the cases found to be smaller compared with the error 
RBF-FD-LS and RBF-FD-C.

Through an example of a drilled sprocket we numerically demonstrated that the 
less skewed boundary stencils allow 
the unfitted RBF-FD-LS method to outperform RBF-FD-LS and RBF-FD-C in the sense of the spatial error distribution, when the 
discretization is built upon large stencils.

Lastly, using a thoracic diaphragm of a human being as a computational domain, we demonstrated that 
the unfitted RBF-FD-LS method can also be used in a three dimensional setup, where the observed convergence trends follow 
$\mathcal O(h^{p-1})$, in the same way as in the two-dimensional cases.

Future work includes using the unfitted RBF-FD-LS method for other elliptic PDEs and investigating a formulation for time-dependent PDEs. 
\refereeThird{A potentially interesting topic is also to find a way to form the unfitted RBF-FD-C method, 
and then compare it to the unfitted RBF-FD-LS method.}

\section*{Acknowledgement}
We thank Alfa Heryudono from University of Massachusetts Dartmouth and Elisabeth Larsson from Uppsala University 
for fruitful discussions. Furthermore we thank Elisabeth Larsson for providing point sets 
over the surface of the diaphragm. 

The first author was supported by the Swedish Research Council,  grant  no.  2016-04849. 
The second author was supported by Center for Interdisciplinary Mathematics, Uppsala University.
\begin{appendices}

    \section{\refereeSecond{Skeweness of the stencils influences the size of the interpolation error bound}}
    \label{sec:appendix:skeweness}
    In this section we provide an insight into the interpolation error behavior over the highly skewed boundary stencils.
    
\begin{figure}[h!]
        \centering
        \includegraphics[width=0.45\linewidth]{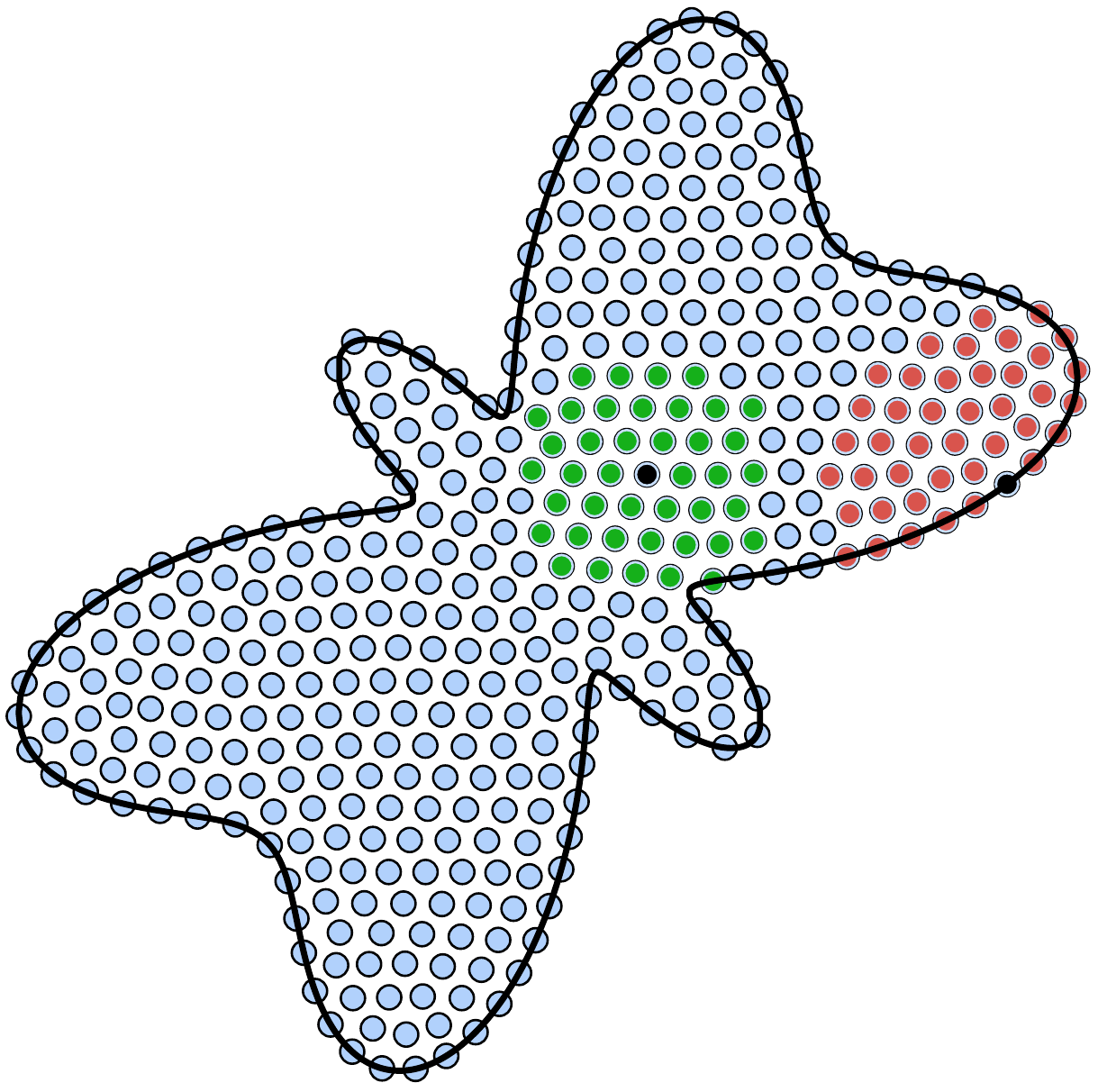}
        \includegraphics[width=0.45\linewidth]{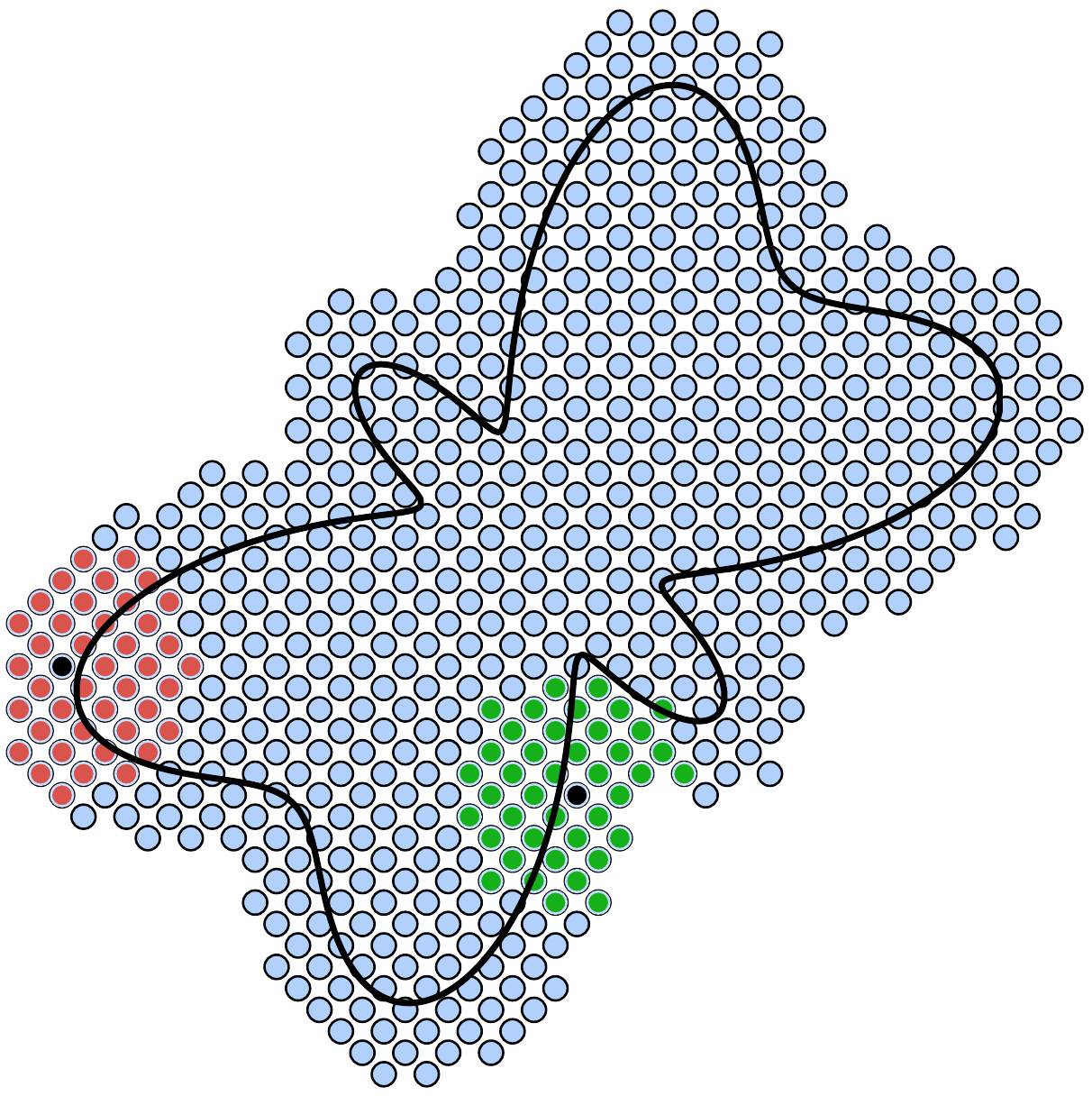}
        \caption{The left image displays two stencils over a set of interpolation points $X$ for the fitted RBF-FD methods, where the red stencil indicates that its interpolation matrix has the largest norm
            $\|\tilde A\|_\infty$ among all stencils over $X$. The green stencil is the stencil with the smallest norm $\|\tilde A\|_\infty$. The right image shows the red stencil and the green stencil chosen according
            to the same criterion, but for the unfitted RBF-FD-LS method.}
        \label{fig:experiments:butterfly:href:stencilMatrixNorms}
    \end{figure}

    We start by making an experiment. We spread interpolation points $X$ over a butterfly domain in the same way as in Section \ref{sec:experiments_butterfly:pointsets}, 
    and fix $p=5$, $h=0.05$. Then we measure the norm of the inverse interpolation matrix of every stencil and pick out the
    stencils with the smallest and the largest norms:
    $$\xi_{\min} = \min_i\, \|\tilde A_{(i)}^{-1}\|_\infty,\quad \xi_{\max} = \max_i\, \|\tilde A_{(i)}^{-1}\|_\infty,\quad i=1,..,N.$$
    We draw the shapes of these two stencils for the fitted RBF-FD-LS method and the unfitted RBF-FD-LS method in Figure \ref{fig:experiments:butterfly:href:stencilMatrixNorms}, where the stencil corresponding to
    $\xi_{\min}$ is colored with green color and the stencil corresponding to $\xi_{\max}$ with red color.
    The rounded values of norms for the fitted RBF-FD-LS method are $\xi_{\min} = 2.3 \cdot 10^2$ and $\xi_{\max} = 2.4 \cdot 10^{4}$. In this case $\xi_{max}$ is around two orders larger in magnitude and corresponds to the highly skewed stencil in
    Figure \ref{fig:experiments:butterfly:href:stencilMatrixNorms}.
    
    \begin{figure}[h!]
        \centering
        \includegraphics[width=0.42\linewidth]{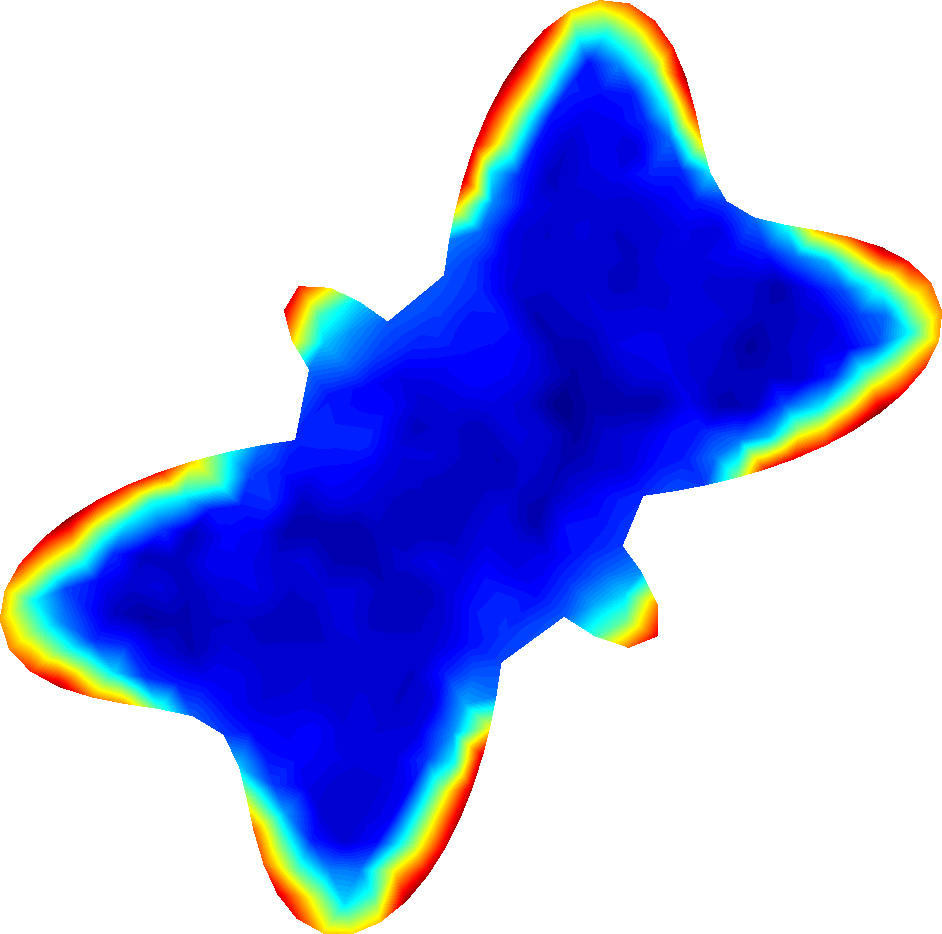}
        \includegraphics[width=0.55\linewidth]{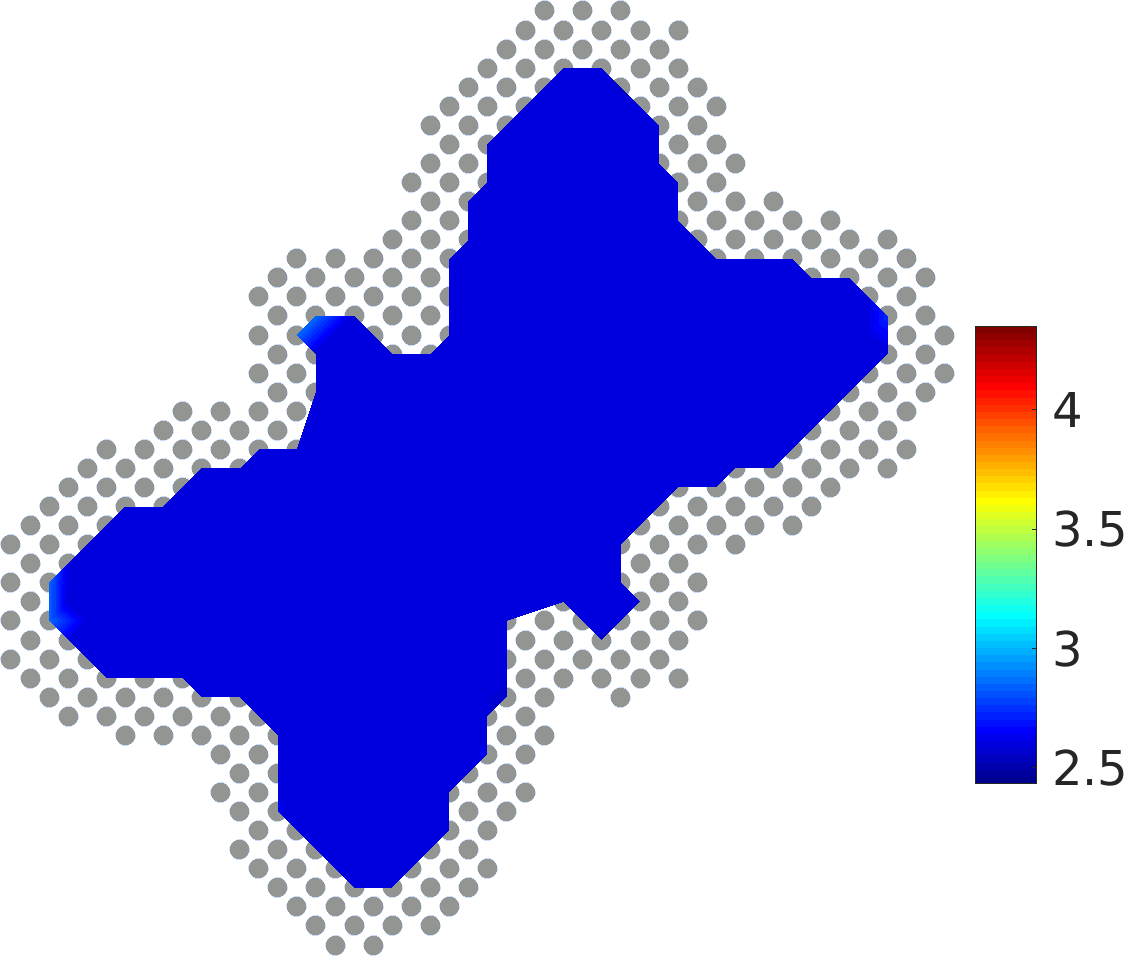}
        \caption{Both images show a spatial distribution of the matrix norm $\|\tilde A_{(i)}^{-1}\|_\infty$ intrinsic to the corresponding stencil center $x_i \in X$, in $\log_{10}$ scale.
        Left image is showing the distribution for the fitted RBF-FD methods, 
        while the right image is showing the distribution for the unfitted RBF-FD-LS method.}
        \label{fig:experiments:butterfly:href:stencilMatrixNorms_field}
    \end{figure}
    For the unfitted RBF-FD-LS method we have
    $\xi_{\min} = 3.7 \cdot 10^2$ and $\xi_{\max} = 7.2 \cdot 10^{2}$, where both norms are similar in magnitude, and both stencils are unskewed, despite both being positioned close to the boundary.
    Thus, we experimentally established that the skeweness of stencils is indeed reflected in the size of $\| \tilde A^{-1} \|_\infty$. Note that we measured $\xi_{\min}$ and $\xi_{\max}$ in $2$-norm and $1$-norm, 
    and the outcome did not change in the interpretation. Another perspective is given in Figure \ref{fig:experiments:butterfly:href:stencilMatrixNorms_field}, where we plot the $\log_{10}$ spatial distribution of the norm 
    $\|\tilde A_{(i)}^{-1}\|_\infty$ for every corresponding stencil center $x_i \in X$. In the fitted RBF-FD case 
    we can see that many boundary stencils have a much larger $\|\tilde A_{(i)}^{-1}\|_\infty$ compared to the interior stencils, while in the unfitted RBF-FD-LS case we can see that the 
    matrix norm distribution is much more evenly distributed.

    Now we focus on the interpolation over
    the $i$-th stencil \eqref{eq:method:stencil_cardinal}. We know that the interpolation error when using a polynomial interpolant in the Lagrange form \cite{Powell, Ibrahimoglu_lebesgue} leads to an error estimate:
    \begin{equation}
        \label{fig:experiments:butterfly:href:stencilMatrixNorms:errorbound}
        \|u_h^{(i)}(z) - u^{(i)}(z)\|_\infty \leq (1 + \Lambda^{(i)})\, \|u^{(i)}(z) - p_*^{(i)}(z)\|_\infty,
    \end{equation}
    where $u^{(i)}(z)$ is an exact function value over the $i$-th stencil, $p_*^{(i)}(z)$ is the best interpolating polynomial over the $i$-th stencil, and:
    \begin{equation}
        \Lambda^{(i)} = \max_z \sum_{k=1}^n |\psi_k^{(i)}(z)| = \max_z \|\underline \psi^{(i)}(z)\|_1,
        \label{fig:experiments:butterfly:href:stencilMatrixNorms:lebesgue}
    \end{equation}
    is a Lebesgue constant over the $i$-th stencil. Here $\underline \psi^{(i)}(z)$ is a vector of local cardinal functions over the $i$-th stencil. For a fixed $u^{(i)}(z)$, the term $p_*^{(i)}(z)$ is fixed and thus the upper limit of the error \eqref{fig:experiments:butterfly:href:stencilMatrixNorms:errorbound}
    depends only on $\Lambda^{(i)}$.
    In order to bound $\Lambda^{(i)}$ we plug the definition of $\underline \psi^{(i)} = \{\psi_k(z)\}_{k=1}^n$ from \eqref{eq:method:stencil_cardinal} into \eqref{fig:experiments:butterfly:href:stencilMatrixNorms:lebesgue}:
    \begin{eqnarray}
        \label{fig:experiments:butterfly:href:stencilMatrixNorms:lebesgueBound}
        \Lambda^{(i)} = \max_z\, \|\underline \psi^{(i)}(z) \|_1 = \max_z\, \| \underline w^{(i)}(z) \|_1 &=& \max_z\, \| b^{(i)}(z)\,  \tilde A_{(i)}^{-1} \|_1 \nonumber \\
        &\leq& \max_z\, \| b^{(i)}(z) \|_1\, \| \tilde A_{(i)}^{-1} \|_1 \nonumber \\
        &\leq& \sqrt{n}\, \max_z\, \| b^{(i)}(z) \|_1\, \| \tilde A_{(i)}^{-1} \|_\infty
    \end{eqnarray}
    From the bound and the experiment in which we found out that the worst $\| \tilde A^{-1} \|_\infty$ among all stencils belong to the skewed stencil,
    we have that for a fixed $p$, $n$, $h$, the Lebesgue constant has a larger upper limit for the skewed stencil compared to the less skewed stencil and thus, the interpolation error also has a larger upper limit 
    over the skewed stencil. 
    
    Finally, we have to understand how \eqref{fig:experiments:butterfly:href:stencilMatrixNorms:lebesgueBound} relates to the global solution defined in \eqref{eq:cardinalF}. 
    We restrict the global solution $u_h(y)$ to one Voronoi cell $\mathcal{V}_i$, where the cell is centered around a point $x_i \in X$. 
    This restriction is denoted by $u_h(y)|_{\mathcal{V}_i}$. In the fitted/unfitted RBF-FD-LS method every $u_h(y)|_{\mathcal{V}_i}$ is approximated using an interpolant over the $i$-th stencil, where $x_i$ has a role of the stencil center point. 
    The interpolation error estimate \eqref{fig:experiments:butterfly:href:stencilMatrixNorms:errorbound} for the $i$-th stencil then applies to every interpolant $u_h(y)|_{\mathcal{V}_i}$: from an interpolation perspective, $u_h(y)|_{\mathcal{V}_i}$ 
    will contain an error dominated by the Lebesgue constant over the $i$-th stencil. It is therefore beneficial to keep $\| \tilde A_{(i)}^{-1} \|_\infty$ as small as possible, i.e., to keep the stencils skewed as little as possible.
    \end{appendices}

\bibliography{references}

\end{document}